\documentclass[preprint,sort&compress,12pt]{elsarticle}
\usepackage{bbm}
\usepackage{amstext}

\usepackage{amsmath}
\usepackage{amssymb}
\usepackage{mathrsfs}
\usepackage{bm}
\usepackage{pstricks}
\usepackage{pst-coil}
\usepackage{pst-3d}
\usepackage{amsfonts}
\usepackage{amsthm}
\usepackage{xcolor}
\allowdisplaybreaks
\newcommand{\mm}{\mathrm}

\newcommand{\be}{\begin{equation}}
\newcommand{\bea}{\begin{equation}\begin{aligned}}
\newcommand{\beas}{\begin{equation*}\begin{aligned}}
\newcommand{\eeas}{\end{aligned}\end{equation*}}
\newcommand{\eea}{\end{aligned}\end{equation}}
\newcommand{\ee}{\end{equation}}

\renewcommand{\div}{{\rm div }}

\begin{document}
\begin{frontmatter}
\title{On Stability and Instability of
Gravity Driven Navier--Stokes--Korteweg Model in Two Dimensions}
\author[FJ,1FJ,2FJ]{Fei Jiang}
\ead{jiangfei0951@163.com}
 \author[SJ]{Fucai Li}
 \ead{fli@nju.edu.cn}
\author[SJ]{Zhipeng Zhang}\ead{zhangzhipeng@nju.edu.cn}
\address[FJ]{School of Mathematics and Statistics, Fuzhou University, Fuzhou, 350108, China.}
\address[SJ]{Department of Mathematics, Nanjing University, Nanjing 210093, China}
 \address[1FJ]{Center for Applied Mathematics of Fujian Province, Fuzhou 350108, China.}
\address[2FJ]{Key Laboratory of Operations Research and Control of Universities in Fujian, Fuzhou 350108, China.}
\begin{abstract}
 Bresch--Desjardins--Gisclon--Sart have derived  that the capillarity can slow  down the growth rate of  Rayleigh--Taylor (RT) instability in the capillary fluids
based on the \emph{linearized} two-dimensional (2D) Navier--Stokes--Korteweg equations  in 2008. Motivated by their linear theory, we further investigate the \emph{nonlinear}  RT problem for  the 2D  incompressible case  in a horizontally periodic   slab domain with  Navier boundary condition, and rigorously verify  that the RT instability can be inhibited by capillarity under our 2D setting. More precisely, if the RT density profile $\bar{\rho}$ satisfies an additional stabilizing condition, then  there is a threshold  $\kappa_{\mm{C}}$ of capillarity coefficient, such that if the  capillarity coefficient $ \kappa $ is bigger than $\kappa_{\mm{C}}$, then the small perturbation solution around the  RT equilibrium state is \emph{algebraically} stable in time.  In particular, if the RT density profile is linear, then the threshold $\kappa_{\mm{C}}$ can be given by the formula $\kappa_{\mm{C}}=g /(\pi^2h^{-2}+L^{-2})\bar{\rho}'$, where $2\pi L$ denotes the  length of a periodic cell of the slab domain in the horizontal direction, and $h$ the height of the slab domain.   In addition, we also provide a nonlinear instability result for $\kappa\in[0, \kappa_{\mm{C}})$. The instability result presents that the capillarity can not inhibit the RT instability, if its strength is too small.
\end{abstract}
\end{frontmatter}
\newtheorem{thm}{Theorem}[section]
\newtheorem{lem}{Lemma}[section]
\newtheorem{pro}{Proposition}[section]
\newtheorem{concl}{Conclusion}[section]
\newtheorem{cor}{Corollary}[section]
\newproof{pf}{Proof}
\newdefinition{rem}{Remark}[section]
\newtheorem{definition}{Definition}[section]

\section{Introduction}\label{introud}

The equilibrium of a heavier fluid on the top of a lighter one, subject to gravity, is usually unstable. In fact, small disturbances acting on the equilibrium will grow and lead to the release of potential energy,
 as the heavier fluid moves down under gravity, and the lighter one is displaced upwards. This phenomenon was first studied
 by Rayleigh \cite{RLIS} and then Taylor \cite{TGTP}, and is called therefore the Rayleigh--Taylor (RT) instability.
In the last decades, the RT instability  {has been} extensively investigated from
 {physical,   numerical, and  mathematical aspects, among others, we mention  \cite{CSHHSCPO,WJH,desjardins2006nonlinear,hateau2005numerical,GBKJSAN} for examples.} It has been also widely analyzed {on how the} physical factors,
such as elasticity \cite{JFJWGCOSdd}, rotation \cite{CSHHSCPO,BKASMMHRJA},
internal surface tension \cite{GYTI2,WYJTIKCT,JJTIWYJ2016ARMA},
magnetic  {fields} \cite{JFJSWWWOA,JFJSJMFMOSERT,WYJ2019ARMA,JFJSARMA2019},  {capillarity \cite{bresch2008instability}} and so on, influence the dynamics of   RT instability.

In this paper, we are interested in the influence of   capillarity on the evolution of  RT instability.
In \cite{bresch2008instability},  Bresch--Desjardins--Gisclon--Sart derived  that the capillarity can slow down the growth rate of  RT instability in the  fluid  endowed with internal capillarity (in the diffuse interface setting)  based on the \emph{linearized}  two-dimensional (2D) Navier--Stokes--Korteweg (NSK) equations.
Motivated by  {their} linear theory, we further investigate the \emph{nonlinear} RT  problem for  the  {2D} incompressible inhomogeneous case  in a horizontal slab domain with {Navier boundary condition}, and rigorously  {prove} a phenomenon that the   RT instability can be  {also} inhibited by capillarity under our setting. Before stating our result  in details, {we {need} to  mathematically formulate the phenomenon of the inhibition of RT instability by capillarity}.

\subsection{Mathematical formulation for the capillary RT problem}\label{subsec:01}
\numberwithin{equation}{section}

A classical model to  describe the dynamics of a fluid endowed with internal capillarity in the presence of a uniform gravitational
field is the following  {general system of} compressible NSK equations:
\begin{align}\label{0101ooo}
\begin{cases}
 \rho_t  +    \mathrm{div}(\rho v) = 0, \\
 \rho v_t + \rho v\cdot \nabla v -\mathcal{A}v + \nabla \mathcal{P}(\rho ) = \mathrm{div}  {K} -\rho g\mathbf{e}^n,
\end{cases}
\end{align}
where ${\rho (x,t)} \in \mathbb{R}^{+}$ and $v(x,t) \in {\mathbb{R}^n}$  denote the density and velocity of the fluid {respectively} at the position $x \in {\mathbb{R}^n}$ with the spacial dimension $n \geqslant 2$ for time $t\in[0, {+\infty} )$. The diffusion operator is $\mathcal{A}v = \mathrm{div}(\mu (\rho )\mathbb{D}v) + \nabla (\lambda (\rho ) {\mathrm{div}\,v})$, and the differential operator $\mathbb{D}$ is defined by $\mathbb{D}v = \nabla v+\nabla v^{\top},$ where the subscript $\top$ denotes the transposition.   $\mathbf{e}^n $ represents the unit vector with the $n$-th component being $1$,  {$g> 0$} the gravitational constant  and $-\rho g \mathbf{e}^n$ the gravity.
The shear viscosity function  $\mu$, the bulk  viscosity function $\lambda$ and the capillarity function $\kappa$ are known smooth functions $\mathbb{R}^+ \to \mathbb{R}$, and satisfy $\mu>0$,  { $2\mu+n \lambda \geqslant 0$} and $\kappa>0$.
The pressure $\mathcal{P}$ is a given function (usually chosen as the Van der Waals pressure), and the general capillary tensor is written as
\begin{align}
 K =  \left( {\rho \mm{div}(\kappa(\rho )\nabla \rho ) + \left( {\kappa(\rho ) - \rho\kappa'(\rho )} \right){{\left| {\nabla \rho } \right|}^2}} /2\right)\mathbb{I} -  {\kappa(\rho )\nabla \rho  \otimes \nabla \rho },
\end{align}
where $\mathbb{I}$ denotes the identity matrix.

In classical hydrodynamics, the interface between two immiscible fluids is modeled as a free boundary  {which} evolves in time. The equations describing
the motion of each fluid are supplemented by boundary conditions
at the free surface  involving the physical properties of the interface.
For instance, in the free-boundary formulation, it is assumed that the interface
has an internal surface tension.
However, when the interfacial thickness is
comparable to the length scale of the phenomena being examined, the free-boundary description breaks down.
Diffuse-interface models provide an alternative description where the surface tension is expressed in its simplest form as $\mm{div}K$,
i.e., the capillary tension which was introduced by  Korteweg in 1901 \cite{korteweg1901forme}. Later, its modern form was derived by Dunn and Serrin \cite{dunn1985thermomechanics}.

 In the physical view, it can serve as a phase transition model to describe the motion of compressible fluid with capillarity effect.
 Owing to its importance in mathematics and physics, there has been profuse  {works} on the mathematical theory of
 {compressible NSK  system, for example,} see  \cite{antonelli2022global,bresch2003some,haspot2011existence} for the global-in-time weak solutions with large initial data, \cite{danchin2001existence,kobayashi2022resolvent,hong2020strong,murata2020global,haspot2020strong}
 for the global-in-time strong solutions with small initial data, \cite{kotschote2010strong,charve2014local,hattori1994solutions,kotschote2008strong}
 for local-in-time strong/classical solutions with large initial data,   \cite{chen2014existence} for stationary solutions, \cite{fanelli2016highly} for highly rotating limit, \cite{hong2022stability,chen2021asymptotic} for the stability  of viscous shock wave,
\cite{bian2014vanishing,germain2016finite,jungel2014asymptotic,charve2013existence}
for the vanishing capillarity limit,  \cite{saito2020maximal} for the  maximal $L^p$--$L^q$ regularity theory,  \cite{kawashima2021lp,charve2018gevrey,bresch2022exponential} for the decay-in-time of global solutions and so on. In addition, the inviscid case is also widely investigated, see \cite{audiard2017global,bresch2019navier,sy2006local,benzoni2007well} and the references cited therein.

To conveniently investigate the influence of   capillarity on the  {flows}, we consider \emph{$\kappa$ to be a positive constant} as in \cite{bresch2008instability}, and thus get
\begin{align}
\mathrm{div}K =\kappa \rho\nabla \Delta \rho.
\end{align}
 Let us further choose an equilibrium $(\bar{\rho},0)$ to \eqref{0101ooo}, where the density profile $\bar{\rho}$ only depends on the $n$-th variable. Then we have the equilibrium state equation
\begin{align}
   \nabla  \mathcal{P}(\bar{\rho} )  = \kappa\bar{\rho}\nabla \Delta \bar{\rho} - g\bar{\rho}\mathbf{e}^n ,  \nonumber
\end{align}
which can be rewritten as an ordinary differential equation on $\bar{\rho}$:
$$ {\frac{\mm{d}\mathcal{P}(\bar{\rho} )}{\mm{d}x_n}=
 \kappa\bar{\rho}\frac{\mm{d}^3\bar{\rho}}{\mm{d}x_n^3} - g\bar{\rho}}.$$

As Bresch--Desjardins--Gisclon--Sart pointed out, how to analyse the stability/instability of the above equilibrium $(\bar{\rho},0)$ is still an interesting (open) problem  \cite{bresch2008instability}.
Due to the essential difficulty arising from the compressibility and the dimensions (in a slab domain),  we  focus on  the 2D incompressible counterpart as in \cite{bresch2008instability}, namely the following 2D incompressible (inhomogeneous) NSK system:
\begin{equation}
\label{0101}
\begin{cases}
\rho_t+ v\cdot \nabla \rho=0,\\
\rho {v}_t+\rho{v}\cdot\nabla {v}+\nabla P-\mm{div}(\mu\nabla {v})=\kappa\rho\nabla \Delta \rho-\rho g\mathbf{e}^2,\\
\div  v = 0,
\end{cases}
\end{equation}
where $P:=P(x,t)$ denotes the kinetic pressure of the fluid. We mention that the well-posdeness problem for the incompressible NSK system has been also widely investigated, see \cite{burtea2017lagrangian} and the references cited therein.

We consider the horizontally periodic motion solutions of \eqref{0101}, and thus define a horizontally periodic domain by
\begin{align}\label{0101a}
\Omega:= 2\pi L\mathbb{T} \times(0,h),
\end{align}
where $\mathbb{T}=\mathbb{R}/\mathbb{Z}$ and $L>0$. \emph{We will see that the stabilizing effect of capillarity depends on
the periodic length $2\pi L$}. For the horizontally periodic domain $\Omega$, the 1D periodic domain
$2\pi L\mathbb{T}\times \{0,h\}$, denoted by $\partial\Omega$,  is customarily regarded as the  boundary of $\Omega$.
For the well-posedness of the system \eqref{0101}, we impose the following initial and boundary conditions:
\begin{gather}
  (\rho,v)|_{t=0}=(\rho^0,v^0),\label{20210031303} \\
 \label{0101b}
v|_{\partial\Omega}\cdot{\mathbf{n}}=0,  \ ((\mathbb{D}v|_{\partial\Omega}){\mathbf{n}})_{\mm{tan}}=0,
\end{gather}
where  {${\mathbf{n}}$ } denotes the outward unit normal vector to
$\partial\Omega$, and the subscript ``$\mm{tan}$" means the tangential component of a vector  {(for example,  $v_\mm{tan}=v-(v\cdot {\mathbf{n}}){\mathbf{n}}$).}
  Here and in what follows, we always use the superscript $0$ to emphasize the initial data.
 We call
\eqref{0101b} the Navier (slip) boundary condition.
 Since $\Omega$ is a slab domain, the Navier boundary condition is equivalent to the boundary condition
\begin{align}
(v_2,\partial_2 v_1)|_{\partial\Omega}=0 .
\label{20220202081737}
\end{align}
We will {explain} the reason for choosing the above Navier boundary condition in Section \ref{subsec:02}.

To investigate the influence of the capillarity on  RT instability,  we shall choose an RT density profile $\bar{\rho}:=\bar{\rho}(x_2)$, which is independent of $x_1$ and satisfies
\begin{gather}
 \label{0102}
\bar{\rho}\in {C^5}(\overline{\Omega}),\quad \inf_{ x\in {\Omega}}\bar{\rho}>0,\\[1mm]
 \label{0102n}\bar{\rho}'|_{x_2=s}>0  \ \   \mbox{ for some} \ \ s\in \{x_2~|~(x_1,x_2)^{{\top}}\in \Omega\},
\end{gather}
where $\bar{\rho}':=\mm{d}\bar{\rho}/\mm{d}x_2$ and $\overline{\Omega}:= {\mathbb{R}}\times [0,h]$.
We remark that the second condition in \eqref{0102}
prevents us from treating vacuum, while the condition in \eqref{0102n} is called the RT condition,
which assures that there is at least a region in which the density is larger with increasing height $x_2$, thus leading to the classical RT instability, see \cite[Theorem 1.2]{JFJSO2014}. \emph{However, we will see that such instability can be inhibited by capillarity}.

With the RT density profile in hand, we further define an RT equilibrium $ (\bar{\rho}, 0 )$ for the system \eqref{0101}.
The pressure profile $\bar{P}$ under the equilibrium state is determined by the hydrostatic relation
\begin{equation}
\nabla \bar{P}=\kappa\bar{\rho}\nabla \Delta \bar{\rho}-\bar{\rho}g  \mathbf{e}^2 \  \ \mbox{ in }\ \  {\Omega},  \label{equcomre}
\end{equation}
i.e.,
$$  \bar{P}' =\kappa\bar{\rho}\bar{\rho}''' -\bar{\rho}g.$$

From now on, we only consider $\mu$ in \eqref{0101} to be a constant  for simplicity.
Denoting the perturbation around the  RT equilibrium by
$$\varrho=\rho -\bar{\rho}, \ v= v- {0}, $$
 and recalling the relations \eqref{equcomre} and
$$
 f\nabla \Delta f= \nabla (f\Delta f+ |\nabla f|^2/2) - \mm{div}(\nabla f\otimes \nabla f) ,
 $$
  we obtain  {the system of  perturbation equations} from \eqref{0101}:
\begin{equation}\label{0103} \begin{cases}
\varrho_t+{  v}\cdot\nabla (\varrho+\bar{\rho})=0, \\
(\varrho+\bar{\rho}){  v}_t+(\varrho+\bar{\rho}){  v}\cdot\nabla
{ v}+\nabla  \beta-\mu \Delta v\\
\qquad =  \kappa\mm{div}(\nabla\bar{\rho} \otimes \nabla \bar{\rho}- \nabla \rho \otimes \nabla \rho)
- \varrho  g \mathbf{e}^2,\\
 \mm{div}v =0,\end{cases}\end{equation}
where $\beta:= P- \bar{P}+\kappa(\bar{\rho}\Delta \bar{\rho}
- \rho\Delta \rho+(|\nabla \bar{\rho}|^2-|\nabla \rho|^2) /2 ) $.
The corresponding initial and boundary conditions read as follows:
\begin{gather} \label{c0104}
(\varrho,v  )|_{t=0}=(\varrho^0,v^0 ) ,  \\
 \label{0105}
 {(v_2,\partial_2 v_1)|_{\partial\Omega}=0}.
\end{gather}
We call the initial-boundary value problem \eqref{0103}--\eqref{0105} {\it the  capillary RT (abbr. CRT)   problem}  for simplicity.
Obviously, to prove the inhibition of the RT instability by capillarity  in a 2D slab domain,
it suffices to verify the stability in time of solutions to the CRT problem
 with some non-trivial initial data.

Due to the nonlinearity of capillary tensor, it seems to be difficult to directly prove the stability of the CRT problem.  In  \cite{LFCZZP}, Li--Zhang established a linear stability result for the CRT problem with a non-slip boundary condition (in place of the Navier boundary condition). Following  Li--Zhang's argument, we also easily obtain a linear stability result for  our  CRT problem, i.e.  there exists  a threshold $\kappa_{\mm{C}}$ such that the  linearized CRT problem (i.e. omitting the nonlinear perturbation terms in \eqref{0103}) is stable under a sharp stability condition
 \begin{align}
\label{2020102241504} \kappa>\kappa_{\mm{C}},
\end{align} where we have defined that
 \begin{align}
\label{2saf01504}
 \kappa_{\mm{C}}:= {\sup_{ { w}\in H_{\sigma}^1}\frac{g\int\bar{\rho}' { w}_2^2 {\mm{d}x}}
{ \int |\bar{\rho}'\nabla  w_2|^2 {\mm{d}x}}}
 \end{align}
 (It should be noted that we have excluded the function $w$ satisfying $ \int |\bar{\rho}'\nabla  w_2|^2 {\mm{d}x}=0$ in the above  definition of the supremum by default)
 and $\bar{\rho}$ should satisfy the stabilizing condition
\begin{align}
\label{2022205071434}
\inf_{ {x_2}\in (0,h)} \{|\bar{\rho}' {(x_2)}|\}>0.
\end{align}
Here and in what follows, $\int:= \int_{(0,2\pi L)\times (0,h)}$ and $H^1_{\sigma}:=\{w\in  {W^{1,2}(\Omega)}~|~\mm{div}w=0,\ w_2|_{\partial\Omega}=0\}$.
\emph{The aim of this paper is to develop new ideas  to further establish the nonlinear stability for the CRT problem under the stability conditions \eqref{2020102241504} and \eqref{2022205071434}.}

\subsection{Reformulation in Lagrangian coordinates}\label{subsec:02}

In \cite{WYJTIKCT}, Wang--Tice--Kim proved that the internal surface tension can inhibit RT instability. Later similar inhibition phenomena  {caused} by magnetic tension  {has been also}  verified in the non-resistive magnetohydrodynamic (MHD) fluids \cite{WYJ2019ARMA,JFJSJMFMOSERT} in Lagrangian coordinates, see also  \cite{JFJWGCOSdd} for  elasticity in the viscoelastic fluids.  As Lin pointed out in \cite{lin2012some}, the stratified fluids with the internal surface tension, viscoelastic fluids, non-resistive MHD fluids and the diffuse-interface model (or NSK model) can also  be regarded  {as} complex fluids with elasticity.  Motivated by the stability results in  \cite{WYJ2019ARMA,JFJSJMFMOSERT} in Lagrangian coordinates,
we also start with reformulation of the CRT problem in Lagrangian coordinates. We mention that Burtea--Charve
\cite{burtea2017lagrangian}  ever studied the existence of global-in-time strong solutions with small initial data for the incompressible  NSK system in Lagrangian coordinates.

Let the flow map $\zeta$ be the solution to the initial value problem:
\begin{equation}
\label{201806122101}
            \begin{cases}
\partial_t \zeta(y,t)=v(\zeta(y,t),t),
\\[1mm]
\zeta(y,0)=\zeta^0(y),
                  \end{cases}
\end{equation}
where the invertible mapping $\zeta^0:=\zeta^0(y)$ maps $\Omega$ to $\Omega$, and satisfies
\begin{gather}
 \label{zeta0inta}
J^0:=\det \nabla \zeta^0=1 ,\\[1mm]
 \label{zeta0inta0}
\zeta^0_2 =y_2\mbox{ on } \partial\Omega.
\end{gather}
Here and in what follows, ``$\det$'' denotes the determinant of a matrix.

We shall see later that the flow map $\zeta$ satisfies that for each fixed $t>0$,
\begin{gather}
  \zeta|_{y_2=r}:{\mathbb{R}\to \mathbb{R}}\,\ \ \mbox{ is a } C^1(\mathbb{R})\mbox{-diffeomorphism mapping for }r=0,\ h; \label{20210301715x}\\
 \zeta:\overline{\Omega}\to \overline{\Omega}\, \ \ \mbox{ is a } C^1(\overline{\Omega})\mbox{-diffeomorphism mapping}.\label{20210301715}
\end{gather}
Since $v$ is divergence-free and satisfies $v_2|_{\partial\Omega}=0$, we can deduce from \eqref{201806122101}--\eqref{zeta0inta0} that
\begin{gather}
 \nonumber
J:=\det\nabla \zeta=1 ,\\
 \nonumber
\zeta_2=y_2\ \  \mbox{ on } \ \ \partial\Omega.
\end{gather}

We define the matrix $\mathcal{A}:=(\mathcal{A}_{ij})_{2\times 2}$ by
\begin{align}\nonumber
\mathcal{A}^{\top}=(\nabla\zeta)^{-1}:= (\partial_j \zeta_i)^{-1}_{2\times 2}.
\end{align}
Then, we further define the differential operators $\nabla_{\mathcal{A}}$, $\mm{div}_{\mathcal{A}}$ and $\mm{curl}_{\mathcal{A}}$ as follows:
for a scalar function $f$ and a vector function $X:=(X_1,X_2)^{\top}$,
\begin{align}
&\nabla_{\mathcal{A}}f:=(\mathcal{A}_{1k}\partial_kf,
\mathcal{A}_{2k}\partial_kf)^{\top},\  \ \mm{div}_{\mathcal{A}}(X_1,X_2)^{\top}:=\mathcal{A}_{lk}\partial_k X_l, \ \ \Delta_{\mathcal{A}}f:=\mm{div}_{\mathcal{A}}\nabla_{\mathcal{A}}f \nonumber
\end{align}
and
\begin{align}
\mm{curl}_{\mathcal{A}}X:=\mathcal{A}_{1k}\partial_{k}X_2-\mathcal{A}_{2k}\partial_{k}X_1,
 \nonumber \end{align}
where we have used  the Einstein summation convention  over repeated indices, and $\partial_k:=\partial_{y_k}$. In particular, $\mm{curl}X:=\mm{curl}_{\mathbb{I}}X$. %

If we further define
\begin{align}
& \eta:=\zeta-y, \nonumber
\end{align}
then $\mathcal{A}^{\top}=(\nabla \eta+\mathbb{I})^{-1}$, in particular,
\begin{align*}
\mathcal{A}=
\begin{pmatrix}
1+\partial_2\eta_2 &-\partial_1\eta_2\\[1mm]
-\partial_2\eta_1 &1+\partial_1\eta_1
        \end{pmatrix}
\end{align*}
and
\begin{align}
\tilde{\mathcal{A}} :=\mathcal{A}-\mathbb{I}=
\begin{pmatrix}
\partial_2\eta_2 \quad-\partial_1\eta_2\\[1mm]
-\partial_2\eta_1\quad\partial_1\eta_1
             \end{pmatrix}. \nonumber
\end{align}
By the above expression of $\tilde{\mathcal{A}}$, it is easy to check that
 \begin{align}
 \label{202204120943}
 {\mm{curl}}_{\tilde{\mathcal{A}}} \eta  =0.
 \end{align}

Defining the Lagrangian unknowns:
\begin{equation*}
(  \vartheta, u , {Q} )(y,t):=(\rho,v,P )(\zeta(y,t),t)\; \mbox{ for } (y,t)\in \Omega \times\mathbb{R}^+_0,
\end{equation*}
then in Lagrangian coordinates, the initial-boundary value problem \eqref{0101}, \eqref{20210031303} and \eqref{20220202081737} can be
rewritten as
\begin{equation}\label{01dsaf16asdfasf00}
\begin{cases}
\zeta_t=u ,  \ \vartheta_t=0 ,\
\div_\mathcal{A}u=0 ,    \\[1mm]
\vartheta u_t+\nabla_{\mathcal{A}}Q
-\kappa \vartheta  \nabla_{\mathcal{A}} \Delta_{\mathcal{A}}\vartheta
+\vartheta g \mathbf{e}^2 = \mu \Delta_{\mathcal{A}}  u, \\[1mm]
(\zeta,\vartheta,u )|_{t=0}=(\zeta^0,\vartheta^0,u^0 ) , \\[1mm]
(\zeta_2-y_2,u_2, \mathcal{A}_{2i}\partial_iu_1)|_{ \partial\Omega} =0  ,
\end{cases}
\end{equation}
where $(\vartheta^0,  u^0 ):=
 (\rho^0(\zeta^0),v^0(\zeta^0) ) $. In addition, the relation \eqref{equcomre} in Lagrangian coordinates takes the form
 \begin{equation}
\label{dstist01}
\nabla_{\mathcal{A}}\bar{P}(\zeta_2) =  \kappa \bar{\rho}(\zeta_2) \nabla_{\mathcal{A}} \Delta_{\mathcal{A}} \bar{\rho}(\zeta_2)    -\bar{\rho}(\zeta_2)g\mathbf{e}^2 .
\end{equation}

Let
\begin{align}
q&= {P}(\zeta)-\bar{P}(\zeta_2)+\kappa(\bar{\rho}(\zeta_2)\Delta_{\mathcal{A}} \bar{\rho}(\zeta_2)
-\bar{\rho}(y_2)\Delta_{\mathcal{A}}\bar{\rho}(y_2)+(|\nabla_{\mathcal{A}} \bar{\rho}(\zeta_2)|^2-|\nabla_{\mathcal{A}} \bar{\rho}(y_2)|^2) /2 ) ,\nonumber \\
\label{202009130836}
 K^{\eta}& = \mm{div}_{\mathcal{A}} (\nabla_{\mathcal{A}}\bar{\rho}(\zeta_2)\otimes  \nabla_{\mathcal{A}}\bar{\rho}(\zeta_2)  - \nabla_{\mathcal{A}}\bar{\rho}(y_2)\otimes \nabla_{\mathcal{A}}\bar{\rho}(y_2)),\\
{G}^{\eta}& = g{(\bar{\rho}(\zeta_2)-\bar{\rho}(y_2)) }\nonumber  .
\end{align}
Now we   additionally assume that
 $\vartheta^0$  satisfies
 \begin{align}
\label{js1}
\vartheta^0  =\bar{\rho}(y_2),
\end{align}
then
\begin{align}
\label{202012280945}
 {\vartheta=\bar{\rho}(y_2)},
\end{align}
and  furthermore the initial-boundary value problem \eqref{01dsaf16asdfasf00}, together with the relations \eqref{dstist01} and {\eqref{202012280945}}, implies that
 \begin{align}\label{01dsaf16asdfasf}
&\begin{cases}
\eta_t=u ,\\[1mm]
\bar{\rho}u_t+\nabla_{\mathcal{A}} q-\mu \Delta_{\mathcal{\mathcal{A}}} u=\kappa\mathcal{K}^{\eta}+ {G}^{\eta}\mathbf{e}^2  ,\\[1mm]
\div_{\mathcal{A}} u=0  , \\[1mm]
(\eta,u)|_{t=0}=(\eta^0,u^0)  ,
\end{cases} \\
&(\eta_2, u_2) |_{\partial\Omega}=0, \label{20202201182345} \\
&( (1+\partial_1\eta_1  )\partial_2u_1 - \partial_2\eta_1 \partial_1u_1  )|_{\partial\Omega}=0.  \label{202345}
\end{align}
We call  the system of equations \eqref{01dsaf16asdfasf}$_1$--\eqref{01dsaf16asdfasf}$_3$ the transformed CRT (abbr. TCRT) system.

It is very interesting  {to point out that  the TCRT system is   similar} to the transformed magnetic RT  (abbr. TMRT)  system, which reads as  {follows} (see \cite{JFJSZYYO} for the relevant physical background)
\begin{align}\label{01dsafdasf16asdfasf}
                \begin{cases}
\eta_t=u ,\\[1mm]
\bar{\rho}u_t+\nabla_{\mathcal{A}} q-\mu \Delta_{\mathcal{\mathcal{A}}} u=\lambda m^2\partial_1^2\eta+ {G}^{\eta}\mathbf{e}^2  ,\\[1mm]
\div_{\mathcal{A}} u=0  .
\end{cases}
\end{align}
Since it is difficult to get the estimates of higher-order normal derivatives of $\eta$
due to the absence the boundary condition of $\eta_1$, Jiang--Jiang--Zhao naturally   imposed  the following boundary condition {on} the linearized {TMRT} problem \cite{JFJSZYYO}:
\begin{align}
& (\eta_2,\partial_2\eta_1, u_2,\partial_2u_1) |_{\partial\Omega}=0,  \label{20safd45}
\end{align}
which automatically implies the boundary conditions \eqref{20202201182345}  and \eqref{202345}. \emph{It is interesting that the boundary conditions \eqref{20202201182345} and \eqref{202345} can also yield  the  one \eqref{20safd45} if the solution $(\eta,u)$ satisfies $\partial_2\eta^0_1|_{\partial\Omega}=0$ and  additional smallness conditions, see Proposition A.1 in \cite{JFJSZYYO}.}
Consequently, they proved that the magnetic tension $\lambda m^2\partial_1^2\eta$ caused by the horizontal magnetic  field can inhibit the RT instablity based on the TMRT system \eqref{01dsafdasf16asdfasf} with the boundary condition \eqref{20safd45}.
 {For  presentation simplicity, we shall call} the initial-boundary value problem of \eqref{01dsaf16asdfasf}$_4$, \eqref{01dsafdasf16asdfasf} and \eqref{20safd45} the TMRT problem.

We mention that the boundary condition
\begin{align}
(\eta_2,\partial_2\eta_1) |_{\partial\Omega}=0
\label{202220118000}
 \end{align}
is called   the characteristic boundary condition. Indeed, if the initial data $\eta^0$ satisfies $(\eta_2^0,\partial_2\eta_1^0) |_{\partial\Omega} \linebreak =0$, then $\eta$ automatically satisfies \eqref{202220118000}
due to \eqref{01dsaf16asdfasf}$_1$ and the boundary condition
\begin{align}
\label{202201180853}
(u_2,\partial_2u_1) |_{\partial\Omega}=0.
 \end{align}
It should be noted that the boundary condition \eqref{202220118000} implies
\begin{align}
\label{2022011091553}
 {\mm{curl}}_{\mathcal{A}}\partial_1^i\eta|_{\partial\Omega} =0\ \  \mbox{ for } \ \ i=0,  1,
\end{align}
 which will play an important role in establishing the higher-order normal estimates for $\eta$ by curl estimates.

Motivated by Jiang--Jiang--Zhao's result in \cite{JFJSZYYO}, {we shall
replace the boundary conditions of \eqref{20202201182345} and \eqref{202345} by \eqref{20safd45}
and also call the initial-boundary value problem of \eqref{01dsaf16asdfasf} and \eqref{20safd45} {\it the TCRT problem}}.
Therefore, the stability  of
the CRT problem reduces to investigating the stability of the TCRT problem. Due to the similarity of the both structures of the TCRT system and the TMRT system, we naturally  expect to follow the Jiang--Jiang--Zhao's argument in \cite{JFJSZYYO}  to verify the stability of the TCRT system. \emph{However, the capillarity  $\kappa\mathcal{K}^{\eta}$ is more complicated than the magnetic tension, and thus new ideas shall be further developed to overcome the difficulties caused by the capillarity.}

\section{Main results and  their proof strategies }\label{subsec:04}

Before stating our main results on  the TCRT problem,  similarly to \cite{JFJSZYYO},  we shall introduce simplified notations which will be used throughout this paper.
\begin{enumerate}[(1)]
  \item Simplified basic notations: $ \langle t\rangle:=1+t$,  $I_a:=(0,a)$ denotes a time interval, in particular, $I_\infty=\mathbb{R}^+:=[0,\infty)$.  $\overline{S}$ denotes the closure of a set $S\subset \mathbb{R}^n$ with $n\geqslant 1$, in particular, $\overline{I_T} =[0,T]$ and $\overline{I_\infty} = \mathbb{R}^+_0$. $\Omega_t:=\Omega\times I_t$. $(u)_{\Omega}$ denotes the mean value of $u$ in a periodic cell $ (0,2\pi L)\times (0,h)$,  {i.e. $(u)_\Omega=\frac{1}{2\pi Lh}\int  u\mathrm{d}{y}$}.
      $a\lesssim b$ means that $a\leqslant cb$ for some constant $c>0$.
      If not stated explicitly, the positive constant $c$ may depend on $\bar{\rho}$, $\mu$, $g$,  $\kappa$ and $\Omega$
      in the TCRT problem, and  vary from line to line.  The letter $\alpha$ always denotes a multi-index with respect to the variable $y$, $|\alpha|=\alpha_1+\alpha_2$ is called the order of the multi-index,
  $\partial^{\alpha}:=\partial_{1}^{\alpha_1} \partial_{2}^{\alpha_2}$ and     $[\partial^{\alpha},\phi]\varphi:=\partial^{\alpha}(\phi\varphi)-\phi\partial^{\alpha}\varphi$. We will use the Kronecker notation   $\delta_{ij}$, i.e., $\delta_{ij}=1$ for $i=j$; else $\delta_{ij}=0$. In addition, we define that
  \begin{align}
 {\underline{\tilde{\mathcal{A}}}_{ij}:={\tilde{\mathcal{A}}}_{ij}}\ \mbox{ for }\ (i,j)\neq(2,1)\ \mbox{ and  } \ {\underline{\tilde{\mathcal{A}}}_{21}:=0}.\nonumber
  \end{align}
  \item  Simplified Banach spaces, norms and semi-norms:
  \begin{align}
&L^p:=L^p (\Omega)=W^{0,p}(\Omega),\ \
{H}^i:=W^{i,2}(\Omega ), \ \ H^{j}_{\mathrm{s}}:=\{w\in {H}^{j}~|~w_2|_{\partial\Omega}=0\},\nonumber \\[1mm]
&H^{j}_\sigma:=\{w\in {H}^{j}_{\mathrm{s}}~|~\div w=0\}, \ \  H^3_{\gamma }:=\{w\in H^3_{\mm{s}}~|~ \|\nabla w\|_2 \leqslant \gamma\},\nonumber\\
& \mathcal{H}^{k}_{\mathrm{s}}:=\{w\in H^{k}_{\mm{s}}~|~ \partial_2w_1|_{\partial\Omega} =0\}, \ \
{H}^{j}_1:=\{w\in H^{j}~|~\det(\nabla w+I)=1\},\nonumber   \\
  & ^0\!{X} : =\{w\in X~|~(\bar{\rho}w_1)_\Omega=0, \   w_1\mbox{ is the first {component} of $w$}\},  \nonumber  \\
&  \underline{X}:=\{w\in X~|~(w)_{\Omega}=0\}, \  {\mathcal{H}}^k_\sigma:= \mathcal{H}^k_{\mm{s}}\cap {H}^1_\sigma, \ \
{ {\mathcal{H}}^{3,\mm{s}}_{1,\gamma}}:= {H}^3_1\cap{ \mathcal{H}^3_{\mm{s}}}\cap H_{\gamma}^3, \nonumber  \\
&
\|\cdot \|_i :=\|\cdot \|_{H^i}, \ \ \|\cdot\|_{l,i}:= \|\partial_{1}^{l}\cdot\|_{i}, \ \
\|\cdot\|_{\underline{l},i}:=\sqrt{\sum_{0\leqslant n \leqslant l}\|\cdot\|_{n,i}^2}, \nonumber
\end{align}
where $1\leqslant p\leqslant \infty$, $i$, $l\geqslant 0$,  $j \geqslant 1$, $k\geqslant 2$, $X$ stands for a general Banach space and
  $\gamma \in(0,1)$ is the constant in Lemma \ref{pro:1221}.
 \item Simplified spaces of functions with values in a Banach space:
\begin{align}
& L^p_TX:=L^p(I_T,X),\nonumber \\
& {\mathcal{U}}_{ T}: =  \{u\in C^0(\overline{I_T},  {\mathcal{H}^2_{\mm{s}}})\cap L^2_T { {H}}^3  ~|~
 u_t\in C^0(\overline{I_T} ,L^2)\cap  L^2_TH^1_{\mm{s}}\}, \nonumber\\
&   \widetilde{\mathfrak{H}}^{1,3}_{\gamma,T}:=\{\eta\in  C^0(\overline{I_T} , { \mathcal{H}^{3}_{ \mathrm{s}}}) ~|~ \eta(t) \in  {^0 {\mathcal{H}}^{3,\mm{s}}_{1,\gamma}}\;\mbox{ for each }t\in \overline{I_T}\} .
\nonumber
\end{align}
\item  Functional of potential energy:
\begin{align}\label{eeE}
E(w):=g\int\bar{\rho}'w_2^2\mm{d}y- \kappa \| \bar{\rho}'\nabla w_2\|_0^2 ,\quad w\in H^1.
\end{align}
  \item Energy and dissipation functionals (generalized):
\begin{align}
& \mathcal{E}:= \|   \eta \|_3^2+\| u\|_2^2+\|u_t\|_0^2+\| q\|_1^2,\nonumber \\
 & \mathcal{D}:= \|  \eta_1\|_{1,2}^2+\|\eta_2\|_{3}^2+\| u \|_3^2+\|u_t\|_1^2+\|q\|_2^2.\nonumber
\end{align}
We call $\mathcal{E}$ the total energy and  $\mathcal{D}$ the dissipation functionals, respectively.
 \item Other notations for decay-in-time:
 \begin{align}
 \mathfrak{E}:=\,&  \langle t\rangle  \|\partial_2^2(\partial_1 \eta_1,\partial_2\eta_2)\|_0^2  + \langle t\rangle^2 ( \|\partial_2(\partial_1\eta_1, \partial_2 \eta_2 )\|_{\underline{1},0}^2
\nonumber \\
&+\| u\|_2^2 + \|q\|_1^2) + \langle t\rangle^3( \|(\partial_2\eta_2 ,\partial_1\eta)\|_{\underline{2},0}^2+\|(\eta_2,u_t)\|^2_0)  ,  \label{2022202180904} \\[1mm]
  \mathfrak{D}:=\,&   \langle t\rangle  (\|\partial_2( \partial_1 \eta_1 ,  \partial_2\eta_2)\|_{\underline{1},0}^2 + \| u\|_{\underline{1},2}^2+\|q\|_{\underline{1},1}^2)+ \langle t\rangle^2  (\|\eta_2\|_{0}^2\nonumber \\
 &+ \|(\partial_2\eta_2,\partial_1\eta)\|_{\underline{2},0}^2  +\|u\|_1^2 )  +
 \langle t \rangle^3 ( \|  \partial_1 u\|_{\underline{1},1}^2 +\|   u_t\|_1^2 ).\label{2022202180904x}
\end{align}
\end{enumerate}

 \subsection{Stability}
Now, we state the stability result for the TCRT problem, which presents that the capillarity can inhibit the RT instability,
if the  capillary  coefficient  is properly large.
\begin{thm}[Stability]\label{thm2}
 {Let $\mu$ and $\kappa$ be positive constants, and $(\eta^0,u^0)\in{^0{\mathcal{H}}^{3,\mm{s}}_{1,\gamma}} \times {^0\mathcal{H}^2_{\mm{s}}}$ and
$\mm{div}_{\mathcal{A}^0}u^0=0$,  where $\mathcal{A}^0 =(\nabla \eta^0+I)^{-\top}$.
Assume further that $\kappa$ and $\bar{\rho}$ satisfy \eqref{0102},} the sharp stability condition \eqref{2020102241504}
  and the stabilizing condition \eqref{2022205071434}.
Then there is a sufficiently small constant $\delta >0$, such that for any $(\eta^0,u^0)$ satisfying
 $$\|(\nabla \eta^0,u^0)\|_2 \leqslant\delta ,$$
the TCRT problem  of \eqref{01dsaf16asdfasf} and \eqref{20safd45} admits a unique global strong solution $(\eta,u,q)$ in the function class $\widetilde{\mathfrak{H}}^{1,3}_{\gamma,\infty}\times {^0\mathcal{U}_{\infty} }\times (C^0(\mathbb{R}^+ _0, \underline{H}^1)\cap L^2_\infty {{H}^2})$. Moreover, the solution enjoys the stability estimate with algebraic decay-in-time: for any $t>0$,
   \begin{align}\label{1.200} \mathcal{E}(t)+ \mathfrak{E} (t)+\int_0^t(\mathcal{D}(\tau)+ \mathfrak{D}(\tau)   )\mm{d} \tau
\lesssim \| (\nabla \eta^0,u^0)\|_2^2.\end{align}
In addition, for any $t>0$,
   \begin{align}
   \label{1.safa200} \langle t\rangle^3 \|u\|_0^2
&\lesssim \| (\nabla \eta^0,u^0)\|_2^2,\\
  \langle t\rangle\|\eta_1(t)-\eta^\infty_1\|_0^2 &\lesssim    \| ( \nabla \eta^0,u^0)\|_2^2   ,\label{1.200xx}
\end{align} where $\eta_1^\infty\in H^2$ only depends on $y_2$.
\end{thm}
 \begin{rem}  It should be noted that  {in some situations $\kappa_{\mm{C}}=+\infty$. For example, we can choose $\bar{\rho}'\geqslant 0 $ such that }there exist four positive constants $\tilde{c}_1$, $\tilde{c}_2$, $s$, $\delta$ and an interval $(y_2^0-\delta,y_2^0+\delta)\subset (0,h)$ such that
 $$
\tilde{c}_1\leqslant  \frac{\bar\rho^\prime}{|y_2-y_2^0|^{2+s}}\leqslant \tilde{c}_2 \ \mbox{ for any }y_2\in  (y_2^0-\delta,y_2^0+\delta),$$
please refer to the proof of \cite[Proposition 2.1]{zhang2022rayleigh} and the construction method in \eqref{202221131926}.
  However, here  we exclude such case $\kappa_{\mm{C}}=+\infty$ by the stabilizing condition \eqref{2022205071434} since we focus on the inhibition effect of capillarity.
 \end{rem}
\begin{rem}We easily find that
\begin{align*}
0<\kappa_{\mm{C}}\leqslant  {g \|
\bar{\rho}'\|_{L^\infty}  \|
(\bar{\rho}')^{-1}\|_{L^\infty} ^2}(\pi^2h^{-2}+L^{-2})^{-1}.
\end{align*}
 {In} particular, if the RT density profile is linear, then $$\kappa_{\mm{C}}=g /(\pi^2h^{-2}+L^{-2})\bar{\rho}',$$ see Lemma \ref{261asdas567}.
 As a by-product, we  observe that the shorter the length of a periodic cell is, the greater the stabilizing effect of capillarity is. Such phenomenon can be also found in the inhibition of RT instability by the internal surface tension \cite{WYJTIKCT,JJTIWYJ2016ARMA}, but {does't} exist in the inhibition of RT instability by the magnetic tension \cite{JFJSZYYO,WYJ2019ARMA}.\end{rem}
\begin{rem}\label{202201271239}
For each fixed $t\in \mathbb{R}_0^+$, the solution $\eta(y,t)$ in Theorem \ref{thm2} belongs to $ H^3_{\gamma}$. Let $\zeta=\eta+y$,
then $\zeta $  satisfies \eqref{20210301715x} and \eqref{20210301715} for each $t\in \mathbb{R}_0^+$ by virtue of Lemma \ref{pro:1221}.
We denote the inverse transformation of $\zeta$ by $\zeta^{-1}$ and {define}
\begin{align}
(\varrho, v, \beta)(x,t):=(\bar{\rho}(y_2)-\bar{\rho}(\zeta_2), u(y,t),  q(y,t))|_{y=\zeta^{-1}(x,t)}.
\label{2022201271250}
\end{align}
 {Then we can easily verify} that $(\varrho,v, \beta)$ is a strong solution of the CRT problem  \eqref{0103}--\eqref{0105}
and enjoys the stability estimate similar to \eqref{1.200}  for sufficiently small $\delta$.
\end{rem}
\begin{rem}
In Theorem \ref{thm2}, we have assumed
\begin{align*}
 (\bar{\rho}\eta^0_1)_{\Omega}=(\bar{\rho}u^0_1)_{\Omega}=0 .
\end{align*}
If $(\bar{\rho}\eta^0_1)_{\Omega}$, $(\bar{\rho} u^0_1)_{\Omega} \neq 0$, we can  define
$\bar{\eta}^0_1:=\eta^0_1-(\bar{\rho}\eta^0_1)_\Omega (\bar{\rho})_\Omega^{-1}$, $\bar{u}^0_1:=u^0_1-(\bar{\rho}u^0_1)_\Omega(\bar{\rho})_\Omega^{-1}$ and
$( \bar{\eta}^0_2,\bar{u}^0_2):=( {\eta}^0_2, {u}^0_2)$. Then, by virtue of Theorem \ref{thm2}, there exists a unique global strong solution $(\bar{\eta},\bar{u}, {q})$ to the TCRT problem   with the  initial data $( \bar{\eta}^0,\bar{u}^0)$. It is easy to verify that
$(\eta_1,\eta_2,u_1,u_2,q):=(\bar{\eta}_1+ t (\bar{\rho}u^0_1)_\Omega(\bar{\rho})_\Omega^{-1} +  (\bar{\rho}\eta^0_1)_\Omega(\bar{\rho})_\Omega^{-1},\bar{\eta}_2, \bar{u}_1+( \bar{\rho}u^0_1)_\Omega(\bar{\rho})_\Omega^{-1},\bar{u}_2, {q})$
{is the} unique strong solution of  the TCRT problem  with the initial data $({\eta}^0,{u}^0)$.
\end{rem}
\begin{rem}
 {Additionally, if} the initial data $(\eta^0,u^0)$ in Theorem \ref{thm2} satisfies the odevity conditions
$ (\eta_1^0,u_1^0)(y_1,y_2)=-(\eta_1^0,u_1^0)(-y_1,y_2)$ and  $(\eta_2^0,u_2^0)(y_1,y_2)=(\eta_2^0,u_2^0)(-y_1,y_2)$,
then we have  $\eta^\infty_1=0$ in \eqref{1.200xx}.
\end{rem}
\begin{rem}  It should be noted that our proof for Theorem \ref{thm2} strongly depends on the two-dimensional structure, and thus can not be extended to the  corresponding 3D case. However, based on the linear 3D stability result in \cite{LFCZZP}, {we believe} that the  inhibition phenomenon of RT instability by capillarity should be verified by the 3D NSK system. Such conjecture will be further investigated in future.
\end{rem}

 Now, we sketch the proof of Theorem \ref{thm2},  which is an adaptation of the argument from \cite{JFJSZYYO} by further introducing   new ideas.
The key step in the existence proof for global small solutions is to derive the \emph{a priori} uniform-in-time
energy estimate \eqref{1.200}. To this end, let $(\eta,u)$ be a solution of the TCRT problem, satisfying that for some $T>0$,
\begin{gather}
 (\bar{\rho} \eta_1)_\Omega\equiv 0\ \mbox{ for any }\ t\in \overline{I_T},\label{apresnew}  \\
 \det(I+\nabla \eta)=1\ \mbox{ in }\ \Omega\times \overline{I_T},\label{aprpiosasfesnew}\\
 \sup_{t\in\overline{I_T}} \|(\nabla \eta,u)(t)\|_2\leqslant {\delta} \in (0,1] . \label{aprpiosesnew}
\end{gather}

Under the above \emph{a priori} assumptions,  Jiang--Jiang--Zhao \cite{JFJSZYYO} derived that the TMRT problem enjoys the following energy inequality
\begin{align}
\label{for:0202sadfn}
\frac{\mm{d}}{\mm{d}t}\tilde{{E}}+ \mathcal{D}\leqslant 0,
\end{align}
where  $\tilde{{E}}$ is equivalent to   {$\mathcal{E}$}
 under a sharp stability condition of the TMRT problem.  Clearly  the above energy inequality immediately implies the global existence of stable solutions.

However, directly following the  Jiang--Jiang--Zhao's argument, we will find that the TCRT  problem enjoys the following energy inequality
\begin{align}
\frac{\mm{d}}{\mm{d}t}\tilde{\mathcal{E}}+ c\mathcal{D}\lesssim\,&    \|  \eta_2\|_{2,1}\| \eta\|_3^2+(
\|\partial_2\eta_2\|_{ {1},0}^{1/2}\|\partial_2^2\eta_2\|_{0}^{1/2}
\nonumber \\
& +\|\partial_2\eta_2\|_{\underline{1},0}
)(
\|\partial_2^2\eta_1\|_{ {1},0}^{1/2}\|\partial_2^3\eta_1\|_{0}^{1/2}+\|\partial_2^2\eta_1\|_{\underline{1},0}
)
\label{for:0202n}
\end{align}
for some energy functional $\tilde{\mathcal{E}}$, which is equivalent to $\mathcal{E}$ under \emph{the   stability conditions \eqref{2020102241504} and \eqref{2022205071434}}. In particular, integrating \eqref{for:0202n} over $(0,t)$ yields the following energy estimate
\begin{align}
\mathcal{E}(t)+ \int_0^t {\mathcal{D}(\tau)}\mm{d}\tau
\lesssim\, & \|(\nabla\eta^0,u^0)\|_2^2+ \int_0^t ( \|  \eta_2\|_{2,1}\| \eta\|_3^2 +(
\|\partial_2\eta_2\|_{ {1},0}^{1/2}\|\partial_2^2\eta_2\|_{0}^{1/2}
\nonumber \\
&+\|\partial_2\eta_2\|_{\underline{1},0}
)(
\|\partial_2^2\eta_1\|_{ {1},0}^{1/2}\|\partial_2^3\eta_1\|_{0}^{1/2}+\|\partial_2^2\eta_1\|_{\underline{1},0}
))\mm{d}\tau . \label{1.200xyx}
\end{align}
We mention that all terms at the right hand side of the inequality \eqref{for:0202n} are related with  the capillarity.
 {Obviously, we can not directly expect  the closure of energy estimates from \eqref{1.200xyx}}, and thus shall further dig new energy estimates from the equations in \eqref{01dsaf16asdfasf}.

To begin with, let us recall the basic energy equality, which is satisfied by the solution $(\varsigma,v)$ for the linearized TCRT   problem:
\begin{align}
\frac{1}{2}\frac{\mm{d}}{\mm{d}t}\left(\|\sqrt{\bar{\rho}}v\|_0^2 -E( \varsigma)\right)+ \|\nabla v\|_0^2=0. \label{1.200sadfgxyx}
\end{align}
For  simplicity, we call $(\varsigma,v)$ the linear solution of the TCRT problem (i.e. a solution of the linearized TCRT problem), which enjoys basic energy estimate by \eqref{1.200sadfgxyx}:
$$\|( \varsigma_2, \nabla \varsigma_2 )(t)\|_0^2+ \|v(t)\|_0^2+\int_0^t\|\nabla v\|_0^2\lesssim \|(\varsigma^0_2, \nabla \varsigma^0_2)\|_0^2+\|v^0 \|_0^2.$$
 Similarly, the  linear solution, also denoted by $(\varsigma,v)$, of the TMRT problem has the following basic energy estimate:
$$\|(\varsigma_2, \partial_1\varsigma)(t)\|_{0}^2+ \|v(t)\|_0^2+\int_0^t\|\nabla v\|_0^2\lesssim \|(\varsigma^0_2,\partial_1\varsigma^0)\|_0^2+\|v^0\|_0^2.$$

Noting that $\mm{div}\varsigma=0$   in the linearized problems, thus it is easy to see  that two norms
$\|\nabla \varsigma_2(t)\|_0$ and $\|\partial_1\varsigma(t)\|_{0}$ are equivalent to each other.
Motivated by the equivalence of the two norms and the other fact that the stable solutions of the linearized TMRT problem enjoy the
algebraic decay-in-time \cite{JFJSZYYO},  {we guess} that the solution $\varsigma$ of the linearized TCRT problem   has the same decay-in-time as well as  the one of the  linearized TMRT problem, i.e.
\begin{align}
  & \langle t\rangle  \|\partial_2^2(\partial_1 \varsigma_1,\partial_2\varsigma_2)\|_0^2  + \langle t\rangle^2  \|
  \partial_2( \partial_1\varsigma_1,\partial_2\varsigma_2 )\|_{\underline{1},0 }^2
 + \langle t\rangle^3(\| \varsigma_2 \|_0^2 +\|(\partial_2\varsigma_2,\partial_1\varsigma)\|_{\underline{2},0}^2 ) \nonumber  \\[1mm]
 &+\int_0^t(\langle\tau\rangle   \|\partial_2(\partial_1\varsigma_1, \partial_2\varsigma_2 )\|_{\underline{1},0}^2 +   \langle \tau\rangle^2  ( \| \varsigma_2 \|_{0}^2
 + \|(\partial_2\varsigma_2,\partial_1\varsigma)\|_{\underline{2},0}^2   ) )\mm{d}\tau \lesssim \| (\nabla \varsigma^0,v^0)\|_2^2 . \label{202220dsfad2180904}
\end{align}
  With the decay-in-time of $\varsigma$ in hand, we  {can further} derive the  decay-in-time of $v$ and the associated pressure function $\sigma$ from the linearized TCRT problem:
   \begin{align}
 & \langle t\rangle^2 (  \| v\|_2^2 + \|\sigma\|_1^2) + \langle t\rangle^3 \|v_t\|^2_0
+ \int_0^t (\langle \tau\rangle  (  \| v\|_{\underline{1},2}^2\nonumber  \\
&+\|\sigma\|_{\underline{1},1}^2) +
 \langle \tau \rangle^2\|v\|_1^2+
 \langle \tau \rangle^3 ( \|  \partial_1 v\|_{\underline{1},1}^2 +\|   v_\tau\|_1^2 )) \mm{d}\tau\lesssim \| (\nabla \varsigma^0,v^0)\|_2^2  .  \label{saf}
\end{align}

Based on the above linear analysis, we naturally believe the solution $(\eta,u)$ of the TCRT problem also possesses the above two decay-in-time estimates as well as the linear solution  and thus introduce the two notations involving decay-in-time  {$\mathfrak{E}$ in \eqref{2022202180904} and   $\mathfrak{D}$  in  \eqref{2022202180904x}.} Moreover, motivated by the decay-in-time of $\eta$, we derive from \eqref{1.200xyx} that
\begin{align}
&\sup_{0\leqslant \tau\leqslant t} \mathcal{E}(\tau)+ \int_0^t{\mathcal{D}(\tau)}\mm{d}\tau\nonumber \\
&\lesssim  \|(\nabla\eta^0,u^0)\|_2^2+ \sup_{0\leqslant \tau\leqslant t} \left(\langle\tau\rangle^{3/2}   \|  \eta_2(\tau)\|_{2,1}  \mathcal{E}(\tau) + \sqrt{\mathfrak{E}(\tau)
\|\partial_2^2\eta_1(\tau)\|_{\underline{1},0} }
 \mathcal{E}^{1/4}(\tau)\right). \label{as1.200xyx}
\end{align}

  Obviously, we shall further  {estimate   $\mathfrak{E}$  and  $\|\partial_2^2\eta_1\|_{\underline{1},0}^{1/2} $, see \eqref{20222010041821} and  \eqref{20200asd46safdxx}, respectively}.
Since the estimate  of $\|\partial_2^2\eta_1\|_{\underline{1},0}^{1/2} $ can be easily obtained by removing the  third-order normal estimates of $\eta$ in \eqref{for:0202n}, thus we next emphatically mention the estimate for $\mathfrak{E}$.

Because Jiang--Jiang--Zhao have established the decay-in-time estimates for the stable solutions of the TMRT problem in   \cite{JFJSZYYO}, we also  expect to follow their argument to derive the desired decay-in-time estimates. However we will face two difficulties:
 \begin{enumerate}[(1)]
   \item Since the linear part of capillarity includes a two-order normal derivative,  the decay-in-time of the solution $(v,q)$  {to} the TCRT problem is slower than the one of the  TMRT problem. Such slower decay-in-time behavior results  in that the nonlinear integrals (i.e. the integrand on the domain $\Omega$ is the multiplication of at least three functions such as $\eta$, $u$ and $q$) involving pressure can not be
             controlled, if we follow directly the approach in \cite{JFJSZYYO} to derive the   decay-in-time estimates involving $\langle t\rangle^{3}$.
             {Thus, we must develop new ideas.}
   \item   {There} always exists  the derivation $\partial_2\eta_1$, which does not decay in time, in some nonlinear terms involving the capillarity. This also results in the trouble of deriving the   decay-in-time  estimates involving $\langle t\rangle^{3}$.
 \end{enumerate}

 To overcome the first difficulty, we exploit the fine properties of structures from the expressions of $\mm{div}\eta$, $\mm{div}_\mathcal{A}u$ and $\partial_1 p$ to deal with the nonlinear integrals  involving pressure, and   capture  more tangential estimates, which enjoy quicker decay-in-time, see the derivations of \eqref{202008241550} and \eqref{202220041}.
As for the second difficulty, we will  utilize the structure of energy integrals, the relation  \eqref{01dsafdasf16asdfasf}$_1$ and the incompressible condition to rewrite the nonlinear integrals involving $\partial_2\eta_1$ into the sum of the other nonlinear integrals, which exclude $\partial_2\eta_1$, and the temporal derivative of  integrals involving  $\partial_2\eta_1$, see the estimate for $I_{8,i} $ in the proof of Lemma \ref{lem:08241445} for examples.  After overcoming the two difficulties, we can build the desired  decay-in-time estimates in $\mathfrak{E}$.

  Adding the energy  estimate of $\mathfrak{E}$ mentioned above to \eqref{as1.200xyx}, and then using Young's inequality and the estimate of $\|\partial_2^2\eta_1\|_{\underline{1},0}^{1/2} $, we can arrive at the following energy estimate
   \begin{align}
   & \sup_{0\leqslant \tau\leqslant t} (\mathcal{E}(\tau)+ \mathfrak{E} (\tau))+\int_0^t(\mathcal{D}(\tau)+ \mathfrak{D}(\tau)   )\mm{d} \tau
\nonumber \\
\lesssim\,& \|(\nabla \eta^0,u^0)\|_2^2  + \sup_{0\leqslant \tau\leqslant t}(\langle\tau\rangle^{3/2}   \|  \eta_2(\tau)\|_{2,1}  \mathcal{E}(\tau) +\langle \tau\rangle^3  \|u(\tau)\|^4_2 )
\nonumber \\
& +\int_0^t ( {\|\eta_2\|_{2,1}\|\eta\|_1\|\eta\|_3}+ \mathcal{I})\mm{d}\tau ,\label{1.2asdf00}\end{align}
 see \eqref{202222100718123} for the   definition of $\mathcal{I}$.

By carefully observing  the expression of the integrand $\mathcal{I}$, we easily  impose the additional \emph{a priori} assumption
\begin{align}
\label{202220201915}
 \mathfrak{A}:=&  \langle t\rangle ( \|\partial_2^2\eta_1\|_{1,0}^2+\|\partial_2^3\eta_2\|_0^2 )  + \langle t\rangle^2 (\|  \partial_2 (\partial_1\eta_1,\partial_2 \eta_2)\|_{\underline{1},0}^2
\nonumber \\
&+\| u\|_2^2 ) + \langle t\rangle^3(\| \eta_2\|_0^2 +\|(\partial_2\eta_2,\partial_1\eta)\|_{\underline{2},0}^2 )     \leqslant \delta^2\ \mbox{ for any } \ t\in \overline{I_T}.
\end{align}
Under  the above \emph{a priori} assumption, the last two terms at the right hand side of the inequality \eqref{1.2asdf00} can be {absorbed} by the two terms at the left hand side, and thus we close the energy estimate (or arrive at \eqref{1.200}).

Thanks to \eqref{1.200} and the unique local (-in-time) solvability of the TCRT problem in Proposition \ref{202102182115}, we further obtain the unique global solvability immediately. Exploiting the parabolic structure of \eqref{01dsafdasf16asdfasf}$_2$ and \eqref{1.200}, we can further improve the decay-in-time rate of $\|u\|_0^2$ and thus get \eqref{1.safa200}. Finally, we easily obtain \eqref{1.200xx} from \eqref{1.200} and \eqref{1.safa200} by an asymptotic analysis method.

  \subsection{Instability}
We point out that  we can not expect a stability result  for the TCRT problem under the condition
$\kappa\in [0,\kappa_{\mm{C}})$. In fact, we have the following instability result.
\begin{thm}[Instability]\label{thm1}
 {Let $\mu$ and $\kappa$ be positive constants, and $\bar{\rho}$ satisfy \eqref{0102} and \eqref{0102n}.} If   $\kappa\in [0,\kappa_{\mm{C}})$,
then the {equilibrium}  $(\bar{\rho},0 )$ is unstable in the Hadamard sense, that is, there are positive constants $\varpi $,
$\epsilon$, $\delta_0$  and $(\tilde{\eta}^0,\eta^\mm{r},\tilde{u}^0,u^\mm{r})\in  {^0\!\mathcal{H}^3_{\mm{s}}}$,
such that for any $\delta\in (0,\delta_0]$ and the initial data
 $$ (\eta^0, u^0):=\delta(\tilde{\eta}^0,\tilde{u}^0)  +\delta^2(\eta^\mm{r},u^\mm{r}), $$
there exists a unique strong solution $(\eta,u,q)$ to the TCRT problem of  \eqref{01dsaf16asdfasf} and \eqref{20safd45}
 on $[0,T^{\max})$,
where $(\eta,u,q)\in \widetilde{\mathfrak{H}}^{1,3}_{\gamma, \tau  }\times {^0\mathcal{U}_{\tau}}
\times (C^0(\overline{I_\tau},\underline{H}^1)\cap L^2_\tau {{H}^2})$ for any given $\tau\in I_{T^{\max}}$,
 and $T^{\max}$ denotes the maximal time of existence of the solution. {Moreover}, the solution satisfies
\begin{align}
&\| \bar{\rho}(y_2)-\bar{\rho}(\eta_2(y, T^\delta)+y_2) \|_{L^1},\ \ \|\chi_i(T^\delta)\|_{L^1}, \ \ \| \partial_1\chi_{i} (T^\delta)\|_{L^1},
\nonumber  \\
& \|\partial_2\chi_i(T^\delta)\|_{{L^1}}, \ \ \| \mathcal{A}_{1k}\partial_k \chi_{i} (T^\delta)\|_{L^1}, \ \ \|\mathcal{A}_{2k}\partial_k\chi_i(T^\delta)\|_{{L^1}}  \geqslant  \epsilon \label{201806012326}
\end{align}
for some escape time $T^\delta:= {\Lambda}^{-1}\mm{ln}({2\epsilon}/{\varpi \delta})\in I_{T^{\max}}$,
where $i=1,2$, and $\chi=\eta$ or $u$.
\end{thm}
 \begin{rem}
It should be noted that  $\kappa_{\mm{C}}>0$ due to $\bar{\rho}$ satisfying \eqref{0102n} and
 $\kappa_{\mm{C}}$ is admitted to be infinite   in the above theorem.
 \end{rem}
\begin{rem}
By the inverse transformation of Lagrangian coordinates in \eqref{2022201271250} and
 the inequalities in \eqref{201806012326}, we easily obtain the corresponding instability in Eulerian coordinates:
for $i=1,2$,
$$\|\varrho(T^\delta)\|_{L^1},\ \|v_i(T^\delta)\|_{L^1},\ \|\partial_1 v_i(T^\delta)\|_{L^1},\
\|\partial_2 v_i(T^\delta)\|_{L^1} \geqslant  \epsilon  . $$
In addition, the above  2D instability result  automatically implies the corresponding 3D one, in which the 3D slab domain is horizontally periodic with respect to the two horizonal components of spacial variables respectively.
\end{rem}

We can follow the  standard  bootstrap instability method in \cite{JFJSZYYO} (interesting reader can refer to the other versions  in \cite{JFJSO2014,GYSWIC,GYSWICNonlinea,FSSWVMNA,GYHCSDDC,JFJSZWC}) to easily establish Theorem \ref{thm1}.  For the completeness we  will briefly sketch the  proof in Section \ref{sec:instable}.

\medskip
The rest of this paper is organized as follows. In Sections \ref{sec:global} and \ref{sec:instable}, we provide the proof of
Theorems \ref{thm2} and  \ref{thm1}, respectively.  Finally, in Appendix, we list some mathematical results, which will be used
in Sections \ref{sec:global}--\ref{sec:instable}.

\section{Proof of Theorem \ref{thm2}}\label{sec:global}
This section is devoted to the proof of Theorem \ref{thm2}. The key step in the proof is to establish the energy estimate with decay-in-time \eqref{1.200}  for the TCRT problem  \eqref{01dsaf16asdfasf} and \eqref{20safd45}.
To this end, let $(\eta,u,q)$ be a solution of the TCRT problem  and satisfy \eqref{apresnew}--\eqref{aprpiosesnew},
where ${\delta} $ is sufficiently small, and the smallness of $\delta$ depends on $\bar{\rho}$,  $\mu$, $g$,  $\kappa$ and $\Omega$.
It should be noted that  $\bar{\rho}$  and $\kappa$ {satisfy} the assumptions in Theorem \ref{thm2}.
 \subsection{Preliminary estimates}
    We start with deriving some preliminary estimates involving $(\eta,u)$.
\begin{lem}
\label{201805141072}
 For any given $t\in \overline{I_T}$, we have:
\begin{enumerate}[(1)]
\item  Estimates of $ {\mm{div}\, { {\eta}}}$:
\begin{align}
 \label{improxtian1}
\|\mm{div} \,{ {\eta}}\|_{i} &\lesssim \|\eta\|_{1,i}\| \eta\|_3 \quad
\mbox{ for }\,\,0\leqslant i\leqslant  2,\\
 \label{improtian1}
\|\mm{div} \,{ {\eta}}\|_{i,0}& \lesssim
\begin{cases}
\|\eta\|_{1+i,0}\| \eta\|_3 &   \mbox{for }\,\,\,i=0, 1;\\
 \|\eta\|_{3,0}\| \eta\|_3+ \|\eta_1\|_{2,1} \|\eta_2\|_{2,1}  &  \mbox{for }\,\,\,i=2, \end{cases} \\
   \|\mm{div}\eta\|_{1,1}& \lesssim \|\eta\|_{2,1} \| \eta\|_3.
 \label{2022210071415}
 \end{align}
 \item Poincar\'e's inequalities for $\eta$, $u$ and $u_t$:  for $i=1$, $2$ and $0\leqslant j\leqslant 2$,
\begin{align}
\|  \eta_i\|_{j,1}& \lesssim \|\nabla \eta_i\|_{j,0 },  \label{202211}\\
\|    u_i\|_{j,1}&\lesssim \| \nabla u_i\|_{j,0}, \label{2022111522311}\\
\|  \partial_t u_i\|_1&\lesssim \|\nabla \partial_t u_i\|_0. \label{2022220111522311}
\end{align}
\item Curl estimates:  {let $0\leqslant i \leqslant 2$}, then for sufficiently small $\delta$,
\begin{align}
\|\eta\|_{i,3-i}\lesssim\|\mm{curl}\eta\|_{i,2-i}. \label{2022202011749}
\end{align}
\item Estimates of $(\partial_1\eta_1,\partial_2\eta_2,\partial_2u_2)$: for sufficiently small $\delta$,
\begin{align}
\label{omessen}
 \| \partial_2^j \eta_2 \|_{i,0}&\lesssim
\|   \eta\|_{1+i,j-1}\quad\mbox{ for }\,\,1\leqslant i+j\leqslant 3\ \mbox{ and }\ 1\leqslant j ,\\
\label{omessetsim122n}
 \| \eta_1 \|_{1+i,j}&\lesssim
\|  \eta_2\|_{i,1+j}\quad \,\,\,\mbox{ for }\,\, 0\leqslant i+j\leqslant 2,\\
\label{omesssafden}
 \| \partial_2 u_2 \|_{0} &\lesssim
\|   u\|_{1,0}+\|\eta_2\|_{2,1}\|\partial_2 u_1\|_0   .
\end{align}
\end{enumerate}
\end{lem}
\begin{rem}
In view of \eqref{202211}, \eqref{omessen}, \eqref{omessetsim122n} and \eqref{202012241002}, it is easy to check that
\begin{align}
&\label{osdfadsa}
 \|  \partial_1\eta_1 \|_{\underline{i},j-1}\lesssim  \|  \eta_2 \|_{\underline{i},j}\lesssim
\|   \eta_2\|_{i,j}\mbox{ for }1\leqslant i+j\leqslant 3\mbox{ and }1\leqslant j   .
\end{align}
\end{rem}
\begin{pf}
 We only estimate  {\eqref{2022210071415} and  \eqref{omessen}--\eqref{omesssafden}} in sequence, since the rest estimates can be found in   Lemmas 2.1 and 2.2 in \cite{JFJSZYYO}  under the conditions \eqref{apresnew}--\eqref{aprpiosesnew}.

 Recalling \eqref{aprpiosasfesnew}, we see that
\begin{align}
&\mm{div}\eta=\partial_2\eta_1\partial_1\eta_2-\partial_1\eta_1\partial_2\eta_2, \label{202209151033}
\end{align}  Differentiating the above identity with  {respect} to $y_1$ yields that
\begin{align}
&\partial_1 \mm{div}\eta=\partial_2\eta_1\partial_1^2\eta_2+\partial_2\partial_1\eta_1\partial_1\eta_2
-\partial_1^2\eta_1\partial_2\eta_2-\partial_1\eta_1\partial_2\partial_1\eta_2.\nonumber
\end{align}
Applying the norm $\|\cdot\|_1$ to the above identity, and then using
\eqref{fgestims}, \eqref{fgessfdims} and \eqref{202012241002}, we immediately get \eqref{2022210071415}.

By \eqref{202209151033}, we have
\begin{align}\nonumber
 \|   \partial_2^j \eta_2\|_{i,0}\lesssim \| \eta_1 \|_{1+i,j-1} +\|\mm{div} \eta \|_{i,j-1} .
\end{align}
Exploiting \eqref{improxtian1}--\eqref{2022210071415}, we further derive \eqref{omessen}  from the above estimate for sufficiently small $\delta$.

Similarly, we have
\begin{align}\nonumber
 \|   \eta_1 \|_{1+i,j}\lesssim \| \eta_2 \|_{i,1+j} +\|\mm{div} \eta \|_{i,j} .
\end{align}
Exploiting \eqref{improxtian1}--\eqref{2022210071415}, we also obtain \eqref{omessetsim122n}.

Recalling  \eqref{01dsafdasf16asdfasf}$_3$, we have
\begin{align}
\label{202222010109}
\partial_2 u_2=-\partial_1u_1-\tilde{\mathcal{A}}_{jk}\partial_k u_j.
 \end{align}
Applying $\|\cdot\|_0$ to the above identity yields
\begin{align}
\label{20222saf2010109}
\|\partial_2 u_2\|_0\lesssim
\| u_1\|_{1,0}+\|\tilde{\mathcal{A}}_{jk}\partial_k u_j\|_0\lesssim
\|u\|_{1,0}+\|\partial_1\eta_2\|_{\underline{1},1}\|\partial_2 u_1\|_0+\|\eta\|_3\|\partial_2 u_2\|_0,
 \end{align}
 which, together with \eqref{osdfadsa}, yields
 \eqref{omesssafden}.
\hfill $\Box$
\end{pf}
\begin{lem}\label{lem:202012242115}
 For any given $t\in \overline{I_T}$, we have
\begin{enumerate}[(1)]
\item Estimates involving the gravity term: for sufficiently small $\delta$,
\begin{align}
\left\| \mathcal{G}\right\|_{\underline{1},0}&\lesssim \|\eta_2\|_0\|\eta_2\|_{ 2,1} , \label{2022011130957}\\
\left\| \mathcal{G}\right\|_{1}&\lesssim \|  \eta_2 \|_{ 0} ,
\label{2022201201821}
\end{align}
where $\mathcal{G}:=G^{\eta} -g \bar{\rho}'\eta_2$.
\item  Estimates involving the capillarity term: let $1\leqslant j$, $k$, $l\leqslant 2$,  for sufficiently small $\delta$,
\begin{align}
\left\| \mathcal{N}^{\eta,1 } _{j,k} \right\|_{i,0}&\lesssim
  \begin{cases}
  \|  \eta_2 \|_{1}\| \eta  \|_3  &\ \mbox{for }\,\,i=0;
  \\ \| \eta_2 \|_{  {i},1}\|  \eta  \|_3  + \|\partial_2 \eta_1 \|_{i, 0} \|  \eta_2 \|_{ {1},2}&\mbox{for }\,\,i=1,\ 2,
 \end{cases}\label{202220927900} \\
\left\| \mathcal{N}^{\eta,2 } _{j,k} \right\|_{i,0}&\lesssim
 \begin{cases}
\|  \eta_2 \|_1\|\eta_2 \|_{3}  &\ \mbox{for }\,\,i=0; \\
  \|  \eta_2 \|_{ {1},2}^2  + \|  \eta_2 \|_{ {i},1} \|  \eta_2 \|_{3}&\ \mbox{for }\,\,i=1,\ 2,
  \end{cases}\label{2022dfdsfs2} \\
\left\|\partial_k\mathcal{N}^{\eta,l } _{j,k} \right\|_{i,0}&\lesssim  \|  \eta_2 \|_{i,2}  \|  \eta_2\|_3+\|  \eta_2 \|_{ i+1,1} \| \eta\|_{3}\quad \ \mbox{ for } \,\, i  =0,\ 1,
\label{2022dsfs2} \\
\left\|\partial_2(\mathcal{N}^{\eta,1 } _{1,1},\mathcal{N}^{\eta,1 } _{2,2} ) \right\|_{1,0}&\lesssim  \|\eta_2\|_{2,1}\|\eta\|_3,
\label{2022dsfasdas2} \\
\left\|\nabla \mathcal{N}^{\eta,l } _{j,k} \right\|_{m,n}
&
\lesssim  \|  \eta_2 \|_{m,2+n}  \|  \eta\|_3+\|  \eta_2 \|_{1,2}  \|  \eta\|_{m,2+n}\quad \mbox{ for }\,\,0\leqslant m+n\leqslant 1, \label{3.21}\\
\left\|\nabla \mathcal{N}^{\eta,2 } _{j,k} \right\|_1
&\lesssim  \|   \eta_2 \|_3^2, \label{2022210220847} \\
\left\| \partial_t \mathcal{N}^{\eta,l } _{j,k} \right\|_{0}
&\lesssim    \|  \eta \|_{3} \| u_2 \|_{1} + \|  \eta_2 \|_{ 1,1} \| u \|_2, \label{2022201010520233} \end{align}
where we have defined that
\begin{align}
 \mathcal{N}^{\eta,1} _{j,k} :=&
  \begin{cases}
  \kappa|\bar{\rho}' |^2( {\tilde{\mathcal{A}}}_{2k} - \underline{\tilde{\mathcal{A}}}_{2k}) \partial_1\eta_2   &\ \mbox{for }\,\,j=1;\\
  \kappa  |\bar{\rho}'   |^2 ( \tilde{\mathcal{A}}_{21}\partial_1\eta_2 ({\mathcal{A}}_{ik} {\mathcal{A}}_{i2}    +  {\mathcal{A}}_{2k} )
    +  ( {\tilde{\mathcal{A}}}_{2k} - \underline{\tilde{\mathcal{A}}}_{2k}) \partial_2\eta_2 )   &\ \mbox{for }\,\,j=2,
  \end{cases} \label{202220926622063}
  \\
\mathcal{N}^{\eta,2} _{j,k} :=\,& \kappa ({\mathcal{A}}_{ik}(|\bar{\rho}'(y_2+\eta_2)|^2  {\mathcal{A}}_{il}\partial_l\eta_2 {\mathcal{A}}_{jm} \partial_m\eta_2 +(|\bar{\rho}'(y_2+\eta_2) |^2 \nonumber \\
 & -|\bar{\rho}' |^2 ) ({\mathcal{A}}_{i2}( {\mathcal{A}}_{j1}\partial_1\eta_2  + 2{\mathcal{A}}_{j2} \partial_2\eta_2 )+{\mathcal{A}}_{i1} {\mathcal{A}}_{j2} \partial_1\eta_2    ) \nonumber \\
 &+ |\bar{\rho}'  |^2  ({\mathcal{A}}_{i2}  ( \underline{\tilde{\mathcal{A}}}_{j1} \partial_1\eta_2  + 2\tilde{\mathcal{A}}_{j2} \partial_2\eta_2 ) +\tilde{\mathcal{A}}_{i2}(\delta_{j1}\partial_1\eta_2 + 2\delta_{j2} \partial_2\eta_2 ) \nonumber \\
 &+({\mathcal{A}}_{i1} \tilde{\mathcal{A}}_{j2}  +
 \underline{\tilde{\mathcal{A}}_{i1}}\delta_{j2}) \partial_1\eta_2   )+ {\mathcal{A}}_{i2}  {\mathcal{A}}_{j2}\mathcal{R}+(|\bar{\rho}'  |^2)'( {\mathcal{A}}_{i2}  \tilde{\mathcal{A}}_{j2}\nonumber \\
 & + \tilde{\mathcal{A}}_{i2}  \delta_{j2} ) \eta_2 ) + \underline{\tilde{\mathcal{A}}}_{ik} ( |\bar{\rho}'  |^2(\delta_{i2}(\delta_{j1}\partial_1\eta_2 + 2\delta_{j2} \partial_2\eta_2)+ \delta_{i1}\delta_{j2} \partial_1\eta_2) \nonumber
  \\
 &+(|\bar{\rho}' |^2)' \delta_{i2} \delta_{j2} \eta_2    ) - \underline{\tilde{\mathcal{A}}}_{jk}\partial_2(|\bar{\rho}' |^2\eta_2 ) ) \label{2022209sadf26622063} \end{align}
 and  \begin{align}\mathcal{R}:={ \int_{0}^{\eta_2}\left(\eta_2(y,t) -
z\right)(|\bar{\rho}' (s)|^2)''|_{s=y_2+z}\mm{d}z}.
\label{20222042015205}
\end{align}
\item  Estimate involving the viscosity term:   for $0\leqslant i\leqslant 2$, $l=1$, $2$  and $k=0$, $1$,
    \begin{align}
\|\mathcal{N}^\mu_{j,l}\|_{ i,0}&\lesssim  (1-\delta_{i0}) \|  \eta\|_{i,1}\|u\|_{\underline{1},2} +\| \eta\|_3\|u\|_{i,1}  , \label{20asd10101643}\\
\|\partial_2\mathcal{N}^\mu_{j,2}\|_{ k,0}&\lesssim \| \eta\|_{1+k,1} \|  u\|_{\underline{1},2}+  \| \eta\|_3\|  u\|_{1+k,1} , \label{20asd101asd01643}\\
\|\mathcal{N}^\mu\|_{k,0}&\lesssim \|\eta\|_{1+k,1}\|u\|_{\underline{1},2}+\| \eta\|_{\underline{1},2 }\|u\|_{1+k,1} ,\label{20222010101643}\\
\|\mathcal{N}^\mu\|_{k}&\lesssim \|\eta\|_{3}\|u\|_{k+2}, \label{202220110230924}
    \end{align}
    where we have defined that
    \begin{align*}
\mathcal{N}^\mu=\, & \partial_l(\mathcal{N}^\mu _{1,l}, \mathcal{N}^\mu _{2,l})^{\top} ,
\\
\mathcal{N}^\mu _{j,1}:  =\,& \mu (\mathcal{A}_{k1}\tilde{\mathcal{A}}_{km}+\tilde{\mathcal{A}}_{m1} )\partial_mu_j
 = \mu ( ( 2\partial_2\eta_2+ (\partial_2\eta_1)^2+ (\partial_2\eta_2)^2 )\partial_1u_j - \Theta\partial_2u_j), \\
\mathcal{N}^\mu _{j,2} : =\,&  \mu (\mathcal{A}_{k2}\tilde{\mathcal{A}}_{km}+\tilde{\mathcal{A}}_{m2} )\partial_mu_j
= \mu ((2\partial_1\eta_1 + (\partial_1\eta_1)^2+(\partial_1\eta_2)^2 )\partial_2u_j  -\Theta\partial_1u_j),\\
\Theta: = \,&  \partial_1\eta_2+\partial_2\eta_1 + \partial_1\eta_2\partial_2\eta_2  + \partial_1\eta_1 \partial_2\eta_1 .
\end{align*}
\end{enumerate}
\end{lem}

\begin{pf}(1) By virtue of \eqref{aprpiosesnew} and Lemma \ref{pro:1221}, {$\zeta=\eta+y$} satisfies the diffeomorphism properties \eqref{20210301715x} and \eqref{20210301715} for sufficiently small $\delta$. Thus,  {for any $y\in\overline{\Omega}$}, $\bar{\rho}^{(m)}(y_2+\eta_2)$  makes sense, and
\begin{align}
\label{20200830asdfa2114}
\bar{\rho}^{(m)}(y_2+\eta_2)-\bar{\rho}^{(m)}(y_2) =  \int_{0}^{\eta_2} \bar{\rho}^{(m+1)}(s)|_{s=y_2+z}\mm{d}z\quad \mbox{ for }\,\,0\leqslant  m\leqslant 4.
\end{align}
Moreover,
\begin{align}
&\label{esmmdforiasdfaasfanfty}
 \sup_{y\in \overline{\Omega}} \sup_{z\in \Psi}\left| \bar{\rho}^{(m+1)}(s)|_{s=y_2+z}\right| \lesssim 1,
\end{align}
where $\Psi:=[0,\eta_2]$ for $\eta_2\geqslant 0$,  {and $[\eta_2,0]$, otherwise}.

Making use of \eqref{20200830asdfa2114} and the relation
\begin{align}
\bar{\rho}(y_2+\eta_2)-\bar{\rho}(y_2) =\bar{\rho}'(y_2)\eta_2+{ \int_{0}^{\eta_2}\left(\eta_2(y,t) -
z\right)\bar{\rho}''(s)|_{s=y_2+z}\mm{d}z}, \nonumber
\end{align}
we find that
\begin{align}
\mathcal{G} & =g \int_{0}^{\eta_2}\left(\eta_2(y,t) -
z\right)\bar{\rho}''(s)|_{s=y_2+z}\mm{d}z ,\nonumber\\
 \partial_1 \mathcal{G}  &  =g  \big(\bar{\rho}'(s)|_{s=y_2+\eta_2}-\bar{\rho}'{(y_2)}\big)
\partial_1 \eta_2 ,\nonumber \\
\partial_2 \mathcal{G} & =g
\left(\int_{0}^{\eta_2} \partial_2\eta_2(y,t) \bar{\rho}''(s)|_{s=y_2+z}\mm{d}z + \int_{0}^{\eta_2}\left(\eta_2(y,t) -
z\right)\bar{\rho}'''(s)|_{s=y_2+z}\mm{d}z\right) .
\nonumber
\end{align}
Thanks to the estimates  \eqref{osdfadsa}, \eqref{20200830asdfa2114}, \eqref{esmmdforiasdfaasfanfty} and \eqref{fgessfdims}, we further deduce \eqref{2022011130957} and \eqref{2022201201821} from the above three identities.

(2) Recalling the definition of $\mathcal{R}$ in \eqref{20222042015205}, and then using \eqref{esmmdforiasdfaasfanfty}, \eqref{fgessfdims} and \eqref{202012241002}, we have
\begin{align}
 \|\mathcal{R}\|_{i,0}\lesssim \|\eta_2\|_{i,0}\|\eta_2\|_2,\  \|\nabla \mathcal{R}\|_{m,n}\lesssim \|\eta_2\|_{ {m},1+n}\|\eta_2\|_2\mbox{ and }   \|  \mathcal{R}_t\|_0\lesssim \|\eta_2\|_2\|u_2\|_0,\label{20222101181354}
\end{align}
where $0\leqslant i\leqslant 2$ and $0\leqslant m+n\leqslant 1$.
  Noting that
\begin{align}
&\Phi:= |\bar{\rho}'(y_2+\eta_2) |^2  -|\bar{\rho}' (y_2)|^2 =
(|\bar{\rho}'(y_2)|^2)'\eta_2+ \mathcal{R}, \label{2022100512045}
\end{align}
thus, exploiting \eqref{20222101181354}, we easily derive from the above identity that
\begin{align}
 \|\Phi\|_{i,0}\lesssim \|\eta_2\|_{ {i},0} , \ \|\nabla \Phi\|_{m,n}\lesssim \|\eta_2\|_{ {m},1+n}   \mbox{ and }  \| \Phi_t\|_0\lesssim \|u_2\|_{0} .\label{2022210181405}
\end{align}

Making use of  \eqref{osdfadsa}, \eqref{20222101181354}, \eqref{2022210181405},   \eqref{fgestims}  and  \eqref{fgessfdims}, we can obtain \eqref{2022dfdsfs2}--\eqref{2022201010520233}.  Similarly, we can also check that \eqref{202220927900}  holds.

(3) Finally, utilizing  \eqref{fgestims}, \eqref{fgessfdims} and \eqref{202012241002} again, we can easily deduce  \eqref{20asd10101643}--\eqref{202220110230924} from {the expressions  of $\mathcal{N}^\mu_{j,1} $, $\mathcal{N}^\mu_{j,2}$ and $\mathcal{N}^\mu $}.
\hfill $\Box$
\end{pf}
\subsection{Tangential estimates}\label{subsec:Horizon}
     This subsection is devoted to establishing the tangential estimates given in the following three lemmas,
     which include the boundedness of the  {horizontal} derivatives of $(\eta,u)$ and the temporal derivative of $u$.
\begin{lem}\label{lem:082sdaf41545}
For sufficiently small $\delta$, it holds that for $0\leqslant i\leqslant 2$,
\begin{align}
&
\frac{\mm{d}}{\mm{d}t}\left(\int\bar{\rho}\partial_1^{i}\eta\cdot\partial_1^{i} u\mm{d}y
+\frac{\mu}{2}\| \nabla \partial_1^{i} \eta\|_{0}^2\right) - E(\partial_1^{i}\eta) \nonumber \\
& \lesssim \| u\|_{i, 0}^2
+ \|\eta\|_{i+1,0}( \|  \eta\|_3\|  u\|_{i,1}+ \|  \eta\|_{i,1}\|  u\|_{\underline{1},2} ) \nonumber \\
&\quad +\delta_{i0} ( \|\eta\|_{1,1}\|  \eta\|_2\| u\|_1 +\| \eta\|_3 (\|  \eta_2 \|_1^2+ \|  \eta_2 \|_1\|\eta_2\|_3 )+\|\eta\|_{1,0}\|\eta\|_3\|\tilde{q}\|_0 )\nonumber \\
 &\quad + (1- \delta_{i0}) (\|\eta_2\|_{2,1}( \|  \eta_2 \|_{2,1}  \|  \eta  \|_3 +\|\eta_2\|_{ {1},2}\|\eta_2\|_3) +\|\eta_2\|_{1,2}^3 + \|\eta_2\|_{2,1}\|\eta_2\|_{1,2} \|\tilde{q}\|_{\underline{1},1}   ),
\label{202008241446}
\end{align}
where ${\tilde{q}:= q-\kappa \partial_2(|\bar{\rho}'|^2\eta_2)}$  {and $E (\cdot)$ is defined as \eqref{eeE}. }
\end{lem}
\begin{pf}
From now on, we define that
\begin{align*}
& \mathcal{N}^{\eta}_{j,k}:= \mathcal{N}^{\eta,1} _{j,k}+\mathcal{N}^{\eta,2} _{j,k}\,\,\  \mbox{ and }\,\,  \ \mathcal{N}^{\eta}:= \partial_k(\mathcal{N}^\eta _{1,k},\mathcal{N}^\eta _{2,k})^{\top},
\end{align*}
see  \eqref{202220926622063} and \eqref{2022209sadf26622063} for the definitions of   $\mathcal{N}^{\eta,1}_{j,k}$ and  $\mathcal{N}^{\eta,2}_{j,k}$, respectively.

Exploiting   the relation
  \begin{align}
\label{202201152005}
\partial_j( \mathcal{A}_{ij}f)= \mathcal{A}_{ij}\partial_jf,
\end{align} we have
\begin{align}
\mu \Delta_{\mathcal{A}} u=\mu \Delta u+ \mathcal{N}^\mu.
\label{20222040201510}
\end{align}

Using the  identities \eqref{2022100512045}, \eqref{202201152005}
and
\begin{align*}
 &\delta_{ik}\partial_k( |\bar{\rho}' |^2(\delta_{i2}(\delta_{j1}\partial_1\eta_2 + 2\delta_{j2} \partial_2\eta_2)+ \delta_{i1}\delta_{j2} \partial_1\eta_2) +(|\bar{\rho}' |^2)' \delta_{i2} \delta_{j2} \eta_2    )\\
 & =  \delta_{2j} \mm{div}(|\bar{\rho}' |^2 \nabla\eta_2 )+ \partial_j\partial_2(|\bar{\rho}'|^2\eta_2 ) , \end{align*}
 we can compute out that, for $j=1$ and $2$,
\begin{align}
\mathcal{K}^{\eta}_j =\, &   \kappa {\mathcal{A}}_{ik}\partial_k(|\bar{\rho}'(y_2+\eta_2)|^2 {\mathcal{A}}_{il} {\mathcal{A}}_{jm}\partial_l(y_2+\eta_2) \partial_m (y_2+\eta_2)  - |\bar{\rho}'|^2 {\mathcal{A}}_{i2} \mathcal{A}_{j2} )\nonumber \\
=\,&\kappa ( \delta_{2j} \mm{div}(|\bar{\rho}' |^2 \nabla\eta_2 )+ {\mathcal{A}}_{jk}\partial_k\partial_2(|\bar{\rho}'|^2\eta_2 )\nonumber \\
 & +   {\mathcal{A}}_{ik}\partial_k(|\bar{\rho}'(y_2+\eta_2)|^2 {\mathcal{A}}_{il}\partial_l\eta_2 {\mathcal{A}}_{jm} \partial_m\eta_2 +(|\bar{\rho}'(y_2+\eta_2) |^2 \nonumber \\
 & -|\bar{\rho}'  |^2 ) ({\mathcal{A}}_{i2}( {\mathcal{A}}_{j1}\partial_1\eta_2  + 2{\mathcal{A}}_{j2} \partial_2\eta_2 )+{\mathcal{A}}_{i1} {\mathcal{A}}_{j2} \partial_1\eta_2    ) \nonumber \\
 &+ |\bar{\rho}'  |^2  ({\mathcal{A}}_{i2}( \tilde{\mathcal{A}}_{j1}\partial_1\eta_2  + 2\tilde{\mathcal{A}}_{j2} \partial_2\eta_2 ) +\tilde{\mathcal{A}}_{i2}(\delta_{j1}\partial_1\eta_2 + 2\delta_{j2} \partial_2\eta_2 ) \nonumber \\
 &+({\mathcal{A}}_{i1} \tilde{\mathcal{A}}_{j2} + \tilde{\mathcal{A}}_{i1}\delta_{j2}) \partial_1\eta_2   )+ {\mathcal{A}}_{i2}  {\mathcal{A}}_{j2}\mathcal{R}+(|\bar{\rho}'  |^2)'( {\mathcal{A}}_{i2}  \tilde{\mathcal{A}}_{j2}\nonumber \\
 & + \tilde{\mathcal{A}}_{i2}  \delta_{j2} ) \eta_2 ) + \tilde{\mathcal{A}}_{ik}\partial_k  ( |\bar{\rho}'  |^2(\delta_{i2}(\delta_{j1}\partial_1\eta_2 + 2\delta_{j2} \partial_2\eta_2)+ \delta_{i1}\delta_{j2} \partial_1\eta_2) \nonumber
  \\
 &+(|\bar{\rho}' |^2)' \delta_{i2} \delta_{j2} \eta_2    )- \tilde{\mathcal{A}}_{jk}\partial_k \partial_2(|\bar{\rho}' |^2\eta_2 )) \nonumber\\
=\, &  \kappa ( \delta_{2j} \mm{div}(|\bar{\rho}' |^2 \nabla\eta_2 )+ {\mathcal{A}}_{jk}\partial_k\partial_2(|\bar{\rho}'|^2\eta_2 ))+ \partial_k \mathcal{N}^{\eta} _{j,k} ,  \label{202204201508} \end{align}
where $\mathcal{R}$  has been defined in \eqref{20222042015205}.

For any given $0\leqslant i\leqslant 2$,
applying $\partial_1^{i}$ to \eqref{01dsaf16asdfasf} and \eqref{20safd45},  and then using \eqref{20222040201510}  {and} \eqref{202204201508}, we get
\begin{equation}\label{01dsaf16asdfasf03n}
\begin{cases}
\partial_1^{i}\eta_t=\partial_1^{i}u ,\\[1mm]
\partial_1^{i}(\bar{\rho}u_t-\mu\Delta  u)\\[1mm]
\qquad =\partial_1^{i}((\kappa \mm{div}(|\bar{\rho}' |^2 \nabla\eta_2 ) +g\bar{\rho}'\eta_2
+ \mathcal{G})\mathbf{e}^2+\mathcal{N}^\mu+  \mathcal{N}^{\eta}-  \nabla_{\mathcal{A}} \tilde{q} ) ,\\[1mm]
[\partial_1^{i}, \mathcal{A}_{kl} ]\partial_l u_k+\mathcal{A}_{kl}\partial_1^{i}\partial_l u_k =0 , \\[1mm]
\partial_1^{i}(\eta_2, u_2,\partial_2(\eta_1, u_1))|_{\partial\Omega} =0.
\end{cases}
\end{equation}
Moreover, by the boundary condition \eqref{20safd45}, we have
\begin{align}
\partial_1^i\mathcal{N}^\mu _{1,2}|_{\partial\Omega}=\partial_1^i\mathcal{N}^{\eta} _{1,2}|_{\partial\Omega}=0.  \label{202220112002138}
\end{align}

Multiplying \eqref{01dsaf16asdfasf03n}$_2$ by $\partial_1^{i}\eta$, and then using  { integration}
by parts,  \eqref{01dsaf16asdfasf03n}$_1$, \eqref{01dsaf16asdfasf03n}$_4$ and the boundary condition of $\partial_1^i\mathcal{N}^{\eta}_{1,2} $ in \eqref{202220112002138}, we get
\begin{align}
&\frac{\mm{d}}{\mm{d}t}\left(\int\bar{\rho}\partial_1^{i}\eta\cdot\partial_1^{i} u\mm{d}y
+\frac{\mu }{2}\| \nabla \partial_1^{i} \eta\|_{0}^2\right)
-E(\partial_1^{i}\eta) =\|\sqrt{\bar{\rho}}\partial_1^{i}u\|_{0}^2+\sum_{j=1}^5I_{j,i},\label{202008241510}
\end{align}
where
\begin{align*}
&I_{1,i}:=\int\partial_1^{i} \mathcal{G}\partial_1^{i} \eta_2\mm{d}y,\ I_{2,i}:=-  \int\partial_1^{i} \mathcal{N}^\mu _{j,1}  \partial_1^{i+1}\eta_j\mm{d}y,\ I_{3,i}:=\int \partial_2\partial_1^{i} \mathcal{N}^\mu _{j,2} \partial_1^{i}\eta_j\mm{d}y\\ &
I_{4,i}:=-\int \partial_1^{i} \mathcal{N}^{\eta}_{j,k}  \partial_k\partial_1^{i}\eta_j\mm{d}y,  \ \  \mbox{ and } \  \ I_{5,i}:=- \int \partial_1^{i}\nabla_{\mathcal{A}}\tilde{q}\cdot\partial_1^{i} \eta\mm{d}y.
\end{align*}
Next we shall estimate the five integrals  $I_{1,i}$ to $I_{5,i}$ successively.

 (1) Using integration by parts, H\"older's inequality and \eqref{2022011130957}, we infer that
\begin{align}
&\label{202008241546}
I_{1,i} \lesssim
\begin{cases}
 \|\eta_2\|_{0}^2 \|\eta_2\|_{ {2},1}  &\mbox{for }i=0;\\
\|\eta_2\|_0\|\eta_2\|_{ {2},1} \|\eta_2\|_{1+i,0} &\mbox{for }i=1,\ 2.
\end{cases}
\end{align}

(2) Employing H\"older's inequality and \eqref{20asd10101643}, we deduce that
\begin{align}
I_{2,i}
\lesssim & \|\eta\|_{1+i,0}(\| \eta\|_{i,1} \|  u\|_{\underline{1},2}+\| \eta\|_3\|  u\|_{i,1 }  ) . \label{2022201161405} \end{align}

(3) By \eqref{20asd101asd01643} and   integration by parts, we also have
\begin{align}
 I_{3,i} \lesssim \|\eta\|_{1+i,0}(   \| \eta\|_{i,1} \|  u\|_{\underline{1},2}+\|  \eta\|_3\|  u\|_{i,1}) ,  \quad i=1,\ 2.   \label{20222011161447} \end{align}
In addition, we have (see (2.26) in \cite{JFJSZYYO})
\begin{align}
I_{3,0}
 \lesssim & \|\eta\|_{1,1}\|  \eta\|_2\| u\|_1 .   \label{2061032}
\end{align}

(4)  Noting that
$ {\mathcal{N}^{\eta,1}_{1,2}=0}$ by \eqref{202220926622063} with $(j,k)=(1,2)$, we have
$$I_{4,i}= -\int  (\partial_1^i\mathcal{N}^{\eta,1}_{1,1}\partial_1^{1+i} \eta_1+\partial_1^i\mathcal{N}^{\eta,1}_{2,k}\partial_k \partial_1^i\eta_2 +\partial_1^i\mathcal{N}^{\eta,2}_{j,k} \partial_k \partial_1^i\eta_j)\mm{d}y.
$$
Thus making use of  \eqref{omessetsim122n}, \eqref{202220927900}  and \eqref{2022dfdsfs2},  we derive from the above identity that
 {\begin{align}
I_{4,i}\lesssim
  \begin{cases}
 \| \eta\|_3( \|  \eta_2 \|_1^2+ \|\eta_2\|_1\|\eta_2\|_3) & \mbox{for}i=0;\\
(  \| \eta_2 \|_{  {i},1}\|   \eta  \|_3  +  \|\partial_2 \eta_1 \|_{
  i, 0}\|  \eta_2 \|_{ {1},2}) \|(\partial_1\eta_1,\nabla \eta_2)\|_{i,0}
 &\\
 \quad  +  \|  \eta\|_{i,1}  ( \|  \eta_2 \|_{ {1},2}^2 +\|  \eta_2 \|_{ {i},1}\|  \eta_2 \|_{3})\lesssim  \|\eta_2\|_{ {i},1}( \|  \eta_2 \|_{ 1 ,2}^2& \\
 \quad +\|\eta_2\|_{ {i-1},2}\|\eta_2\|_3+\|  \eta_2 \|_{ {i},1}  \|   \eta  \|_3 )+ \|\eta_2\|_{ {i-1},2} \|  \eta_2 \|_{ 1 ,2}^2 &\mbox{for }i=1,\ 2.
  \end{cases}  \label{20222210052105}
  \end{align}}

(5)
Finally, we control the integral term $I_{5,i}$. Recalling
$$\mm{div}_{\tilde{\mathcal{A}}}\eta=2(\partial_1\eta_1 \partial_2\eta_2-\partial_1\eta_2\partial_2 \eta_1),$$
 we integrate {$I_{5,i}$} by parts, and then make use of \eqref{improtian1}, \eqref{osdfadsa},    \eqref{01dsaf16asdfasf03n}$_4$  and
 \eqref{fgessfdims} to arrive at
\begin{align}
I_{5,0} =\int \mm{div}_{\tilde{\mathcal{A}}}\eta \tilde{q}\mm{d}y+ \int \mm{div}\eta \tilde{q}\mm{d}y
\lesssim   \| \eta\|_{1,0} \|\eta\|_3\|  \tilde{q}\|_0,    \label{201910040902n02}
\end{align}
and for $i=1$, $2$,
\begin{align}
I_{5,i}
=&\int\partial_1^{i+1} \eta \cdot\partial_1^{i-1}\nabla_{\tilde{\mathcal{A}}}\tilde{q}\mm{d}y
+\int\partial_1^{i}\mm{div} \eta\partial_1^{i}\tilde{q}\mm{d}y\nonumber \\
=& \int( \partial_1^{i+1} \eta_k \partial_1^{i-1} ({\tilde{\mathcal{A}}}_{k2}\partial_2 \tilde{q})+ \partial_1^{i+1} \eta_1 \partial_1^{i-1}({\tilde{\mathcal{A}}}_{11}\partial_1 \tilde{q} ))\mm{d}y
-\int\partial_1^{i}(\partial_1\eta_1\partial_2\eta_2)\partial_1^{i}\tilde{q}\mm{d}y\nonumber \\ &-\int\partial_1^{i+1} \eta_2[\partial_1^{i-1},\partial_2\eta_1]\partial_1\tilde{q} \mm{d}y
+\int[\partial_1^{i},\partial_2\eta_1]\partial_1\eta_2\partial_1^{i}\tilde{q}\mm{d}y\nonumber \\
\lesssim &\|\tilde{q}\|_{\underline{1},1}( \| \eta\|_{1+i,0}(\|\partial_1\eta\|_{\underline{i-1},1}+\|\eta_2
\|_{{i-1},2})+ \|\partial_1\eta\|_{\underline{1},1}\| \partial_1\eta_2\|_{\underline{i-1},1})\nonumber  \\
\lesssim &\|\tilde{q}\|_{\underline{1},1}( \| \eta_2\|_{i,1}\|\eta_2\|_{ {i-1},2}+\|\eta_2\|_{i,1} \| \eta_2\|_{1,2}). \label{202008241550}
\end{align}

Consequently, thanks  to the   estimates  \eqref{202008241546}--\eqref{202008241550} and \eqref{202012241002}, we obtain \eqref{202008241446} from \eqref{202008241510}. This completes the proof.
\hfill $\Box$
\end{pf}
\begin{lem}\label{lem:08241445}
For sufficiently small $\delta$, it holds that for $0\leqslant i\leqslant 2$,
\begin{align}
& \frac{\mm{d}}{\mm{d}t}\left(\|\sqrt{ \bar{\rho} }  u\|^2_{i,0}
-E(\partial_1^{i} \eta) +\frac{\kappa}{2}\|\bar{\rho}'  \partial_2\eta_1  \partial_1^{1+i}\eta_2\|^2_0  \right.\nonumber \\
  &\left.-\kappa\int |\bar{\rho}'  |^2   \partial_2\eta_1 \partial_1^{1+i}\eta_2\partial_2\partial_1^{i}\eta_2   \mm{d}y \right)+ c\|   u\|_{i,1}^2 \nonumber \\
&\lesssim
 \delta_{i0} \|  \eta_2 \|_1 \| \eta_2\|_3\|  u\|_1 +
 (1-\delta_{i0})( \|\eta_2\|_{ {2},1}( \|   \eta_2  \|_3 \|  u \|_{2,1}+
     \|   \eta_2  \|_{2,1} \|  u \|_{\underline{1},2} )
     \nonumber \\ & \quad+\|u\|_{2,1} (\|  \eta_2\|_{ {1},2}^2 +\|\eta \|_{3,0}\|u\|_{3} +\|  \eta\|_{2,1}\|u\|_{\underline{1},2})
\nonumber \\
 &\quad + \| \tilde{q}\|_{\underline{1},1}( \|  \eta\|_{2,1}\|  u\|_{2,1}+\|  \eta\|_{3,0 }(\| \eta_2 \|_{ {1},2} \|u\|_{\underline{1},2 }+ \|  \eta\|_3\| u\|_{ 2,1} ))\nonumber \\
 &\quad+\|  \eta\|_{3,0}\|u\|_{\underline{1},2 }( \|  \eta_2 \|_{2,1} + \| u_t\|_{1}
 + \|  \eta_2 \|_{1,2}\|\eta_2\|_3 + \|\eta\|_{2,1}\|u\|_{\underline{1},2}
 ) ).
  \label{202008241448}
\end{align}
\end{lem}
\begin{pf}
Multiplying \eqref{01dsaf16asdfasf03n}$_2$ by $\partial_1^{i}u$ in $L^2$, and then using  integration by parts, \eqref{01dsaf16asdfasf03n}$_1$,
 \eqref{01dsaf16asdfasf03n}$_4$ and  \eqref{202220112002138}, we have
\begin{align}
& \frac{1}{2}\frac{\mm{d}}{\mm{d}t}\left(\|\sqrt{\bar{\rho}} u\|_{i,0}^2
-E(\partial_1^{i}\eta)\right)+\mu \|  \nabla u\|_{i,0}^2 \nonumber\\
 &=\int\partial_1^{i} \mathcal{G}\partial_1^{i} u_2\mm{d}y-  \int \partial_1^{i} \mathcal{N}^\mu _{j,l}\partial_l \partial_1^{i}u_j\mm{d}y
 \nonumber \\
 &\quad -\int \partial_1^{i} \mathcal{N}^{\eta}_{j,k}  \partial_k\partial_1^{i} u_j \mm{d}y-\int \partial_1^{i}\nabla_{\mathcal{A}}\tilde{q } \cdot\partial_1^{i}u\mm{d}y=:\sum_{j=6}^9 I_{j,i}.\label{201510n}
\end{align}
Next we shall estimate the four integrals  {$I_{6,i}$ to $I_{9,i}$ in turn.}

(1)  Analogously to \eqref{202008241546} and \eqref{2022201161405}, we obtain
\begin{align}
& \label{202008241624}
I_{6,i} \lesssim \begin{cases}
\|\eta_2\|_{0}  \|\eta_2\|_{ {2},1} \|u_2\|_{0} &\mbox{for }\,\,i=0;\\
\|\eta_2\|_0\|\eta_2\|_{ {2},1}   \| u_2\|_{1+i,0} &\mbox{for }\,\,i=1,\ 2
\end{cases}
\end{align}
and
\begin{align}
I_{7,i}
 \lesssim & (1-\delta_{i0}) \|  \eta\|_{i,1}\|u\|_{\underline{1},2}\|u\|_{i,1} +
 \|\eta\|_3\|u\|_{i,1}^2.  \label{2022202181730}
\end{align}

(2)
We rewrite $I_{8,i} $ as follows
\begin{align}
I_{8,i}
=\,& \kappa\int|\bar{\rho}' |^2( \partial_1^{i}( \partial_2\eta_1  \partial_1\eta_2 )\partial_1^{1+i} u_1+
  \partial_1^{i} (   \partial_2\eta_1\partial_2\eta_2-\partial_2\eta_1\partial_1\eta_2 (2\partial_2\eta_1 \nonumber \\
  &+ \partial_1\eta_2(1+\partial_2\eta_2)+\partial_1\eta_1\partial_2\eta_1 )) \partial_1^{1+i} u_2
  +\partial_1^{i}(
\partial_2\eta_1\partial_1\eta_2 ( 2+3\partial_1\eta_1\nonumber \\
&+|\partial_1(\eta_1,\eta_2)|^2))\partial_2\partial_1^{i} u_2)
     \mm{d}y  -\int \partial_1^{i} \mathcal{N}^{\eta,2}_{j,k}  \partial_k\partial_1^{i} u_j \mm{d}y\nonumber \\
   =\,& I_{8,i,1} + I_{8,i,2},  \label{202209282108}
\end{align}
where we have defined that
\begin{align}
 I_{8,i,1}=\,  &\kappa\int\big\{|\bar{\rho}' |^2( \partial_1^{i}( \partial_2\eta_1  \partial_1\eta_2 )(\partial_1^{1+i}u_1+2\partial_2\partial_1^{i}u_2 )\nonumber \\
  &\quad \ \    +
  \partial_1^{i} (\partial_2\eta_1 \partial_2\eta_2  -2|\partial_2\eta_1|^2 \partial_1\eta_2  ) \partial_1^{ 1+i} u_2 )\big\}\mm{d}y\nonumber \end{align}
  and\begin{align}
 I_{8,i,2} :=\, &\kappa
    \int\big\{|\bar{\rho}'   |^2(  \partial_1^{i}(
\partial_2\eta_1\partial_1\eta_2 ( 3\partial_1\eta_1+|\partial_1(\eta_1,\eta_2)|^2))\partial_2\partial_1^{i} u_2\nonumber \\
 & \qquad -  \partial_1^{i} [   \partial_2\eta_1 |\partial_1\eta_2|^2 (1+\partial_2\eta_2) + \partial_1\eta_1 |\partial_2\eta_1|^2\partial_1\eta_2 ] \partial_1^{1+i} u_2)\big\}
     \mm{d}y  -\int \partial_1^{i} \mathcal{N}^{\eta,2}_{j,k}  \partial_k\partial_1^{i} u_j \mm{d}y.\nonumber
\end{align}

Exploiting \eqref{202222010109} and \eqref{01dsaf16asdfasf03n}$_1$, the integral $I_{8,i,1}$ can be further rewritten as follows
\begin{align}
 I_{8,i,1}=\,  &\kappa\frac{\mm{d}}{\mm{d}t}\int |\bar{\rho}' |^2(  \partial_2\eta_1 \partial_1^{1+i}\eta_2\partial_2\partial_1^{i}\eta_2 -|\partial_2\eta_1  \partial_1^{1+i}\eta_2|^2/2     )\mm{d}y
    +\tilde{I}_{8,i,1}, \label{202209291241} \end{align}
 where we have defined that
 \begin{align}
 \tilde{I}_{8,i,1}:=\kappa\int   &\big\{|\bar{\rho}'  |^2(  [\partial_1^{i},  \partial_2\eta_1]  \partial_1\eta_2 (\partial_1^{1+i}u_1+2\partial_2\partial_1^{i}u_2 ) \nonumber \\
  &+
(  [\partial_1^{i}, \partial_2\eta_1] \partial_2\eta_2-2[\partial_1^{i}, |\partial_2\eta_1 |^2 ]\partial_1\eta_2 ) \partial_1^{1+i} u_2 \nonumber \\ & -  \partial_2\eta_1  \partial_1^{1+i}\eta_2 (\partial_1^i(\tilde{\mathcal{A}}_{j2}\partial_2 u_j +\tilde{\mathcal{A}}_{11}\partial_1 u_1) -[\partial_1^i,\partial_2\eta_1 ]\partial_1 u_2)  \nonumber \\
&+   \partial_2\eta_1\partial_2u_1|\partial_1^{1+i}\eta_2|^2-\partial_2u_1 \partial_1^{1+i}\eta_2\partial_2\partial_1^{i}\eta_2     )\big\}\mm{d}y . \nonumber
 \end{align}
 {Making use of  \eqref{osdfadsa} and \eqref{2022dfdsfs2}, we easily obtain that}
 \begin{align}
 \tilde{I}_{8,i,1}+ I_{8,i,2}\lesssim  &
  \begin{cases} \|  \eta_2 \|_1 \| \eta_2\|_3 \|  u\|_1 &\mbox{for }\ i=0, \\
    \|\eta_2\|_{ {i},1}( \|   \eta_2  \|_3 \|  u \|_{i,1}+
     \|   \eta_2  \|_{i,1} \|  u \|_{\underline{1},2} )
      +\| \eta_2\|_{ 1,2}^2\|  u \|_{ {i},1} &\mbox{for }\ i=1,\ 2.
    \end{cases}    \label{20220929sadf1241}
 \end{align}

(3)   We integrate $\tilde{I}_{9,i}$ by parts, and then use \eqref{202201152005} and \eqref{01dsaf16asdfasf03n}$_3$ to get that
  \begin{align}
 I_{9,i}
=&\begin{cases}
0&\mbox{for }i=0;\\
  I_{9,i,1}+  I_{9,i,2}
&\mbox{for }i=1,\ 2,\end{cases}\label{20200n}
\end{align}
where we have defined that \begin{align*}
 & I_{9,i,1}:=-\int  [\partial_1^{i},\mathcal{A}_{kl}]\partial_l \tilde{q}\partial_1^{i}u_k \mm{d}y,\\
 & I_{9,i,2}:=-\int [\partial_1^{i}, \mathcal{A}_{kl} ]\partial_l u_k\partial_1^{i} \tilde{q}  \mm{d}y.
\end{align*}

It is easy to derive that for $i=1$, $2$,
 \begin{align}
  I_{9,i,1}
  \lesssim  \| \eta\|_{i,1}\|  u\|_{i,1}\|\tilde{q}\|_{\underline{1},1 }   \label{2022210601525}
\end{align}
and
\begin{align}
 I_{9,i,2}
 =\,& \int \big\{([\partial_1^{i},\partial_1 \eta_1  ]  \partial_1 u_1+[\partial_1^{i},\partial_2 \eta_1  ]\partial_1 u_2+[\partial_1^{i},\partial_1 \eta_2 ]\partial_2 u_1\nonumber \\
 &\quad -[\partial_1^{i}, \partial_2\eta_2]\partial_1 u_1+[\partial_1^{i},\partial_1 \eta_1  ]( \tilde{\mathcal{A}}_{kl}\partial_lu_k) )\partial_1^{i} \tilde{q} \big\} \mm{d}y  \nonumber \\
  \leqslant\,  &c( \|  \eta\|_{i,1}\|  u\|_{2,1}+\|  \eta\|_{1+i,0 }
  (\|  \eta_2\|_{ 1,2}\| u\|_{\underline{1}, 2}+\|  \eta\|_{ 3}\| u\|_{2,1}))\| \tilde{q}\|_{i,0}+ \tilde{I}_{9,i,2}, \label{20222010061540}
\end{align}
where we have defined that
$$
\tilde{I}_{9,i,2} :=\int  \partial_1^{1+i} \eta_2\partial_2 u_1 \partial_1^{i} \tilde{q}  \mm{d}y.
$$

By \eqref{01dsaf16asdfasf03n}$_2$, for $i=1$  and $2$,
\begin{align}
\partial_1^{i}\tilde{ q  }
= \partial_1^{i-1}(\mu\Delta  u_1  - \bar{\rho}\partial_tu_1+\mathcal{N}^\mu_1+ \mathcal{N}^{\eta}_1-   \tilde{\mathcal{A}}_{1l}\partial_l \tilde{q}).\nonumber
\end{align}
Making use of \eqref{osdfadsa}, \eqref{2022dsfs2}, \eqref{20222010101643}, \eqref{fgessfdims} and    the above identity, we have
\begin{align}
\tilde{I}_{9,i,2} =\,&\int   \partial_1^{1+i}  \eta_2\partial_2 u_1
 \partial_1^{i-1}(\mu\Delta  u_1  - \bar{\rho}\partial_tu_1+\mathcal{N}^\mu_1+  \mathcal{N}^{\eta}_1-   \tilde{\mathcal{A}}_{1l}\partial_l \tilde{q}) \mm{d}y\nonumber \\
 \lesssim \, & \|  \eta\|_{1+i,0}\|u\|_{\underline{1},2 }(\| u\|_{1+i,0}+\| u_t\|_{i-1}
 + \|\eta\|_{i,1}\|u\|_{\underline{1},2} + \|\eta\|_{\underline{1},2}\|u\|_{i,1} \nonumber \\
 &+ \|  \eta_2 \|_{i-1,2}(\|\eta_2\|_3+ \|\tilde{q}\|_{\underline{1},1} )+ \|  \eta_2 \|_{i,1}\|\eta\|_3
 )+\left|\int  \partial_1^{1+i}  \eta_2  \partial_2 u_1
\partial_2^2 \partial_1^{i-1}   u_1 \mm{d}y\right|. \nonumber
\end{align}
In addition, using   integration by parts and \eqref{fgessfdims}, we can {deduce} that
\begin{align}
{\Big|}\int  \partial_1^{1+i}  \eta_2  \partial_2 u_1
\partial_2^2 \partial_1^{i-1}   u_1 \mm{d}y{\Big|}= & \frac{1}{2}{\Big|}\int  (\partial_2 \partial_1^{ i} \eta_2\partial_1^{ i}  (\partial_2u_1)^2-\delta_{i2} \partial_1^{1+i} \eta_2 \partial_2\partial_1^{i-1}u_1
\partial_2^2  u_1) \mm{d}y{\Big|}\nonumber \\
\lesssim&  \|\eta_2\|_{i,1}\|u\|_{i,1}\|u\|_{\underline{1},2} +\|\eta_2\|_{1+i,0}\|u\|_{i,1}\|u\|_{3}. \nonumber
\end{align}

Now putting \eqref{2022210601525} and \eqref{20222010061540} into \eqref{20200n}, and then using the above two estimates, we get
\begin{align}
 I_{9,i} \lesssim \,  &(1-\delta_{i0})(\| \eta\|_{i,1}\|  u\|_{i,1}\|\tilde{q}\|_{\underline{1},1 }   +( \|  \eta\|_{i,1}\|  u\|_{2,1}+ \|  \eta\|_{1+i,0 } (\|  \eta_2\|_{ 1,2}\| u\|_{\underline{1}, 2}\nonumber \\
 &+\|  \eta\|_{ 3}\| u\|_{2,1}))\| \tilde{q}\|_{i,0}+\|  \eta\|_{1+i,0}\|u\|_{\underline{1},2 }(\| u\|_{1+i,0}+\| u_t\|_{i-1}
 + \|\eta\|_{i,1}\|u\|_{\underline{1},2} \nonumber \\
 & + \|\eta\|_{\underline{1},2}\|u\|_{i,1} + \|  \eta_2 \|_{i-1,2}(\|\eta_2\|_3+\|\tilde{q}\|_{\underline{1},1})+\|  \eta_2 \|_{i,1}\|\eta\|_3  ) \nonumber \\
 &+\|\eta_2\|_{i,1}\|u\|_{i,1}\|u\|_{\underline{1},2} +\|\eta_2\|_{1+i,0}\|u\|_{i,1}\|u\|_{3} ) \label{202220041} .
 \end{align}

Therefore, substituting  \eqref{202209291241}  into \eqref{201510n}, and then making use of \eqref{aprpiosesnew},
 \eqref{2022111522311}, \eqref{202008241624}, \eqref{2022202181730}, \eqref{20220929sadf1241}, \eqref{202220041} and \eqref{202012241002}, we conclude \eqref{202008241448}
for sufficiently small $\delta$. This completes the proof. \hfill $\Box$
\end{pf}
\begin{lem}\label{2019100216355nnn}
For sufficiently small $\delta$, we have
\begin{align}
&\frac{\mm{d}}{\mm{d}t}\left(\|\nabla_{\mathcal{A}} u \|^2_0+ \kappa\int|\bar{\rho}' |^2 \nabla\eta_2 \cdot\nabla u_2  \mm{d}y \right)+
c\|   u_t\|_0^2\nonumber \\
&\lesssim \|(\eta_2,  u_2 )\|_{1}^2+\|u\|_2^3+\|u\|_2^2\|\tilde{q}\|_1+ {\|\eta_2\|_1\|u_t\|_1}
\label{202221235}
\end{align}
and
\begin{align}
&\frac{\mm{d}}{\mm{d}t}\left(\|\sqrt{\bar{\rho}}\psi \|^2_0- E(  u) \right)+
c\|   u_t\|_1^2\nonumber \\
&\lesssim(\|  \eta_2\|_2 ^2+\|u\|_2^2)\|u\|_2^2 + \|  \eta \|_{3}^2(\|u\|_{1,0}^2+\|\eta_2\|_{2,1}^2\|u\|_1^2)\nonumber \\
 &\quad + (\|u\|_{1,0}+\|\eta_2\|_{2,1}\|u\|_1)\|u\|_2^2 +\|u\|_2^4,
\label{20222saf201121235}
\end{align}
where $\psi :=u_t- {u} \cdot \nabla_{\mathcal{A}} u $.
\end{lem}
\begin{pf}
(1)
We rewrite \eqref{01dsaf16asdfasf}$_2$ {as follows,}
  \begin{align}\label{01sdsdffasf}\bar{\rho}u_t+\nabla_{\mathcal{A}} \tilde{q} -\mu \Delta_{\mathcal{\mathcal{A}}} u=(\kappa \mm{div}(|\bar{\rho}' |^2 \nabla\eta_2 ) +g\bar{\rho}'\eta_2
+ \mathcal{G})\mathbf{e}^2 + \mathcal{N}^{\eta} .\end{align}
In view of \eqref{01dsaf16asdfasf}$_3$, we find that
$$\div_{\mathcal{A}} u_t=-\div_{\mathcal{A}_t} u. $$  Multiplying \eqref{01sdsdffasf} by $u_t$ in $L^2$, and then utilizing the above relation and \eqref{20safd45}, we infer that
  \begin{align}\label{01sdfasf}
\frac{\mm{d}}{\mm{d}t}\left(\frac{\mu}{2} \|\nabla_{\mathcal{\mathcal{A}}} u\|_0^2
+\kappa\int|\bar{\rho}' |^2 \nabla\eta_2 \cdot\nabla u_2 \mm{d}y
\right)+\|\sqrt{\bar{\rho}}u_t\|_0^2  =J_1 +J_2 ,
\end{align}
where
\begin{align*}
 J_1:=\, &\int ( g\bar{\rho}' \eta_2 \partial_t u_2+\mathcal{G}  \partial_t u_2+\mu \nabla_{\mathcal{\mathcal{A}}} u:\nabla_{\mathcal{\mathcal{A}}_t} u - \mm{div}_{\mathcal{A}_t  }u  \tilde{q})\mm{d}y,\\
    J_2:=\, &\int(\kappa  |\bar{\rho}'  \nabla u_2|^2    -\mathcal{N}^{\eta}_{j,k}  \partial_k  \partial_t u_j)\mm{d}y.
   \end{align*}

Thanks to H\"older's inequality, the embedding inequality of $L^3\hookrightarrow H^1$, \eqref{2022011130957}, \eqref{202220927900} and \eqref{2022dfdsfs2}, we have
\begin{align*}
&J_1+J_2 \lesssim (\|\eta_2\|_0+ \| \eta_2\|_1 )\| u_t\|_0 +\|u_2\|_1^2+ \|u\|_2^3+\|u\|_2^2\|\tilde{q}\|_1 + {\|\eta_2\|_1\|u_t\|_1}.
   \end{align*}
Putting the  {above estimate} into \eqref{01sdfasf}, and using Young's inequality, we get \eqref{202221235}.

(2) Let
\begin{align*}
   & J_3:=-\int
 \partial_t \mathcal{N}^{\eta}_{j,k}  \partial_k\psi_j \mm{d}y,\nonumber \\
&J_4:=\int
(( \mu \mathcal{A}_{il} \partial_l( \mathcal{A}_{ik} \partial_k  u )  +(\mathcal{G}+g  \bar{\rho}'\eta_2+ \kappa \mm{div}(|\bar{\rho}' |^2 \nabla\eta_2 ) ) \mathbf{e}^2+
 \mathcal{N}^{\eta}
 \\
&\qquad-  \bar{\rho} u\cdot\nabla_{\mathcal{A}}    u   -\bar{\rho}\psi) \cdot ( u\cdot\nabla_{\mathcal{A}} \psi)-   \partial_t (  \bar{\rho} u\cdot\nabla_{\mathcal{A}}    u ) \cdot \psi ) \mm{d}y,
 \\
&J_5:=   \int( ( g\bar{\rho}'(s)|_{s=y_2+\eta_2}u_2 +\kappa  \mm{div}(|\bar{\rho}' |^2 \nabla u_2 ))  \psi_2-   \mu
  \partial_t ( \mathcal{A}_{il} \mathcal{A}_{ik} \partial_k  u ) \cdot \partial_l \psi )\mm{d}y .
   \end{align*}
Thus we can get the following identity  from \eqref{01dsaf16asdfasf}   (referring to (2.38) in \cite{JFJSZYYO})
\begin{align}
 \frac{1}{2}\frac{\mm{d}}{\mm{d}t}\|\sqrt{\bar{\rho}}\psi \|^2_0
= J_{3}+J_{4}+J_{5} .  \label{62053}
  \end{align}

Note that the integral term $J_5$ can be further rewritten {as follows:}
\begin{align}
 J_5=\frac{1}{2}\frac{\mm{d}}{\mm{d}t}  E(  u)-
 \mu \|  \nabla u_t\|_0^2+\tilde{J}_5,  \label{2022101262053}
  \end{align}
  where
 \begin{align*}
 \tilde{J}_5:=\, & \int\big\{ ( g(\bar{\rho}'(y_2+\eta_2) - \bar{\rho}'(y_2))u_2\partial_t u_2-   ( g \bar{\rho}'(y_2+\eta_2)u_2 + \kappa  \mm{div}(|\bar{\rho}' |^2 \nabla u_2 )) u  \cdot \nabla_{\mathcal{A}} u_2  \\
 &\quad + \mu ( \partial_t ( \mathcal{A}_{il} \mathcal{A}_{ik} \partial_k  u) \cdot \partial_l(u \cdot \nabla_{\mathcal{A}} u)-\partial_t(\mathcal{A}_{il}\tilde{\mathcal{A}}_{ik} \partial_k  u +  \tilde{\mathcal{A}}_{il}\partial_i  u  )\cdot \partial_l \partial_t u ) )\big\}\mm{d}y.
  \end{align*}
We make use of \eqref{01dsaf16asdfasf}$_1$,  \eqref{2022201201821},  \eqref{2022dsfs2}, \eqref{2022201010520233} and \eqref{fgestims}
to deduce that
\begin{align*}
J_3 +
J_{4}+\tilde{J}_{5}  \lesssim &\,( \| \eta_2\|_2+\|u\|_2 +\|u_t\|_1)\|u\|_2 \|u_t\|_1 + \|  \eta \|_{3} \| u_2 \|_{1}   \|u_t\|_1
\\
&+ \|u_2\|_1\|u\|_2^2 +(\|\eta_2\|_2+\|u\|_2)\|u\|_2^3+\| \eta\|_3\|u_t\|_1^2.
\end{align*}
 Inserting \eqref{2022101262053} into \eqref{62053}, and using the above estimate, \eqref{2022111522311},  \eqref{2022220111522311}, \eqref{omesssafden}   and Young's inequality, we obtain \eqref{20222saf201121235} for sufficiently small $\delta$.
This completes the proof. \hfill $\Box$
\end{pf}

\subsection{Stabilizing estimates }\label{subseasdc:Vor}
This subsection is devoted to establishing the stabilizing estimates.
    \begin{lem}\label{lem:08250749}
We have
\begin{align}
\label{202008250745}
\|\eta_2\|_{i,1}^2
&\lesssim-E(\partial_1^{i}\eta)  \quad\mbox{ for }\ \,\,0\leqslant i\leqslant 2,
\\
\label{20250745}
\| u_2\|_{1}^2
&\lesssim-E(u)+\|  \eta\|_3\| u\|_1^2.
\end{align}
\end{lem}
\begin{pf}By virtue of the definition of $\kappa_{\mm{C}}$, it is easy to see that
\begin{align*}
-g\int \bar{\rho}'w_2^2\mm{d}y\geqslant- \kappa_{\mm{C}}   \|\bar{\rho}'\nabla  w_2 \|^2_0\quad \mbox{ for any }\,\,w\in H_{\sigma}^1,
\end{align*}
 which, together with the stability conditions \eqref{2022205071434} and $\kappa>\kappa_{\mm{C}}$, implies that
\begin{align}\label{202008250814}
\| w_2\|_{1}^2\lesssim  (\kappa-\kappa_{\mm{C}}){\|\bar{\rho}'\nabla w_2\|_0^2}\leqslant -E(w)\quad\mbox{ for any }\,\, w\in H_{\sigma}^1.
\end{align}

 In order to apply the above estimate to $\eta$, we have to modify $\eta$ because of $\mm{div}\eta\not =0$.
 Let $0\leqslant i\leqslant 2$ be given. Then there is a Bogovskii's operator
 $\mathcal{B} : \partial_1^{i}\mm{div}\eta\in \underline{L}^2 \to H^1_0$, such that ({see \cite[Lemma A.5]{JFJSZWC}})
\begin{align}
 \mathrm{div}\mathcal{B}(\partial_1^{i}\mm{div} \eta)=\partial_1^{i}\mm{div}\eta\  \ \mbox{ and }\ \
 \|{\mathcal{B}}(\partial_1^{i}\mm{div}\eta)\|_1  \lesssim \|\partial_1^{i}\mm{div}\eta\|_0.
 \label{2022202191518}
\end{align}

Now, we use $\partial_1^i\eta -\mathcal{B}(\partial_1^{i}\mm{div}\eta)$ to rewrite ${E}(\partial_1^{i}\eta)$ as
\begin{align}\label{201912191626}
 {E}(\partial_1^{i}\eta) = {E}( \partial_1^{i} \eta-\mathcal{B}(\partial_1^{i}\mm{div} \eta))
 -{E}(\mathcal{B}(\partial_1^{i}\mm{div} \eta) ) -J_{6,i},
\end{align}
where
\begin{align*}
J_{6,i}: = 2 \kappa\int |\bar{\rho}'|^2\nabla \partial_1^{i}{\eta_2}\cdot   \nabla \mathcal{B}_2(\partial_1^{i}\mm{div} \eta)\mm{d}y
-2g\int \bar{\rho}'\partial_1^{i}\eta_2{\mathcal{B}_2}(\partial_1^{i}\mm{div} \eta)\mm{d}y.
\end{align*}
Recalling $\partial_1^{i} \eta- \mathcal{B}(\partial_1^{i}\mm{div} \eta) \in H_{\sigma}^1$, we use \eqref{202008250814} to get
\begin{align}\nonumber
\| \partial_1^{i}\eta_2- \mathcal{B}_2(\partial_1^{i}\mm{div} \eta)\|_{ 1}^2 \lesssim-E(\partial_1^{i} \eta-\mathcal{B}(\partial_1^{i}\mm{div} \eta)),
\end{align}
which, together with \eqref{201912191626} and Young's inequality, gives
\begin{align}\label{202008250814nn}
\| \eta_2\|_{i,1}^2 \lesssim \| \mathcal{B}_2(\partial_1^{i}\mm{div} \eta)\|_{ 1}^2-{E}(\mathcal{B}(\partial_1^{i}\mm{div} \eta) ) -E(\partial_1^{i}\eta) -J_{6,i}.
\end{align}

Using the estimate in \eqref{2022202191518}, we find that
\begin{align}
  &-{E}(\mathcal{B}(\partial_1^{i}\mm{div} \eta) ) - J_{6,i} +\|\mathcal{B}_2(\partial_1^{i}\mm{div} \eta)\|_{1}^2  \nonumber \\
\lesssim\, &\|\eta_2\|_{i,1}\|\mathcal{B}_2(\partial_1^{i}\mm{div} \eta)\|_{ 1 } +\|\mathcal{B}_2(\partial_1^{i}\mm{div} \eta)\|_{1}^2\nonumber \\
 \lesssim\, &  \|\eta_2\|_{i,1} \|\mm{div}\eta\|_{ i,0} +\|\mm{div}\eta\|_{i, 0}^2.\nonumber
\end{align}
Finally, putting the above estimate into \eqref{202008250814nn}, and making use of \eqref{improtian1}, \eqref{osdfadsa} and Young's inequality, we obtain \eqref{202008250745}.

 On the other hand, we can argue, similarly to the derivation of \eqref{202008250745} with $i=0$,
to establish \eqref{20250745} by employing the relation
\begin{align}
\label{2022202121921}
\mm{div}u=-\mm{div}_{\tilde{\mathcal{A}}}u
\end{align}
and the estimate
$$  \|\mm{div}_{\tilde{\mathcal{A}}}u\|_0\lesssim \|  \eta\|_3\|u\|_1. $$
This completes the proof.
\hfill $\Box$
\end{pf}

\subsection{Normal estimates of $( {u},q)$  and  equivalence estimates}\label{subsec:Vor}
In this subsection, we derive the normal estimates of $( {u},q)$ by the regularity theory of Stokes equations.
 \begin{lem}
For sufficiently small $\delta$, we have
\begin{align}
\label{omessetsim}
&\|  u\|_{\underline{i},2+j} +\|\tilde{q}\|_{\underline{i},1+j} \lesssim \|   \eta_2\|_{\underline{i},2+j}+\|u_t \|_{\underline{i},j}
\  \ \mbox{ for }\  \ 0\leqslant i+j\leqslant  1.
\end{align}
\end{lem}
\begin{pf}By \eqref{202201152005}, \eqref{01dsaf16asdfasf03n}$_2$ and \eqref{2022202121921},
we see that $u$ satisfies the following Stokes problem: for $i=0,1$,
\begin{equation}
\label{Stokesequson}
\begin{cases}
\partial_1^i ( \nabla \tilde{q}-\mu\Delta  u)
 =   \partial_1^i((\kappa \mm{div}(|\bar{\rho}' |^2 \nabla\eta_2  )    +g\bar{\rho}'\eta_2
+ \mathcal{G})\mathbf{e}^2- \bar{\rho}u_t+\mathcal{N}^\mu+\mathcal{N}^{\eta}-  \nabla_{\tilde{\mathcal{A}}} \tilde{q} ) ,\\
 \partial_1^i \div u=-  \partial_1^i\div({\tilde{\mathcal{A}}}^{\top} u) ,\\
 \partial_1^i (u_2,
 \partial_2   u_1) |_{\partial\Omega} =0.
\end{cases}
\end{equation}
We apply the regularity estimate \eqref{202201122130} to the above Stokes problem to get
\begin{equation}   \label{ometsim}
  \|  u \|_{i,2+j} +\| \tilde{q}\|_{i,1+j} \lesssim \|    \mm{div}(|\bar{\rho}' |^2 \nabla\eta_2 )    ,\eta_2,u_t)\|_{i,j} +\|  \mathcal{N}^{\eta}  \|_{i,j}+J_7, \end{equation}
where $0\leqslant i+j\leqslant 1$ and
$$J_7:=\|  ( \mathcal{G}, \mathcal{N}^\mu, \nabla_{\tilde{\mathcal{A}}} \tilde{q},{\tilde{\mathcal{A}}}^{\top} u)  \|_{i,j}  +\| \div({\tilde{\mathcal{A}}}^{\top} u)\|_{i,1+j}.$$

By \eqref{2022dsfs2}--{\eqref{3.21}}, we have
\begin{align}
\|  \mathcal{N}^{\eta}  \|_{i,j}\lesssim \|  \eta_2 \|_{i,2+j}\| \eta\|_3.
\end{align}
Thanks to \eqref{2022011130957}, \eqref{2022201201821},  \eqref{20222010101643} and \eqref{202220110230924}, it is easy to  see that
\begin{align}
 &J_7\lesssim \| \eta_2\|_0 +\|\eta\|_3\| \tilde{q}\|_{\underline{i},1+j}+  \| \eta\|_3\| u\|_{\underline{i},2+j} .\nonumber
 \end{align}
Inserting the above two inequalities into \eqref{ometsim} yields the desired estimate \eqref{omessetsim}. This completes the proof.
  \hfill $\Box$
\end{pf}

 \begin{lem}
For sufficiently small $\delta$, we have
\begin{align}
&\label{202012252005}\sqrt{\mathcal{E}}\,\mbox{ and }\, \| (\nabla \eta,u)\|_2\, \mbox{ are equivalent to each other} ,
\end{align}
where the  equivalent coefficients in \eqref{202012252005} are independent of $\delta$.
\end{lem}
\begin{pf} To show \eqref{202012252005}, recalling the definition  of $\mathcal{E}$, we use \eqref{202211} and the Stokes estimate
\eqref{omessetsim} with $i=j=0$ to conclude
\begin{align}
& \|(\nabla \eta,u)\|_2^2 \lesssim \mathcal{E}\lesssim \|(\nabla \eta,u)\|_2^2+\|u_t\|_0^2.
\end{align}
 Obviously, to complete the proof, it suffices to show
 \begin{align}
\|u_t\|_0  \lesssim   \| \nabla  \eta \|_2+\| u\|_2 . \label{20222012120226}
\end{align}
 To  this end, we multiply \eqref{01dsaf16asdfasf}$_2$ by $u_t$ in $L^2$ to obtain
\begin{align}
\|\sqrt{ \bar{\rho} }u_{t}\|_{0}^2  =  & \int   (((\kappa \mm{div}(|\bar{\rho}' |^2 \nabla\eta_2  )    + g  \bar{\rho}'\eta_2+\mathcal{G})\mathbf{e}^2
+\mu\Delta_{\mathcal{A}}  u)\cdot u_t + \nabla \tilde{q} \cdot({\mathcal{A}}_t^{\top}u )) \mathrm{d}y \nonumber\\
&+ \int    \mathcal{N}^{\eta} \cdot u_t \mathrm{d}y =: J_8+J_9. \label{20222012120232}
\end{align}

Utilizing \eqref{2022dsfs2}, we have
\begin{align*}
 J_9   \lesssim &  \| \mathcal{N}^{\eta}\|_0\|u_t\|_0 \lesssim \|  \eta_2 \|_2\|  \eta \|_3 \|u_t\|_0. \end{align*}
In addition, it holds
\begin{align*}
 J_8   \lesssim & ( \|( \eta_2,u)\|_2+ \| u\|_2^2 )\|u_t \|_0 + {\|\eta_2\|_{2}} \| u\|_2^2 . \end{align*}
    Putting the above two estimates into \eqref{20222012120232}, and using
 \eqref{aprpiosesnew} and Young's inequality, we arrive at \eqref{20222012120226}.
 This completes the proof. \hfill$\Box$
\end{pf}

\subsection{Curl estimates for $\eta$}\label{subec:Vor}
    In this subsection we establish the curl estimates of $\eta$ for the normal estimates of $\eta$.
\begin{lem}\label{2055nnn}
We have
\begin{align}
&\label{202005021600}
 \frac{\mm{d}}{\mm{d}t}  \left(\frac{\mu }{2}
 \|\nabla_{\mathcal{A}} \partial_1^i \mm{curl} \eta \|_0^2+
 \int  \partial_1^i \mm{curl}   \eta
  \partial_1^i \mm{curl}_{\mathcal{A}} \left(\bar{\rho}u\right)
 \mm{d}y\right)
 + c\|  \nabla^2 \eta_2 \|^2_{i,0}\nonumber\\
&\lesssim
 \|\eta_2\|_{i,1}(\|\eta_2\|_1+\|\eta_2\|_{1+i,1}) +\|u\|_{ i,1}\|u\|_{\underline{i},1} + \| \eta\|_{2i,0}\| \eta_2\|_3^2 \nonumber \\
 &\quad+ \|  \eta\|_3(\|\eta_2\|_{i,2}^2+ \|  \eta\|_{2i,1}  \| \eta_2\|_{2,1}  +\delta_{i 0}\|\eta\|_{1,2}\|  u\|_{ 2}+\|u\|_2^2)\nonumber\\
 &\quad+ \delta_{i 1}(\|\eta\|_{1,2}^2\|  u\|_{\underline{1},2}
 + \|\eta\|_{1,2}\|\eta\|_3\|  u\|_{2,1} )
\quad \mbox{ for }\,\,\, i=0,   1
\end{align}
 and
    \begin{align}
 &\frac{\mm{d}}{\mm{d}t} \left(
 \frac{\mu }{2}\|   \nabla \partial_2 \mm{curl} \eta\|^2_0+\mathcal{I}_1\right)
+c \|  \partial_2^2 \eta  \|_{1,0}^2 \nonumber \\
&
 \lesssim  \| \eta_2\|_{2}^2 + \| \eta_2   \|_{2,1}^ 2  +\|u\|_1\|  u\|_3+ (
\|\partial_2\eta_2\|_{ {1},0}^{1/2}\|\partial_2^2\eta_2\|_{0}^{1/2}\nonumber \\
 &\quad
+\|\partial_2\eta_2\|_{\underline{1},0}
)(
\|\partial_2^2\eta_1\|_{ {1},0}^{1/2}\|\partial_2^3\eta_1\|_{0}^{1/2}
+\|\partial_2^2\eta_1\|_{\underline{1},0}
) + \|\eta_2\|_{2,1}\|\eta\|_3^2 +\sqrt{\mathcal{E}}\mathcal{D}, \label{202201211416}  \end{align}
           where we have defined that
\begin{align*}
\mathcal{I}_1:=\, &\int ( \mu   \partial_2^2  ( \partial_1\eta_1\partial_2 \eta_1) \partial_2^2\mm{curl} \eta-\partial_2^2\mm{curl}  \eta
 \mm{curl}_{\mathcal{A}}(\bar{\rho}  u))  \mm{d}y +\mathcal{I}_2,\\
\mathcal{I}_2:=\,& \int   \partial_2^2  \mm{curl} \eta (   \partial_2  (\partial_2\eta_1\mm{curl}_{\mathcal{A}}\partial_1 \eta  )
+\partial_2\eta_1 (  \partial_2\mm{curl}_{\mathcal{A}} \partial_1 \eta-\partial_2\eta_1 \mm{curl}_{\mathcal{A}} \partial_1^2 \eta) ) \mm{d}y.
 \end{align*}

\end{lem}
\begin{pf} We can directly compute out that
$$\mm{curl}_{\mathcal{A}}\left({G^{\eta}} \mathbf{e}^2\right)=g\mm{curl}_{\mathcal{A}}\left(- \bar{\rho} \mathbf{e}^2\right)
=-g \mathcal{A}_{1j}\partial_j \bar{\rho} =g  \bar{\rho}'\partial_1\eta_2 .
$$
Therefore, applying $\mm{curl}_{\mathcal{A}}$ to \eqref{01dsaf16asdfasf}$_2$  and keeping in mind that
$$\mm{curl}_{\mathcal{A}} \Delta_{\mathcal{A}}u=\Delta_{\mathcal{A}}\mm{curl}_{\mathcal{A}} u , $$
we get
\begin{align}
&\partial_t\mm{curl}_{\mathcal{A}}( \bar{\rho} u) -\mu \Delta_{\mathcal{A}} \mm{curl}_{\mathcal{A}}u-\kappa   \mathcal{A}_{1j} \partial_j  \mm{div}(|\bar{\rho}' |^2 \nabla\eta_2) \nonumber \\
&   = g\bar{\rho}'\partial_1 \eta_2
+  \mm{curl}_{\mathcal{A}_t}( \bar{\rho} u)+ \mm{curl}_{\mathcal{A}}\mathcal{N}^{\eta} .\label{201910072117}
\end{align}

Let $\alpha$ be a multi-index. An application of $\partial^{\alpha}$ to \eqref{201910072117} yields
\begin{align}
& \partial_t\partial^{\alpha}\mm{curl}_{\mathcal{A}}( \bar{\rho} u)
-\mu \partial^{\alpha}\Delta_{\mathcal{A}}\mm{curl}_{\mathcal{A}} u- \kappa\partial^{\alpha}     \mm{div}(|\bar{\rho}' |^2 \nabla \partial_1\eta_2 )   \nonumber \\
&  = g\partial^{\alpha}( \bar{\rho}'\partial_1\eta_2)   + \partial^{\alpha}\mm{curl}_{\mathcal{A}_t}( \bar{\rho} u)+ \kappa\partial^{\alpha}(
 \tilde{\mathcal{A}}_{1j} \partial_j  \mm{div}(|\bar{\rho}' |^2 \nabla\eta_2    )+\mm{curl}_{\mathcal{A}}\mathcal{N}^{\eta}).\label{202005021542}
\end{align}

(1) If we take $\alpha=(i,0)$ ($0\leqslant i\leqslant 1$), and then multiply \eqref{202005021542} by $ {\partial_1^i\mm{curl}\eta}$
in $L^2$,   we  can obtain
 \begin{align}
&\frac{\mm{d}}{\mm{d}t} \int\partial_1^i \mm{curl}_{\mathcal{A}} (\bar{\rho} u)\partial_1^i\mm{curl}  \eta  \mm{d}y =\sum_{j=1}^5K_{j,i} , \label{f202005031720}
\end{align}
where
\begin{align}
K_{1,i} :=& \kappa \int  \partial_1^i( \tilde{\mathcal{A}}_{1j}
 \partial_j  \mm{div}(|\bar{\rho}' |^2 \nabla\eta_2 )  +\mm{curl}_{\mathcal{A}}\mathcal{N}^{\eta}    ) \partial_1^i \mm{curl}  \eta  \mm{d}y,
\nonumber\\
K_{2,i}:=&  \kappa\int \partial_1^i  ( \mm{div}(|\bar{\rho}' |^2 \nabla \partial_1\eta_2 )) \partial_1^i\mm{curl} \eta\mm{d}y,\nonumber \\
K_{3,i}:=& \mu  \int \partial^{i}_1\Delta_{\mathcal{A}}\mm{curl}_{\mathcal{A}} u   \partial_1^i\mm{curl} \eta \mm{d}y,\nonumber \\
 K_{4,i}:=&
\int  \partial_1^i \mm{curl}_{\mathcal{A}_t}( \bar{\rho} u) \partial_1^i\mm{curl} \eta \mm{d}y, \nonumber \\
 K_{5,i}:=&\int \left(\partial_1^i
  \mm{curl}_{\mathcal{A}} (\bar{\rho} u )\partial_1^i\mm{curl}  u  - g  \bar{\rho}'\partial_1^{i}\eta_2
   \partial_1^{1+i} \mm{curl} \eta\right)\mm{d}y .\nonumber
\end{align}
Next we shall estimate the above five integrals in sequence.

Thanks to  \eqref{osdfadsa}     and \eqref{2022dsfasdas2}--\eqref{2022210220847}, we can estimate that\begin{align}
K_{1,i}=& \kappa (-1)^{ i}  \int
( \tilde{\mathcal{A}}_{1j}
 \partial_j  \mm{div}(|\bar{\rho}' |^2 \nabla\eta_2 )+{\mathcal{A}}_{11} \partial_2\partial_1\mathcal{N}^{\eta,1}_{2,2}+
{\mathcal{A}}_{11}\partial_1^2\mathcal{N}^{\eta,1}_{2,1}+{\mathcal{A}}_{12}\partial_2 \partial_k\mathcal{N}^{\eta,1}_{2,k} \nonumber\\
 &
-{\mathcal{A}}_{21}\partial_1^2\mathcal{N}^{\eta,1}_{1,1}
-{\mathcal{A}}_{22}\partial_2 \partial_1\mathcal{N}^{\eta,1}_{1,1}+  {\mathcal{A}}_{1l}\partial_l \partial_k \mathcal{N}^{\eta,2} _{2,k}-{\mathcal{A}}_{2l}\partial_l \partial_k \mathcal{N}^{\eta,2} _{1,k})  \partial_1^{2i}  \mm{curl} \eta\mm{d}y \nonumber \\
\lesssim &\| \eta\|_{2i,1}( \| \eta_2\|_{2,1}   \| \eta\|_3
 +\|  \eta_2\|_3^2  ) .
\label{20222091710032}
\end{align}

Making use of \eqref{improxtian1}, \eqref{2022210071415}, integration by parts, and the boundary condition in \eqref{01dsaf16asdfasf03n}$_4$, we have
\begin{align}
K_{2,i}=\,&  \kappa \int   |\bar{\rho}' |^2 \nabla \partial_1^{1+i}  \eta_2   \cdot \nabla \partial_2 \partial_1^{i}\eta_1 \mm{d}y- \kappa \|\bar{\rho}' \nabla   \eta_2 \|_{1+i, 0}^2  \nonumber \\
=\,& \kappa\int (    |\bar{\rho}' |^2 \nabla \partial_2\partial_1^i\eta_2\cdot\nabla \partial_1^i \mm{div}\eta +\partial_2  |\bar{\rho}' |^2 \nabla \partial_1^{i}\eta_2\cdot\nabla \partial_1^{1+i} \eta_1 )\mm{d}y- \kappa \|\bar{\rho}' \nabla^2   \eta_2 \|_{i,0}^2 \nonumber \\
\leqslant\,&  c \| \eta_2\|_{i,1}\|\eta\|_{1+i,1}+\|  \eta\|_3\| \eta\|_{1+i,1}\| \eta_2\|_{i,2}- \kappa \|\bar{\rho}' \nabla^2   \eta_2 \|_{i,0}^2. \label{2022asfd210181911}
\end{align}

By   integration by parts, we get
\begin{align}
K_{3,i} =
\sum_{j=1}^3 K_{3,i,j} -\frac{\mu }{2}\frac{\mm{d}}{\mm{d}t} \| \nabla_{\mathcal{A}}\partial_1^i \mm{curl} \eta \|_0^2 ,  \label{2022202041810}
  \end{align}
 where
 \begin{align*}
 K_{3,i,1} :=\,&\mu \int \nabla_{\mathcal{A}} \partial_1^i  \mm{curl} \eta  \cdot \nabla_{\mathcal{A}}\partial_1^i \mm{curl}_{\mathcal{A}_t} \eta \mm{d}y,\\
 K_{3,i,2}:=\,&
\mu \int  \nabla_{\mathcal{A}} \partial_1^i \mm{curl} \eta \cdot  \nabla_{\mathcal{A}_t}\partial_1^i \mm{curl}  \eta \mm{d}y,\\
 K_{3,i,3}:=\,&\begin{cases}
0&\mbox{for }\,\,i=0;\\
-\mu \int  \partial_1 ( \mathcal{A}_{kl}  \mathcal{A}_{kn}) \partial_n \mm{curl}_{\mathcal{A}}u\cdot   \partial_l \partial_1\mm{curl} \eta \mm{d}y&\mbox{for }\,\,i=1.
\end{cases}\end{align*}
From  \eqref{fgessfdims}, we easily get
\begin{align}
 \sum_{j=1}^3  K_{3,1,j}
\lesssim  \|\eta\|_{1,2}^2\| u\|_{\underline{1},2}+\|\eta\|_{1,2} \|\eta\|_3\| u\|_{2,1}.  \nonumber
\end{align}
In addition,
 \begin{align}  K_{3,0,1} + K_{3,0,2}
 \lesssim  \| \eta\|_2\|\eta\|_{1,2}\| u\|_2.
  \nonumber
\end{align}
Thanks to the above two estimates, we infer from \eqref{2022202041810} that
\begin{align}
K_{3,i}\leqslant& c (\delta_{i0}\|\eta\|_{1,2} \|\eta\|_2\| u\|_{2}+\delta_{i1}( \|\eta\|_{1,2}^2\| u\|_{\underline{1},2}\nonumber \\
&+\|\eta\|_{1,2} \|\eta\|_3\| u\|_{2,1}))-\frac{\mu }{2}\frac{\mm{d}}{\mm{d}t}
\|\nabla_{\mathcal{A}}\partial_1^i \mm{curl} \eta \|_0^2.
\label{2022202081350}
\end{align}

 Using \eqref{fgessfdims} again, we have
\begin{align}
K_{4,i}+K_{5,i}  \lesssim   \|u\|_{i,1}\|u\|_{\underline{i},1}+ \| \eta_2\|_{i,0} \| \eta \|_{1+i,1}+ \| \eta\|_3 \|  u\|_{2}^2    .
\label{20222020813fas50}
\end{align}
Inserting the estimates \eqref{20222091710032}, \eqref{2022asfd210181911}, \eqref{2022202081350} and \eqref{20222020813fas50} into \eqref{f202005031720}, and then employing \eqref{osdfadsa}, we arrive at \eqref{202005021600}.

(2) Now, we turn to the derivation of \eqref{202201211416}.  Multiplying \eqref{201910072117} by
${-\partial_2^2\mm{curl}\eta}$,
  integrating the results by parts, and   using  \eqref{202204120943}, \eqref{202220118000} and $\partial_1\mm{curl}_{\mathcal{A}}u|_{\partial\Omega}=0$
give
 \begin{align}
 &\frac{\mm{d}}{\mm{d}t} \left(  \frac{\mu }{2}\| \nabla \partial_2\mm{curl} \eta\|^2_0
 -\int \partial_2^2\mm{curl}  \eta \mm{curl}_{\mathcal{A}}(\bar{\rho}  u) \mm{d}y\right)
  = K_{6}+K_{7}+K_{8}, \label{f202005720}
\end{align}
where
\begin{align}
K_6 :=\,& {-\kappa \int (
\mm{curl}_{\mathcal{A}}\mathcal{N}^{\eta}+  \tilde{\mathcal{A}}_{1j} \partial_j   \mm{div}(|\bar{\rho}'|^2\nabla\eta_2))\partial_2^2\mm{curl} \eta\mm{d}y},\nonumber \\
 K_7 :=\,&  -\kappa  \int    \mm{div}(|\bar{\rho}' |^2 \nabla \partial_1\eta_2 ) \partial_2^2 \mm{curl}  \eta
    \mm{d}y,
\nonumber\\
 K_8  :=\, & \int\big\{ ( \mu \nabla \partial_2   \mm{curl} \eta \cdot \nabla \partial_2\mm{curl}_{\mathcal{A}_t}\eta
  +(\mu  ( \Delta-\Delta_{\mathcal{A}}) \mm{curl}_{\mathcal{A}} u\nonumber\\
& \quad - \mm{curl}_{\mathcal{A}_t} (\bar{\rho} u   ) )\partial_2^2 \mm{curl} \eta {-}
\mm{curl}_{\mathcal{A}} (\bar{\rho} u  )\partial_2^2 \mm{curl}  u  - g \bar{\rho}' \partial_2\eta_2
 \partial_1  \partial_2 \mm{curl}  \eta
  ) \big\}\mm{d}y .
\nonumber
\end{align}
Next, we shall estimate the above three integrals.

Similarly to \eqref{20222091710032}, we can estimate that
\begin{align}
 K_6 =\,&- \kappa \int
({\mathcal{A}}_{11} \partial_2\partial_1\mathcal{N}^{\eta,1}_{2,2} +
{\mathcal{A}}_{11}\partial_1^2\mathcal{N}^{\eta,1}_{2,1}+{\mathcal{A}}_{12}\partial_2 \partial_k\mathcal{N}^{\eta,1}_{2,k}
-{\mathcal{A}}_{21}\partial_1^2\mathcal{N}^{\eta,1}_{1,1} \nonumber
-{\mathcal{A}}_{22}\partial_2 \partial_1\mathcal{N}^{\eta,1}_{1,1} \\
&+{\mathcal{A}}_{1l}\partial_l \partial_k \mathcal{N}^{\eta,2} _{2,k}-{\mathcal{A}}_{2l}\partial_l \partial_k \mathcal{N}^{\eta,2} _{1,k} +  \tilde{\mathcal{A}}_{1j} \partial_j   \mm{div}(|\bar{\rho}'|^2\nabla\eta_2))\partial_2^2\mm{curl} \eta \mm{d}y\nonumber \\
\lesssim \,&  \| \eta\|_3 (\|  \eta_2\|_{2,1} \| \eta\|_3+\|  \eta_2\|_3^2). \label{2022209140110}
\end{align}

Exploiting the integral  {$K_{7}$} by parts and \eqref{202108261406}, we have
\begin{align}
K_7  = &\kappa\int ((\partial_2(|\bar{\rho}' |^2   \partial_1^2\eta_2  + 2 \bar{\rho}' \bar{\rho}'' \partial_2 \eta_2 ) -  2 \bar{\rho}' \bar{\rho}'' \partial_2   \partial_1 \eta_1  )\partial_2^2 \partial_1\eta_1  \nonumber \\
    & -( 2 \bar{\rho}' \bar{\rho}'' \partial_2\partial_1\eta_2 +  |\bar{\rho}' |^2    \partial_1^3\eta_2)\partial_2^2\partial_1 \eta_2 +    |\bar{\rho}' |^2    \partial_2\partial_1\mm{div}\eta \partial_2^3\eta_1
   ) \mm{d}y \nonumber\\
   &-\kappa\|\bar{\rho}'  \partial_2^2 \eta  \|_{1,0}^2- {2\kappa\int_{\partial\Omega}{\mathbf{n}}_2
  \bar{\rho}' \bar{\rho}'' \partial_2\eta_2\partial_1\partial_2^2\eta_1 \mm{d}y_1} \nonumber\\
  \lesssim &  \|\partial_2^2\eta_1\|_{1,0}(\|\partial_2\eta_1\|_{1,0}+\|(\eta_2,\partial_2\eta_2  )\|_{2,0}+\|(\partial_2\eta_2, \partial_2^2\eta_2 )\|_{0} )  \nonumber \\
  & + (
\|\partial_2\eta_2\|_{ {1},0}^{1/2}\|\partial_2^2\eta_2\|_{0}^{1/2}
+\|\partial_2\eta_2\|_{\underline{1},0}
)(
\|\partial_2^2\eta_1\|_{ {1},0}^{1/2}\|\partial_2^3\eta_1\|_{0}^{1/2}
+\|\partial_2^2\eta_1\|_{\underline{1},0}
) \nonumber \\
&+ \|\partial_2^2\eta_2\|_{1,0}(\|\partial_2\eta_2\|_{1,0}+\|\eta_2  \|_{3,0})
-\|  \partial_2^2 \eta  \|_{1,0}^2  +  \sqrt{\mathcal{E}}\mathcal{D} ,
\label{20240110} \end{align}
where $ {\mathbf{n}}_2$ denotes the second component of the outward unit normal vector $ {\mathbf{n}} $ to
$\partial\Omega$.
In addition (see the derivation of (2.41) in \cite{JFJSZWC}),
\begin{align}
 K_8  {\leqslant }& - \frac{\mm{d}}{\mm{d}t}\left( \mu\int  \partial_2^2  ( \partial_1\eta_1\partial_2 \eta_1)
\partial_2^2\mm{curl} \eta  \mm{d}y+\mathcal{I}_1\right)\nonumber \\
&+ {c(\| \partial_2\eta_1\|_{1,0}\|\partial_2^2\eta_2\|_0
+\| \eta_2\|_{2,1}\|\eta_2\|_{1}}+\|u\|_1\|  u\|_3+\sqrt{\mathcal{E}}\mathcal{D}).
\label{2022201191657}
\end{align}

Finally, substituting \eqref{2022209140110}--\eqref{2022201191657}  into \eqref{f202005720}, and then using the estimate  {\eqref{osdfadsa}} and Young's inequality, we arrive at \eqref{202201211416}. This completes the proof.
\hfill$\Box$
\end{pf}

\subsection{Total energy estimates}\label{202220201221244}
Now, we are in a position to establish the \emph{a priori} stability estimate  \eqref{1.200}  under the assumptions \eqref{apresnew}--\eqref{aprpiosesnew}.

We can infer from  Lemmas \ref{lem:082sdaf41545}--\ref{lem:08241445} and the estimate \eqref{20222saf201121235} that
 there are two constants $c$ and (suitably large) $\tilde{\chi}$, such that for any $\chi\geqslant\tilde{\chi}$ and
any sufficiently small $\delta\in (0,1]$,  the following tangential energy inequality holds:
  \begin{align}&\label{202008250856n0}
\frac{\mm{d}}{\mm{d}t} {\mathcal{E}}_{\mm{tan}}
+c  {\mathcal{D}}_{\mm{tan}}  \lesssim \chi\sqrt{\mathcal{E}}\mathcal{D} ,
\end{align}
where
$$
\begin{aligned}
&\begin{aligned}
{\mathcal{E}}_{\mm{tan}} :=\,& \chi
 \sum_{0\leqslant i \leqslant 2}\left(\frac{\kappa}{2}\|\bar{\rho}'  \partial_2\eta_1  \partial_1^{1+i}\eta_2\|^2_0   -\kappa\int |\bar{\rho}'  |^2   \partial_2\eta_1 \partial_1^{1+i}\eta_2\partial_2\partial_1^{i}\eta_2   \mm{d}y -E(\partial_1^{i} \eta)\right) \\
&+\chi   \|\sqrt{ \bar{\rho} }  u\|^2_{\underline{2},0} +\sum_{0\leqslant i \leqslant 2} {\bigg{(}}\int\bar{\rho}\partial_1^{i}\eta\cdot\partial_1^{i} u\mm{d}y
+\frac{\mu}{2}\| \nabla   \eta\|_{i,0}^2{\bigg{)}} +\|\sqrt{\bar{\rho}}\psi\|_0^2 - E(  u),\end{aligned}\\
  &  {\mathcal{D}}_{\mm{tan}} :=\chi\|  u\|_{\underline{2},1}^2 -\sum_{0\leqslant i\leqslant 2 }E(\partial_1^{i}\eta) + \|   u_t\|_1^2  .
\end{aligned}
$$

On the other hand,  if we make use of Lemma \ref{lem:08250749}, \eqref{202211}, \eqref{omessetsim} and Young's inequality,
we obtain that for any sufficiently large $\chi$ and any sufficiently small $\delta\in (0,\chi^{-1})$,
\begin{align}
&\| \eta_1\|_{\underline{2},1}^2 +\|u\|_2^2+ \|u_t \|_0^2+\|\tilde{q} \|_1^2+ \chi (  \|  \eta_2 \|^2_{\underline{2},1} +\| u \|^2_{\underline{2},0}  )\lesssim  {\mathcal{E}}_{\mm{tan}} \lesssim \chi \mathcal{E}
\label{20221057}
\end{align}
and
\begin{align}
   \| \eta_2\|_{\underline{2},1}^2+ \chi\|  u\|_{\underline{2},1}^2  +\|u_t\|_1^2
   \lesssim  {\mathcal{D}}_{\mm{tan}} .  \label{2022202012047}
\end{align}

Thus  we can further derive
the total energy inequality from Lemma  \ref{2055nnn} and \eqref{202008250856n0}:
  \begin{align}
\frac{\mm{d}}{\mm{d}t} \tilde{\mathcal{E}} +c \tilde{\mathcal{D}} \lesssim\,&
  \| \eta_2\|_{2}^2 +   \| u\|_1\| u\|_3+ (
\|\partial_2\eta_2\|_{ {1},0}^{1/2}\|\partial_2^2\eta_2\|_{0}^{1/2}+\|\partial_2\eta_2\|_{\underline{1},0}
)(
\|\partial_2^2\eta_1\|_{ {1},0}^{1/2}\|\partial_2^3\eta_1\|_{0}^{1/2} \nonumber\\
&
+\|\partial_2^2\eta_1\|_{\underline{1},0}
)+\|\eta_2\|_{2,1}\|\eta\|_3^2+  {\chi^3} \sqrt{\mathcal{E}}\mathcal{D} ,\label{202008250856}
\end{align}
where we have defined that
\begin{align*}
&
\tilde{\mathcal{E}}:=
{\mathcal{I}_3 + \frac{\mu}{2} \| ( {\chi^{\frac{1}{2}}} \nabla_{\mathcal{A}} \mm{curl}  \eta, {\chi^\frac{1}{2}}\nabla_{\mathcal{A}}\partial_1
\mm{curl}  \eta,\nabla  \partial_2 \mm{curl} \eta)\|^2_0+\chi^2  {\mathcal{E}}_{\mm{tan}} },     \\
& \tilde{\mathcal{D}} := \|  \partial_2^2 \eta  \|_{1,0}^2 + \chi  \|  \nabla^2 \eta_2 \|^2_{\underline{1},0} + \chi^2  {\mathcal{D}}_{\mm{tan}}  ,  \\
&  {\mathcal{I}_3} :=  {\mathcal{I}_1}  + \chi\int( \mm{curl} \eta\,\mm{curl}_{\mathcal{A}}(\bar{\rho} u )
 - \partial_1^2\mm{curl} \eta\,\mm{curl}_{\mathcal{A}}(\bar{\rho} u) )\mm{d}y   .
\end{align*}
Moreover, by \eqref{omessen}, \eqref{omessetsim122n} and \eqref{2022202012047}, we see that for any sufficiently small
$\delta\in (0,\chi^{-1})$,
\begin{align}
 {\mathcal{D}} \lesssim \,&\|    \eta_1\|_{1,2}^2+ \|    \eta_2\|_{3}^2+\| u\|_{3}^2+\|\tilde{q}\|_2^2+\chi \|\partial_2^2    \eta_2\|_{\underline{1},0}^2+  \chi^2(\| \eta_2\|_{\underline{2},1}^2
 + \chi\|  u\|_{\underline{2},1}^2+\|   u_t\|_1^2) \lesssim  \tilde{\mathcal{D}} . \label{20x012047}
\end{align}

Noticing that
 \begin{align}
 \|   \mm{curl}  \eta\|^2_2
 \lesssim\,&\|\nabla \eta\|_0^2 + \|  \nabla \mm{curl}  \eta\|^2_{1}\nonumber \\
 \lesssim\, &\|(\nabla \eta, \nabla_{\mathcal{A}} \mm{curl}  \eta,
  \nabla_{\mathcal{A}}\partial_1 \mm{curl} \eta,
 \nabla \partial_2 \mm{curl} \eta, \nabla_{\tilde{\mathcal{A}}}  \mm{curl}  \eta, \nabla_{\tilde{\mathcal{A}}}\partial_1\mm{curl} \eta      )\|^2_0  \nonumber \\
\lesssim\,& \|(\nabla \eta,\nabla_{\mathcal{A}} \mm{curl}  \eta, \nabla_{\mathcal{A}} \partial_1 \mm{curl}  \eta , \nabla \partial_2 \mm{curl} \eta)\|^2_0 + \| \nabla \eta\|^4_2  ,
\label{202202172007}
\end{align}
 we get from \eqref{2022202011749} that
\begin{align}
 \|  \eta \|_3^2 \lesssim \|(\nabla \eta,\nabla_{\mathcal{A}} \mm{curl} \eta, \nabla_{\mathcal{A}}
 \partial_1 \mm{curl}  \eta , \nabla \partial_2 \mm{curl} \eta)\|^2_0 . \nonumber
\end{align}
 In addition, it is easy to verify that
 \begin{align}
&| {\mathcal{I}_3}| \lesssim \chi
\| \eta\|_3\|u\|_1+\|\eta\|_3^3. \nonumber
\end{align}
 From the above two estimates, \eqref{202012252005}, \eqref{20221057} and the interpolation inequality \eqref{201807291850}, it follows that  for any sufficiently large $\chi$ and
any sufficiently small $\delta$,
\begin{align}
&\mathcal{E}\lesssim \|\eta\|_3^2+ \chi^2  ( \|u\|_2^2+ \|u_t \|_0^2+\|\tilde{q}\|_1^2 ) + \chi^3\|u\|_{\underline{2},0}^2
\lesssim   \tilde{\mathcal{E}}\lesssim \chi^3\mathcal{E}\lesssim \chi^3\|(\nabla \eta,u)\|_2^2.
\label{2027}
\end{align}

With the help of \eqref{20x012047}, the interpolation inequality \eqref{201807291850} and Young's inequality, we further deduce \eqref{for:0202n} for any sufficiently small $\delta$  from \eqref{202008250856} with some appropriately
large $\chi$,
where $\tilde{\mathcal{E}}$ satisfies \eqref{2027}. Integrating \eqref{for:0202n}  over $(0,t)$, we arrive at \eqref{1.200xyx}
which implies \eqref{as1.200xyx}.
\subsection{Decay-in-time estimates}\label{202244}

This subsection is devoted to the derivation of the  decay-in-time estimates in $\mathfrak{E}$. By integration by parts, we get from Lemmas \ref{lem:082sdaf41545} and \ref{lem:08241445} that for $i=1, 2$,
\begin{align}
& \frac{\mm{d}}{\mm{d}t}\mathcal{E}_{1,\mm{D}}+c\sum_{i=1}^2 \langle t\rangle^i\big\{\gamma \|   u\|_{i,1}^2 -  E(\partial_1^{i}\eta) \big\} \nonumber \\
&  \lesssim
 \sum_{i=1}^2\langle t\rangle^{i-1}  (\| ( \nabla \partial_1^{i} \eta,\partial_1^{i-1} u  ) \|_{0}^2+ \gamma (\| u\|^2_{i,0}-E(\partial_1^{i}\eta) ))+ \mathcal{I}_4,
\label{2020046}
\end{align}
where $\gamma \geqslant 1$ is an appropriately large constant (may depending on $\bar{\rho}$, $\mu$, $g$, $\kappa$ and $\Omega$),
\begin{align*}
 \mathcal{E}_{1,\mm{D}}:=\,&\sum_{i=1}^2 \langle t\rangle^i\left\{ \int\bar{\rho}\partial_1^{i}\eta\cdot\partial_1^{i} u\mm{d}y
+\frac{\mu}{2}\| \nabla \partial_1^{i} \eta\|_{0}^2 + \gamma  \left(\|\sqrt{ \bar{\rho} }  u\|^2_{i,0}
-E(\partial_1^i \eta) \right)\right.  \nonumber \\
&  \qquad + \left.\gamma \kappa\left(\frac{1}{2}\|\bar{\rho}'  \partial_2\eta_1  \partial_1^{1+i}\eta_2\|^2_0   -\int |\bar{\rho}'  |^2   \partial_2\eta_1 \partial_1^{1+i}\eta_2\partial_2\partial_1^{i}\eta_2   \mm{d}y\right)\right\}  \nonumber    \end{align*}
and
\begin{align*}
 \mathcal{I}_4:=\,&
\langle t\rangle^2\big\{
 \|\eta\|_{3,0}( \|  \eta\|_3\|  u\|_{2,1}+ \|  \eta\|_{ {2,1}}\|  u\|_{\underline{1},2} ) +  \|\eta_2\|_{2,1}( \|  \eta_2 \|_{2,1}  \|  \eta  \|_3 \nonumber \\
 &\quad+\|\eta_2\|_{ {1},2}\|\eta_2\|_3) +\|\eta_2\|_{1,2}^3 +  \|\eta_2\|_{2,1}\|\eta_2\|_{1,2} \|\tilde{q}\|_{\underline{1},1}\big\}
 \nonumber \\
& \quad+ \gamma \langle t\rangle^3   \|\eta_2\|_{ {2},1} ( \|   \eta_2  \|_3 \|  u \|_{2,1}+
     \|   \eta_2  \|_{2,1} \|  u \|_{\underline{1},2} )
     \nonumber \\ & \quad+ (\|  \eta_2\|_{ {1},2}^2 +\|  \eta\|_{2,1}\|u\|_{\underline{1},2}  +\|\eta \|_{3,0}\|u\|_{3})\|u\|_{2,1}
\nonumber \\
 &\quad +  \|  \eta\|_{2,1}\|  u\|_{2,1}+\|  \eta\|_{3,0 }( \| \eta_2 \|_{ {1},2} \|u\|_{\underline{1},2 }+\|  \eta\|_{3}\| u\|_{ 2,1}))\| \tilde{q}\|_{\underline{1},1}\nonumber \\
 &\quad+\|  \eta\|_{3,0}\|u\|_{\underline{1},2 }(\|  \eta_2 \|_{2,1} +  \| u_t\|_{1}
 + \|  \eta_2 \|_{1,2}\|\eta_2\|_3 + \|\eta\|_{2,1}\|u\|_{\underline{1},2}
   )
 .   \end{align*}

If we make use of  \eqref{osdfadsa} and  \eqref{202008250745}, we can deduce from \eqref{202008241448} with $i=2$,
\eqref{202005021600} with $i=1$ and \eqref{2020046} that
\begin{align}
& \frac{\mm{d}}{\mm{d}t}\left(\mathcal{E}_{\mm{D}}+\gamma^2 \langle t\rangle
 \int\partial_1\mm{curl}    \eta
 \partial_1\mm{curl}_{\mathcal{A}} \left(\bar{\rho}u\right)
 \mm{d}y\right)  +c \mathcal{D}_{\mm{D}} \nonumber\\
&\lesssim  \gamma   \big\{\gamma ( \|    \eta\|_{1,2}^2+ \|u\|_{\underline{1},1}^2) + \langle t\rangle( \| ( \nabla \partial_1^2 \eta,\partial_1 u  ) \|_{0}^2 +     \gamma (\|\eta_2\|_{2,1}^2 +\|   u\|^2_{2,0}\nonumber\\
&\quad +\|  u\|_{1,1}\|u\|_{2})  )\big\}+ \langle t\rangle^2 (\|   u\|^2_{2,0} +\|\eta_2\|_{2,1}^2  ) + \gamma   \mathcal{I}_4+ \mathcal{I}_5,
\label{2020046xx}
\end{align}
where we have defined that
\begin{align*}
\mathcal{E}_{\mm{D}}:=\,&\gamma \mathcal{E}_{1,\mm{D}}+\mu  \gamma^2 \langle t\rangle
 \|\nabla_{\mathcal{A}} \partial_1\mm{curl} \eta \|_0^2 /{2}  \nonumber \\
&\quad  + \langle t\rangle^{3}\left(\|\sqrt{ \bar{\rho} }  u\|^2_{2,0}
-E(\partial_1^{2} \eta)+ \frac{\kappa}{2}\|\bar{\rho}'  \partial_2\eta_1  \partial_1^3\eta_2\|^2_0   -\kappa\int |\bar{\rho}'  |^2   \partial_2\eta_1 \partial_1^3\eta_2\partial_2\partial_1^2\eta_2   \mm{d}y \right)\nonumber \\
\mathcal{D}_{\mm{D}}:=\, &\gamma\sum_{i=1}^2 \langle t\rangle^i \| {( \eta_2,  \sqrt{\gamma}u)}\|_{ {i},1}^2 + \langle t\rangle^3 \|   u\|_{2,1}^2 +\gamma^2 \langle t\rangle\|\nabla^2 \eta_2\|^2_{1,0} ,\nonumber \\
\mathcal{I}_5:=\,&\gamma^2 \langle t\rangle\big\{\|  \eta\|_3(\|\eta_2\|_{1,2}^ 2+ \|\eta\|_{2,1}^2+\|\eta\|_{1,2}\|u\|_{2,1}
+\|u\|_2^2 )+\|\eta\|_{1,2}^2\|u\|_{\underline{1},2}+\|\eta\|_{2,0}\|\eta_2\|_3^2\big\}.
\end{align*}

In the same manner as used in the derivation of \eqref{202202172007}, we have
\begin{align}
 \|   \mm{curl}  \eta\|^2_{1,1}
 \lesssim&   \|(\nabla \partial_1\eta, \nabla_{\mathcal{A}} \partial_1 \mm{curl} \eta  )\|^2_0 + \| \eta\|^2_{1,2}   \| \nabla \eta\|^2_2, \nonumber
\end{align}
which, combined with \eqref{2022202011749}, implies that
\begin{align}
 \|\eta\|_{1,2}^2\lesssim\|(\nabla \partial_1\eta, \nabla_{\mathcal{A}}\partial_1\mm{curl} \eta )\|^2_0. \nonumber
\end{align}
In view of the above estimate,   \eqref{202008250745} and \eqref{202012241002}, we further find that
 for any given sufficiently large $\gamma $,
\begin{align}
\begin{cases}
\langle t\rangle   \| \partial_2^2 \eta\|_{1,0}^2  + \langle t\rangle^2 \| \partial_2\partial_1\eta\|_{\underline{1},0}^2
+ \langle t\rangle^3  \|\eta_2 \|_{\underline{2},1}^2 \\
 \qquad -c \gamma^2 \langle t\rangle^3 \|\eta_2\|_{2,1}^2 \|\eta\|_3\lesssim  \mathcal{E}_{\mm{D}}\lesssim \gamma^2\langle t\rangle^3\|(\nabla\eta,u)\|_2^2 , \\[1mm]
 \gamma^2\langle t\rangle  \| \nabla^2 \eta_2\|_{1,0}^2 +\gamma \langle t\rangle^2    \|  \eta_2\|_{\underline{2},1}^2+ (\gamma^2\langle t\rangle^2+ \langle t\rangle^3) \| u\|_{ 2,1}^2\lesssim  \mathcal{D}_{\mm{D}}.
 \end{cases}
\label{20222020171604}
\end{align}

 On the other hand, we easily see that
$$
 {\left|\int\partial_1\mm{curl}   \eta  \partial_1\mm{curl}_{\mathcal{A}} \left(\bar{\rho}u\right)
 \mm{d}y\right|\lesssim  \|\eta\|_{1,1}\| u\|_2}.
$$
Hence, integrating \eqref{2020046xx} with some suitably large $\gamma$ over $(0,t)$, and using then \eqref{as1.200xyx}, \eqref{omessen}, \eqref{osdfadsa}, \eqref{20222020171604}, \eqref{202012241002},
the above estimate and Young's inequality, we infer that for any sufficiently small $\delta$,
\begin{align}
&\langle t\rangle \|\partial_2^2( \partial_1\eta_1,\partial_2\eta_2)\|_0^2  + \langle t\rangle^2  \|
  \partial_2( \partial_1\eta_1,\partial_2\eta_2 )\|_{\underline{1},0 }^2
 + \langle t\rangle^3(\| \eta_2 \|_0^2 +\|(\partial_2\eta_2,\partial_1\eta)\|_{\underline{2},0}^2 ) \nonumber  \\[1mm]
 &+\int_0^t\big\{\langle\tau\rangle   \|\partial_2(\partial_1\eta_1, \partial_2\eta_2 )\|_{\underline{1},0}^2 +   \langle \tau\rangle^2  ( \| \eta_2 \|_{0}^2
 + \|(\partial_2\eta_2,\partial_1\eta)\|_{\underline{2},0}^2   ) + \langle \tau\rangle^3 \|\partial_1 u\|_{\underline{1},1}^2\big\}\mm{d}\tau \nonumber \\
\lesssim  &\|(\nabla\eta^0,u^0)\|_2^2  +\int_0^t(   \|u\|_{2}^2 +    \|\eta \|_{1,2}^2
   +  \mathcal{I}_4+ \mathcal{I}_5)\mm{d}\tau \nonumber \\
   &+ \sup_{0\leqslant \tau\leqslant t} \left(\langle\tau\rangle^{3/2}   \|  \eta_2(\tau)\|_{2,1}  \mathcal{E}(\tau) + \sqrt{\mathfrak{E}(\tau)
\|\partial_2^2\eta_1(\tau)\|_{\underline{1},0} }
 \mathcal{E}^{1/4}(\tau)\right). \label{2022235}
\end{align}

 It easily follows from \eqref{202008241448}, \eqref{202221235} and \eqref{20222saf201121235}  that for any sufficiently large $\alpha $ (may depend on  $\bar{\rho}$, $\mu$, $g$, $\kappa$ and $\Omega$),
\begin{align}
&\frac{\mm{d}}{\mm{d}t}\left(\sum_{i=1}^2\alpha^{4-i} \langle t\rangle^i\left(\|\sqrt{ \bar{\rho} }  u\|^2_{0}
-E( \eta) +\frac{\kappa}{2}\|\bar{\rho}' \partial_2\eta_1  \partial_1\eta_2\|^2_0-\kappa\int |\bar{\rho}' |^2   \partial_2\eta_1 \partial_1\eta_2\partial_2\eta_2   \mm{d}y\right)\right.\nonumber\\
& \left.+    \alpha\langle t\rangle^2\left(\|\nabla_{\mathcal{A}} u \|^2_0+ \kappa\int  |\bar{\rho}' |^2 \nabla\eta_2 \cdot\nabla u_2 \mm{d}y \right)+\langle t\rangle^3 (\|\sqrt{\bar{\rho}}\psi\|^2_0- E(  u))\right)\nonumber \\
& +
c( \alpha\langle t\rangle^2(\alpha\|u\|_1^2+\|   u_t\|_0^2)+\langle t\rangle^3 \|   u_{t}\|_1^2)\nonumber \\
&\lesssim \alpha^3(\|\eta_2\|_1^2+\|u\|_0^2 + \langle t \rangle \|  \eta_2 \|_1 \| \eta_2\|_3\|  u\|_1) +\langle t\rangle^2(\alpha ( \| \eta_2  \|_{1}^2\nonumber \\
&\quad+\|u\|_2^3+\|u\|_2^2\|\tilde{q}\|_1+\alpha\|\eta_2\|_1\|\eta_2\|_3\|u\|_1) +\|u\cdot\nabla_{\mathcal{A}}u\|_0^2 ) +  \langle t\rangle^3((\|  \eta_2\|_2 ^2 \nonumber \\
&\quad+\|u\|_2^2)\|u\|_2^2 + {\|\eta\|^2_3(\|u\|^2_{1,0}+\|\eta_2\|^2_{2,1}\|u\|^2_1)}
+(\|u\|_{1,0}+\|\eta_2\|_{2,1}\|u\|_1)\|u\|_2^2  +\|u\|_2^4  ).\nonumber
\end{align}

If we integrate the above estimate over $(0,t)$, make use of \eqref{aprpiosesnew}, \eqref{2022111522311}, \eqref{202012252005}, Lemma \ref{lem:08250749} and the product estimate \eqref{fgestims},
we derive from the above estimate that
 \begin{align}
&   \langle t\rangle^2\|  u \|^2_1 +\langle t\rangle^3 \|u_t\|^2_0  +
c\int_0^t( \langle \tau\rangle^2 \|u\|_1^2  +   \langle    \tau \rangle^3 \| u_{\tau}\|_1^2 )\mm{d}\tau\nonumber \\
&\lesssim   \| (\nabla\eta^0, u^0) \|^2_2 +\langle t\rangle^3 \|u\|^4_2  + \int_0^t ( \langle \tau\rangle^{2}    \| \eta_2 \|_{1}^2 +\|u\|_0^2+\mathcal{I}_6) \mm{d}\tau,
\label{202225}
\end{align}
where we have defined that
\begin{align*}\mathcal{I}_6:= &  \|  \eta_2 \|_1^2\| \eta_2\|_3^2  +\langle \tau\rangle^{2} ( \|\eta_2\|_1\|\eta_2\|_3\|u\|_1+\|u\|_2^3+\|u\|_2^2\|\tilde{q}\|_1) + \langle \tau\rangle^3((\|  \eta_2\|_2 ^2 +\|u\|_2^2)\|u\|_2^2\\
& +\|  \eta \|_{3}^2 (\|u\|_{1,0}^2+\|\eta_2\|_{2,1}^2\|u\|_1^2)  + (\|u\|_{1,0}+\|\eta_2\|_{2,1}\|u\|_1)\|u\|_2^2  +\|u\|_2^4 ).
\end{align*}

  Employing  \eqref{omessetsim}
and Young's inequality, we further get from \eqref{202225} with some sufficiently large $\alpha $ that
 \begin{align}
&\langle t\rangle^2(\|u\|_2^2+ \|(q,\tilde{q})\|_1^2)+ \langle t\rangle^3\|u_t\|^2_0
\nonumber \\
&  +
c\int_0^t(\langle \tau \rangle( \|u\|_{\underline{1},2}^2+\|(q,\tilde{q})\|_{\underline{1},1}^2)+ \langle \tau \rangle^2 \|u\| _1^2
+\langle \tau \rangle^3 \|   u_{\tau}\|_1^2)\mm{d}\tau   \nonumber \\
&\lesssim  \| (\nabla\eta^0, u^0) \|^2_2+\langle t\rangle^2\|\eta_2\|_2^2+\langle t\rangle^3 \|u\|^4_2 \nonumber \\
&\qquad+ \int_0^t ( \|u\|_0^2+    \langle \tau\rangle     \| \eta_2 \|_{ {1},2}^2+    \langle \tau\rangle^{2}    \| \eta_2 \|_{1}^2+ \mathcal{I}_6  ) \mm{d}\tau,
\nonumber
\end{align}
which, together with \eqref{2022235}, yields
 \begin{align}
\mathfrak{E} (t) +
 \int_0^t \mathfrak{D} {(\tau)}\mm{d}\tau \lesssim  & \| (\nabla\eta^0, u^0) \|^2_2 +\langle t\rangle^3  \|u\|^4_2 +\int_0^t(  \|    \eta\|_{1,2}^2+\|u\|_2^2+\mathcal{I} ) \mm{d}\tau\nonumber \\
 &+ \sup_{0\leqslant \tau\leqslant t} \left(\langle\tau\rangle^{3/2}   \|  \eta_2(\tau)\|_{2,1}  \mathcal{E}(\tau) + \sqrt{\mathfrak{E}(\tau)
\|\partial_2^2\eta_1(\tau)\|_{\underline{1},0} }
 \mathcal{E}^{1/4}(\tau)\right), \label{20222010041821}
\end{align}
where we have defined that
\begin{align}\mathcal{I}:=\mathcal{I}_4+\mathcal{I}_5+  \mathcal{I}_6.
\label{202222100718123}
\end{align}

\subsection{Energy estimate for $\|\partial_2^2\eta_1\|_{\underline{1},0}^2$}
Now we further establish the estimate of   $\|\partial_2^2\eta_1\|_{\underline{1},0}^2$ in \eqref{as1.200xyx}.
To begin with, we can derive the following estimate  from
\eqref{202005021600} and  \eqref{202008250856n0} with $ {\mathcal{E}}_{\mm{tan}}$ and $ {\mathcal{D}}_{\mm{tan}}$ satisfying  \eqref{20221057} and  \eqref{2022202012047}:
 \begin{align}
& \frac{\mm{d}}{\mm{d}t}\left( \frac{\mu}{2} \| ( \nabla_{\mathcal{A}} \mm{curl}  \eta, \nabla_{\mathcal{A}}\partial_1
\mm{curl}  \eta )\|^2_0+\chi   \mathcal{E}_{\mm{tan}} +\sum_{i=0}^1
 \int  \partial_1^i \mm{curl}   \eta
  \partial_1^i \mm{curl}_{\mathcal{A}} \left(\bar{\rho}u\right)
 \mm{d}y  \right) \nonumber \\
 &\qquad \lesssim  \chi^2\sqrt{\mathcal{E}}\mathcal{D}+\|\eta_2\|_{ {2},1}\|\eta\|_1\|\eta\|_3 .
\label{20200sdfa46xx}
\end{align}
It is easy to check that
\begin{align}
\|\eta\|_{\underline{1},2}^2\lesssim\, &  \frac{\mu}{2} \| ( \nabla_{\mathcal{A}} \mm{curl}  \eta, \nabla_{\mathcal{A}}\partial_1
\mm{curl}  \eta )\|^2_0+\chi   \mathcal{E}_{\mm{tan}}
\nonumber \\
 &+\sum_{i=0}^1
 \int  \partial_1^i \mm{curl}   \eta
  \partial_1^i \mm{curl}_{\mathcal{A}} \left(\bar{\rho}u\right)
 \mm{d}y
\lesssim \chi^2 \|(\nabla \eta,u)\|_2^2.\nonumber
\end{align}
Integrating \eqref{20200sdfa46xx} over $(0,t)$, and then using the above relation, we arrive at
 \begin{align}
&
\|\partial_2^2\eta_1\|_{\underline{1},0}^2 \lesssim \|(\nabla\eta^0,u^0)\|_2^2+\int_0^t(\|\eta_2\|_{ {2},1}\|\eta\|_1\|\eta\|_3 +\sqrt{\mathcal{E}}\mathcal{D})\mm{d}\tau.
\label{20200asd46safdxx}
\end{align}

\subsection{Stability estimates}

 Putting \eqref{as1.200xyx} and \eqref{20222010041821} together,  and then using \eqref{202012252005} and Young's inequality, we obtain
\begin{align}
& \sup_{0\leqslant \tau\leqslant t}(\mathfrak{E}  (\tau)+ {\mathcal{E}}(\tau)) + \int_0^t(\mathfrak{D}(\tau)+ {\mathcal{D}}(\tau))\mm{d}\tau \nonumber \\
&\lesssim \|(\nabla \eta^0,u^0)\|_2^2+  \sup_{0\leqslant \tau\leqslant t}
(\langle\tau\rangle^{3/2}\|  \eta_2(\tau) \|_{2,1}  \mathcal{E}(\tau)+
\|\partial_2^2\eta_1(\tau)\|_{1,0}^2)
+\langle \tau\rangle^3 \|u\|^4_2  )+\int_0^t  \mathcal{I}  \mm{d}\tau  , \nonumber
\end{align}
which, together with \eqref{20200asd46safdxx}, yields \eqref{1.2asdf00}.

Now we assume that $(\eta,u)$ satisfies \eqref{202220201915},   we immediately derive from the estimate \eqref{1.2asdf00} satisfied by $(\eta,u,q)$ that
\begin{align}
\sup_{0\leqslant \tau\leqslant t} (\mathfrak{E} (\tau)+ {\mathcal{E}}(\tau))+ \int_0^t(\mathfrak{D}(\tau)+ {\mathcal{D}}(\tau))\mm{d}\tau \leqslant\tilde{c} \|(\nabla\eta^0,u^0)\|_2^2 ,\label{estemsfsadfaalasn0}
\end{align}
where the positive constant $\tilde{c}\geqslant 1$  depends on $\bar{\rho}$, $\mu$, $g$, $\kappa$   and $\Omega$.

Finally, we sum up the \emph{a priori} stability estimates with decay-in-time as follows.
\begin{pro}Let $(\eta,u,q)$ be a solution to the TCRT problem of \eqref{01dsaf16asdfasf} and \eqref{20safd45}, and satisfy \eqref{apresnew}--\eqref{aprpiosesnew}. If  $\bar{\rho}$  and $\kappa$ satisfy the assumptions in Theorem \ref{thm2}, {then} there is a constant $\delta_1$,
depending possibly on $\bar{\rho}$, $\mu$, $g$, $\kappa$ and $\Omega$, such that the solution $(\eta,u,q)$ enjoys
the estimate \eqref{estemsfsadfaalasn0} for any $\delta\leqslant \delta_1/\sqrt{2}$.
\end{pro}

\subsection{Proof  of Theorem \ref{thm2}}\label{subsec:08}

In this subsection, we prove Theorem \ref{thm2}. First, we state a local well-posedness result for the TCRT problem.
\begin{pro}\label{202102182115}
Let $b>0$ be a constant and $\gamma >0$ the same constant as in Lemma \ref{pro:1221}. Assume that $\bar{\rho}$ satisfies \eqref{0102},
$( \eta^0,u^0)\in\mathcal{H}^{3}_{\mm{s}}\times {\mathcal{H}^2_{\mm{s}}}$, $\|(\nabla \eta^0,u^0)\|_2\leqslant b$ and $\mm{div}_{\mathcal{A}^0}u^0=0$, where $\mathcal{A}^0:=(\nabla\eta^0+\mathbb{I})^{-\top }$.
Then, there is a sufficiently small constant $\delta_2\leqslant \gamma /2$, such that if $\eta^0$ satisfies
\begin{align}
\|\nabla \eta^0\|_2\leqslant \delta_2, \nonumber
\end{align}
 the TCRT problem of \eqref{01dsaf16asdfasf} and \eqref{20safd45} admits a unique local-in-time classical solution
$( \eta, u,q)\in {C}^0(\overline{I_{T}}, \mathcal{H}^{3}_{\mm{s}} )\times {\mathcal{U}}_T \times (C^0(\overline{I_T},
\underline{H}^1)\cap L^2_T{H}^2)$ for some $T>0$. Moreover, $(\eta,u) $ satisfies
\footnote{ {Here the uniqueness means that if there is another solution
$(\tilde{u}, \tilde{\eta},\tilde{q})\in {C}^0 (\overline{I_{T }},\mathcal{H}^{3}_{\mm{s}} )\times \mathcal{U}_{T }
\times (C^0(\overline{I_T} ,\underline{H}^1)\cap L^2_T {{H}^2})$
 satisfying $0<\inf_{(y,t)\in \Omega_T} \det(\nabla \tilde{\eta}+I)$, then
 $(\tilde{\eta},\tilde{u},\tilde{q})=(\eta,u,q)$ by virtue of the smallness condition
 ``$\sup_{t\in \overline{I_T}}\|\nabla \eta\|_2\leqslant 2\delta_2$''. In addition,
  we have, by the fact ``$\sup\nolimits_{t\in \overline{I_T}}\|\nabla \eta\|_2\leqslant  \gamma $" and  Lemma \ref{pro:1221} that
$$\inf\nolimits_{(y,t)\in \Omega_T} \det(\nabla \eta+I)\geqslant 1/4.$$}}
$$ \sup\nolimits_{t\in  {I_T}} \| \nabla \eta\|_2\leqslant 2\delta_2 , $$
where $\delta_2$ and $T$ may depend on $\mu $,  $g$,   $\bar{\rho}$ and $\Omega$, while $T$ further depends on $b$.
\end{pro}
\begin{pf}{Since} Proposition \ref{202102182115} can be easily proved by the standard iteration method as in  \cite{JFJSZYYO}, thus we omit the trivial proof.
\hfill $\Box$
\end{pf}
\begin{rem}\label{2022202021617}
If the initial data $(\eta^0,u^0)$ in Proposition \ref{202102182115} additionally satisfies $\det (\nabla \eta^0+\mathbb{I})=1$ and $( \bar{\rho}\eta_1^0)_\Omega= ( \bar{\rho}u_1^0)_\Omega= 0$, then $(\eta,u)\in \widetilde{\mathfrak{H}}^{1,3}_{\gamma,T}\times  {^0\mathcal{U}_T}$.
\end{rem}

 Due to \eqref{estemsfsadfaalasn0} and Proposition \ref{202102182115}, we can easily establish the global solvability
in Theorem \ref{thm2}. Next, we briefly describe the proof.

Let $(\eta^0,u^0)$ satisfy the assumptions in Theorem \ref{thm2} with $\delta=\min\left\{ \delta_1, \delta_2  \right\}/ \sqrt{2\tilde{c}}$,
where  $\tilde{c}$ is the  constant  as in  \eqref{estemsfsadfaalasn0}.
Noting that
\begin{align}
\|(\nabla \eta^0,u^0)\|_2 \leqslant\delta \leqslant \delta_2/\sqrt{2\tilde{c}} <\delta_2, \nonumber
\end{align}
thus, by virtue of Proposition \ref{202102182115} and Remark \ref{2022202021617}, there exists a unique local solution $(\eta,u, {q})$ to the TCRT problem of
\eqref{01dsaf16asdfasf} and   \eqref{20safd45} with the maximal existence time $T^{\max}$, which satisfies
\begin{itemize}
  \item for any $a\in I_{T^{\max}}$,
the solution $(\eta,u, {q})$ belongs to $\widetilde{\mathfrak{H}}^{1,3}_{\gamma,a}\times {^0\mathcal{U}}_{a} \times ( C^0(\overline{I_a }, \underline{H}^1)\cap L^2_a {{H}^2})$, and $ \sup\nolimits_{t\in  {I_a}} \|\nabla \eta\|_2\leqslant 2\delta_2$;
  \item $\limsup_{t\to T^{\max} }\|\nabla \eta (t)\|_2 >\delta_2$,\ \ or\ \ $\limsup_{t\to T^{\max} }\|(\nabla \eta,u)( t)\|_2=\infty$\ \  if\ \  $T^{\max}<\infty$.
\end{itemize}

Let
\begin{align}  \nonumber
T^{*}:=&\sup \{\tau\in I_{T^{\max}}~ |~  \|(\nabla \eta,u)(t)\|_2^2 +\mathfrak{A}
\leqslant 2 \tilde{c}\delta^2 \ \mbox{ for any }
\ t\leqslant\tau \},\nonumber
\end{align}
where  $\mathfrak{A}$ is defined in \eqref{202220201915}.
Then, it is easy to see that the definition of $T^*$ makes sense. Thus, to show the existence of a global solution, it suffices to verify $T^*=\infty$. We shall prove this by contradiction below.

Assume $T^*<\infty$, then by Proposition \ref{202102182115}, we have
\begin{align}T^*\in (0,T^{\max})
\label{20222022519850}
 \end{align}
 and
$$ (\|(\nabla \eta,u)(t)\|_2^2 + \mathfrak{A} )|_{t=T^*} = 2 \tilde{c} \delta^2 . $$

 Noting that
\begin{equation}
\label{201911262sadf202}
\sup\nolimits_{t\in\overline{I_{T^{*}}}}(\|(\nabla \eta,u)(t) \|_2^2  + \mathfrak{A})   \leqslant  2 \tilde{c} \delta^2 \leqslant  {\delta^2_1},
\end{equation}
 we make use of \eqref{201911262sadf202} and a standard regularization method, and follow the same arguments
as used in the derivation of \eqref{estemsfsadfaalasn0} to deduce that
\begin{align*}
 {{ \sup\nolimits_{t\in \overline{I_{T^{*}}}}}(\mathcal{E}(t)+\mathfrak{E} (t))}+ \int_0^{T^{*}} ( \mathcal{D}(\tau) + \mathfrak{D}(\tau))\mm{d}\tau \leqslant \tilde{c}
 \|(\nabla \eta^0,u^0)\|_2^2 \leqslant   {\tilde{c} \delta^2}.
\end{align*}
In particular,
\begin{align}
{ \sup\nolimits_{t\in \overline{I_{T^{*}}}}(\|(\nabla \eta,u)(t) \|_2^2  + \mathfrak{A})  \leqslant\tilde{c} \delta^2 .}  \label{2020103261534}
\end{align}

By \eqref{20222022519850}, \eqref{2020103261534} and the strong continuity $(\nabla\eta,u)\in C^0([0,T^{\max}), H^2)$, we see that
  there is a  {constant} $\tilde{T}\in (T^*,T^{\max})$, such that
\begin{align}
 \sup\nolimits_{t\in \overline{I_{\tilde{T} }}} (\|(\nabla \eta,u)(t) \|_2^2+ \mathfrak{A}) \leqslant 2\tilde{c}\delta^2  , \nonumber
\end{align}
which contradicts with the definition of $T^*$. Hence, $T^*=\infty$ and thus $T^{\max}=\infty$.
This completes the proof of the existence of a global solution. The uniqueness of the global solution is obvious due to
the uniqueness result of the local solutions in Proposition \ref{202102182115}
and the fact $\sup_{t\geqslant 0}\|\nabla \eta\|_2\leqslant 2 \delta_2$.

Recalling the derivation of  \eqref{estemsfsadfaalasn0}, we easily find that
 the global solution $(\eta,u)$ enjoys \eqref{1.200}  by a standard regularization method.
To complete the proof of Theorem \ref{thm2}, {it remains} to derive the quicker decay-in-time estimate \eqref{1.safa200} for the perturbation velocity.

Noting $(\eta,u)$ satisfies \eqref{aprpiosesnew}, thus we can exploit \eqref{01dsaf16asdfasf}$_1$ and Young's inequality to   derive from \eqref{202008241448} with $i=0$ that
\begin{align}
& \frac{\mm{d}}{\mm{d}t} \left(\|\sqrt{ \bar{\rho} }  u\|^2_{ 0}
   -\kappa\int|\bar{\rho}'|^2\partial_2\eta_1\partial_1\eta_2\partial_2\eta_2\mm{d}y\right) + c   \|   u \|_1^2   \lesssim
  \|  \eta_2 \|_1^2. \nonumber
\end{align}
 Obviously we further deduce from the above inequality that, for some positive constant $\nu$ (may depend on  $\bar{\rho}$, $\mu$, $g$, $\kappa$ and $\Omega$),
\begin{align}
&\frac{\mm{d}}{\mm{d}t}\left(e^{\nu t}\left( \|\sqrt{ \bar{\rho}}u\|^2_{0}
  -\kappa\int|\bar{\rho}'|^2\partial_2\eta_1\partial_1\eta_2\partial_2\eta_2\mm{d}y\right)\right)+c e^{\nu t} \|   u\|_{ 1}^2 \lesssim  e^{\nu t} \|\eta_2\|_1^2 ,\nonumber
\end{align}
which yields that
  \begin{align}
  &\|  u\|_0^2
 +  \int_0^t e^{\nu (\tau-t) } \|   u \|_1^2 \mm{d}\tau  \lesssim
 \|  u^0\|_0^2  e^{-\nu t } +\int_0^t e^{\nu (\tau-t) }   \|\eta_2\|_1^2 \mm{d}\tau.
\label{qwessebdafgsfassgsd}
 \end{align}  Thanks to the inequalities
$$ e^{-\nu t } \lesssim\langle t \rangle^{-3}\mbox{ and } \int_0^t e^{\nu (\tau-t) } \langle \tau \rangle^{-3}\mm{d}\tau\lesssim  \langle t \rangle^{-3},$$
we immediately deduce \eqref{1.safa200} from \eqref{1.200} and \eqref{qwessebdafgsfassgsd}. Finally we easily obtain \eqref{1.200xx}
 from \eqref{1.200} and  \eqref{1.safa200} by employing an asymptotic analysis method, see \cite{JFJSZYYO} for the derivation. This completes the proof of Theorem \ref{thm2}.

\section{Proof of Theorem \ref{thm1}}\label{sec:instable}
The existence of the RT instability solutions have been widely investigated, see \cite{JFJSO2014,JFJSZWC,JFJSZYYO} for examples. We can  {also} establish the instability result for the TCRT problem of {\eqref{01dsaf16asdfasf} and \eqref{20safd45}}
in Theorem \ref{thm1} by following the same proof  {framework
as that in} \cite{JFJSZYYO}. Next, we sketch the proof for the sake of completeness.
In what follows, the fixed positive constant $c_i^I$ for $i\geqslant 1$ may depend on $\bar{\rho}$, $\mu$, $g$, $\kappa$
and the domain $\Omega$.

To begin with, we introduce the instability result for the linearized  CRT problem under the instability condition $\kappa\in[0,\kappa_{\mm{C}})$.
\begin{pro}\label{pro:08252100}
Let $\mu >0$ and $\bar{\rho}$ satisfy \eqref{0102} and \eqref{0102n}. If $\kappa\in[0,\kappa_{\mm{C}})$, then the {zero} solution
 is unstable to the linearized  capillary RT problem
\begin{equation}\label{01dsaf16asdfasf0101}
                              \begin{cases}
\eta_t=u,   \\[1mm]
\bar{\rho}u_t+\nabla \tilde{q}-\mu  \Delta u=(\kappa \mm{div}(|\bar{\rho}' |^2 \nabla\eta_2)+g\bar{\rho}'\eta_2)\mathbf{e}^2   ,\\[1mm]
\div u=0 , \\[1mm]
(\eta_2,\partial_2\eta_1,u_2,\partial_2u_1)|_{\partial\Omega}=0.
\end{cases}
\end{equation} That is, there is an unstable solution
$(\eta, u,  {\tilde{q}}):=e^{ \Lambda   t}(w/ \Lambda ,w, \beta )$
 to the above problem \eqref{01dsaf16asdfasf0101}, where
 \begin{equation}\nonumber
 (w, \beta )\in {^0\mathcal{H}^5_{\mm{s}}}\times\underline{H}^5
 \end{equation}
 solves  the boundary-value problem  \begin{equation*}
                              \begin{cases}
\Lambda ^2\bar{\rho}w+ \Lambda  \nabla  \beta
- \Lambda\mu  \Delta w=(\kappa \mm{div}(|\bar{\rho}' |^2 \nabla w_2)+g\bar{\rho}'w_2)\mathbf{e}^2
 ,\\[1mm]
\div w=0  , \\[1mm]
(w_2,\ \partial_2 w_1)|_{\partial\Omega}=0
\end{cases}
\end{equation*}
  with some growth rate $ \Lambda>0 $,  satisfying
\begin{equation}
E(v)\leqslant  {\Lambda^2}\|\sqrt{\bar{\rho}}v\|_0^2+ {\Lambda }\mu \|\nabla v\|_0^2 \ \   \mbox{ for any }\  \ v\in H_{\sigma}^{1}. \label{Lambdard} \end{equation}
In addition,
\begin{align}  \label{201602081445MH}
 \int \bar{\rho}'|w_2|^2\mm{d}y\|  w_i\|_{0}\|\partial_1w_i\|_{0}\|\partial_2w_i\|_0\neq 0\;\; \mbox{ for }\ \; i=1,\ 2. \end{align}
\end{pro}
\begin{pf} We can use Lemma \ref{20222011161567} and the modified variational {method in} \cite[Proposition 3.1]{JFJSZWC}
 to show Proposition \ref{pro:08252100}, and hence, we omit the proof here.
\hfill$\Box$
\end{pf}

Based on the linear RT instability, we turn to investigating the nonlinear RT instability. To this end, we first establish the following Gronwall-type energy inequality for the solutions of the TCRT problem.
\begin{pro}  \label{pro:0301n0845}
Let $\Lambda>0$ be the same as in Proposition \ref{01dsaf16asdfasf0101} and $(\eta,u,q)$ be the local solution constructed
by Proposition \ref{202102182115} with initial data $(\eta^0,u^0)\in {^0{\mathcal{H}}^{3,\mm{s}}_{1,\gamma}} \times {^0\mathcal{H}^2_{\mm{s}}}$.
 Assume that there exists a constant $\delta_1^I$ (may depend on $\bar{\rho}$, $\mu$, $g$, $\kappa$ and $\Omega$) {such that} $\|(\nabla \eta,u)\|_2\leqslant\delta_1^I$
in some time interval $ I_{\tilde{T}}\subset I_T$ where $I_T$ is the existence time interval of $(\eta,u,q)$. Then,  for some constant $ {c}_1^I>0$,
$(\eta,u,q)$ satisfies the Gronwall-type energy inequality:
 for all $t\in I_{\tilde{T}}$,
\begin{align}
&\Xi (t) +\underline{\|\eta\|}_3^2+ \int_0^t(\|  u\|_{\underline{2},1}^2  +\|u_t\|_1^2)(\tau)\mm{d}\tau\nonumber \\
&\leqslant {c}_1^I\|(\nabla \eta^0, u^0)\|_2^2 +\int_0^t\left({c}_1^I\|(\eta_2,u_2) \|_{0}^2+\Lambda \underline{\|\eta\|}_3^2 \right)\mm{d}\tau \label{2019120521430850} \end{align}
with
\begin{align}
 \mathcal{E} \lesssim \Xi (t)+ \underline{\|\eta\|}_3^2
 \  \ \mbox{ and }   \ \ \|\eta\|_3\lesssim \underline{\|\eta\|}_3 \lesssim  \|\eta\|_3,\label{20222121161226}
\end{align}
 Moreover, we have
\begin{align}\mathcal{D} \lesssim \|   \eta_2\|_{3}^2+\|u_t \|_{1}^2.
\label{2022211172001}
\end{align}
\end{pro}
\begin{pf}
Let  $(\eta,u,q)$ be the local solution given by Proposition  \ref{202102182115}. Then,
$(\eta,u)\in \widetilde{\mathfrak{H}}^{1,3}_{\gamma,T}\times {^0\mathcal{U}_T}$. We further assume
\begin{align*}
\|(\nabla \eta,u)\|^2_2\leqslant \delta\in (0,1]\ \  \mbox{ for all }\  \ t\in I_{\tilde{T}}\subset I_T.
\end{align*}

Recalling the derivation of \eqref{202008250856n0} and using the regularity of $(\eta,u,q)$, we conclude that
 for sufficiently small $\delta$,
\begin{align}
&\frac{\mm{d}}{\mm{d}t}\Xi_{\mm{tan}}+c  \mathfrak{D}_{\mm{tan}} \lesssim  \|(\eta_2,u_2)\|^2_{\underline{2},0} +\chi\sqrt{\mathcal{E}}\mathcal{D}    , \label{201sdf0940}
\end{align}
where $ {\Xi_{\mm{tan}}}$ is defined  by $\tilde{\mathcal{E}}$ with $ {g=0}$ and  $\chi$ is a constant such that
\begin{align}
&\| \eta \|_{\underline{2},1}^2 +\|u\|_2^2+ \|u_t \|_0^2+\|\tilde{q} \|_1^2+ \chi   \| u \|^2_{\underline{2},0}   \lesssim  \Xi_{\mm{tan}} \lesssim \chi \mathcal{E}
\label{202sdaf21057}
\end{align}
and
\begin{align}
   \chi\|  u\|_{\underline{2},1}^2  +\|u_t\|_1^2
   \lesssim  {\mathcal{D}}_{\mm{tan}} .  \label{2022saf202012047}
\end{align}

Noting that $(\eta,u)$ also satisfies \eqref{omessen} and \eqref{omessetsim}, we derive form the relation $\eta_t=u$ that
 \begin{align}
\frac{\mm{d}}{\mm{d}t}\|\eta\|_{\underline{1},2}^2
 &  \lesssim \|\eta\|_{\underline{1},2}\| u\|_{\underline{1},2}
   \lesssim \|\eta\|_{\underline{1},2}(\|\eta\|_{\underline{2},1} +\| u_t\|_{1}),  \label{2022sdfasx02012047}\\
  \frac{\mm{d}}{\mm{d}t}\|\eta\|_{3}^2
  & \lesssim \|\eta\|_{3}\| u\|_{3}\lesssim\|\eta\|_{3}(\|\eta\|_{1,2}+\| u_t\|_1).   \label{22012047}
\end{align}
We immediately derive from \eqref{201sdf0940}, \eqref{2022sdfasx02012047} and \eqref{22012047} that\begin{align}
&\frac{\mm{d}}{\mm{d}t}(\Xi_{\mm{tan}}+ a_1\|\eta\|_{\underline{1},2}^2 + a_1^2\|\eta\|_{3}^2 ) +c  \mathfrak{D}_{\mm{tan}}\nonumber \\
 &\lesssim  \|(\eta_2,u_2)\|^2_{\underline{2},0} +\chi\sqrt{\mathcal{E}}\mathcal{D}  + a_1 \|\eta\|_{\underline{1},2}(\|\eta\|_{\underline{2},1} +\| u_t\|_{1}) + a_1^2\|\eta\|_{3}(\|\eta\|_{1,2}+\| u_t\|_1). \label{201sdf0asdf940}
\end{align}

 {By integration by parts and Cauchy--Schwarz's inequality},
we have that for any $\varepsilon\in (0,1]$,
\begin{align}
 \| \chi_2\|_{k,0}  \lesssim
 \begin{cases}
 \varepsilon^{-1}  \|\chi_2\|_0 + \varepsilon \| \chi_2\|_{2,0}   & \hbox{for }\  \ k=1; \\
 \varepsilon^{-1}\| \chi_2\|_{1,0} + \varepsilon \| \chi _2\|_{3,0 } & \hbox{for } \ \  k=2,
 \end{cases}
  \label{202102181813}
\end{align}
where $\chi=\eta$ or $u$. In addition, $(\eta,u)$ also satisfies the estimate \eqref{omessetsim}, thus
\begin{align}
 \mathcal{D}   \lesssim  \|   \eta_2\|_{3}^2+\|u_t \|_{1}^2 . \label{2010940}
\end{align}
Similarly to \eqref{202012252005}, we also have
\begin{align}
\sqrt{\mathcal{E}}\mbox{ and } \| (\nabla \eta,u)\|_2 \mbox{ are equivalent to each other} , \label{201asf0940}
\end{align} Therefore, with the help of \eqref{202102181813}--\eqref{201asf0940}, we deduce \eqref{2019120521430850}
from \eqref{201sdf0asdf940} for sufficiently small $a_1$ and $\delta$, where $\Xi$, $\underline{\|\eta\|}$ and  $\mathcal{D}$ satisfy \eqref{20222121161226} and \eqref{2022211172001}. The proof is complete.
\hfill$\Box$
\end{pf}

For any given $\delta>0$, let
\begin{equation}\label{0501}
\left(\eta^\mm{a}, u^\mm{a}, q^\mm{a}\right)=\delta e^{ \Lambda  t } (\tilde{\eta}^0, \tilde{u}^0, \tilde{q}^0),
\end{equation}
where $ (\tilde{\eta}^0, \tilde{u}^0, \tilde{q}^0):=(w/\Lambda ,w,\beta  )$, and $(w,\beta ,\Lambda)$ is given by Proposition \ref{pro:08252100}.
Then, $\left(\eta^\mm{a},u^\mm{a},q^\mm{a}\right)$ is also a solution to the linearized CRT problem \eqref{01dsaf16asdfasf0101},
and enjoys the estimate: for any $j\geqslant 0$,
\begin{equation}
\label{appesimtsofu1857}
\|\partial_{t}^j(\eta^\mm{a}, u^\mm{a})\|_3+\|\partial_{t}^j q^\mm{a} \|_2=\Lambda^j\delta e^{\Lambda t}(\|(\tilde{\eta}^0,\tilde{u}^0)\|_3
+\|\tilde{q}^0\|_2)\lesssim \Lambda^j \delta e^{\Lambda t}.
\end{equation}
In addition, we have by \eqref{201602081445MH} that
\begin{eqnarray}\label{n05022052}
\|\bar{\rho}'\chi_2\|_0\|\chi _i\|_{0}\|\partial_1\chi _i\|_{0}\|\partial_2\chi _i\|_{0} >0,\quad i=1,2,
\end{eqnarray}
where $\chi  =\tilde{\eta}^0 $ or $\tilde{u}^0 $.

 Since the initial data of the solution $( {\eta}^{\mm{a}},{u}^{\mm{a}},{q}^{\mm{a}})$ to the linearized CRT problem
may not satisfy the necessary compatibility conditions required by the initial data of the corresponding (nonlinear) TCRT problem.
So, we have to modify the initial data of the linearized problem as done in \cite[Proposition 5.1]{JFJSZWC},
such that the modified initial data approximate the original initial data of the linearized problem,
and satisfy the compatibility conditions for the corresponding nonlinear problem. More precisely,
\begin{pro}\label{pro:0101}
Let $(\tilde{\eta}^0, \tilde{u}^0 ):=(w/\Lambda , w)$  be the same as in \eqref{0501}. Then, there
is a constant $\delta_2^I \in (0,1]$, such that for any $\delta\in(0, \delta_2^I ]$, there exists
$(\eta^{\mm{r}}, u^{\mm{r}})\in {^0\mathcal{H}^3_{\mm{s}}}$ enjoying the following properties:
\begin{enumerate}
\item[(1)] The modified initial data
\begin{align}\nonumber
({\eta}^{\delta}_{0},{u}^{\delta}_{0 } )
:=\delta(\tilde{\eta}^0,\tilde{u}^0 )+\delta^2(\eta^{\mm{r}},u^{\mm{r}} )
\end{align}
belongs to
${^0{\mathcal{H}}^{3,\mm{s}}_{1,\gamma}} \times {^0\mathcal{H}^3_{\mm{s}}}$ and satisfies
the compatibility condition
\begin{align}\nonumber
&\mm{div}_{\mathcal{A}_{0 }^{\delta}}u_{0 }^{\delta}=0  \mbox{ in } \Omega,
\end{align}  where  $\mathcal{A}^{\delta}_0$ is defined as $\mathcal{A}$ with $\eta^{\delta}_{0 }$ in place of $\eta$.
\item[(2)]
Uniform-in-$\delta$ estimate:
\begin{align}
\label{202103281107}
 \| (\eta^{\mm{r}},u^{\mm{r}}  )\|_3 \leqslant {c}_2^I,
\end{align}
where the positive constant $ {c}_2^I $ is independent of $\delta$.
\end{enumerate}
\end{pro}
\begin{pf}
Thanks to Lemma \ref{20222011161567}, we can obtain Proposition \ref{pro:0101}
by employing the same arguments used {in} \cite[Proposition 5.1]{JFJSZWC}, and hence, we omit the detailed proof here.
\hfill$\Box$
\end{pf}

Now,  let $({\eta}^{\delta}_{0},{u}^{\delta}_0 )\in {^0{\mathcal{H}}^{3,\mm{s}}_{1,\gamma}} \times {^0\mathcal{H}^3_{\mm{s}}} $ be constructed in  Proposition \ref{pro:0101},
\begin{align}
&\label{201912041727}
 {c}_3^I= {\|(\tilde{\eta}^0,\tilde{u}^0)\|_3}  + {c}_2^I>0
 \end{align}and
 \begin{align}
&\label{201912041705}
 {\delta}_0= \frac{1}{2c_3^I } \min\left\{\gamma,{\delta_1^I}, {{\delta}_2} , 2c_3^I {\delta}_2^I \right\}\leqslant 1 .
\end{align}
From now on, we assume that $\delta\leqslant {\delta}_0$. Since $\delta\leqslant\delta_2^I$, we can use Proposition \ref{pro:0101} to construct
$(\eta_0^\delta, u_0^\delta)$ that satisfies
$$  {\|(\eta_0^\delta,u_0^\delta)\|_3}
\leqslant c_3^I\delta\leqslant  {\delta}_2 .$$
By virtue of Proposition \ref{202102182115}, there exists a unique solution $(\eta, u, q)$ of the TCRT problem  of
\eqref{01dsaf16asdfasf} and \eqref{20safd45} with the initial data $({\eta}_0^\delta, {u}_0^\delta)$ in place of $(\eta^0, u^0)$,
where $(\eta, u, q)\in \widetilde{\mathfrak{H}}^{1,3}_{\gamma,\tau}\times  {^0\mathcal{U}_{\tau}}
\times ( C^0(\overline{I_\tau}, \underline{H}^1)\cap L^2_{\tau} {H}^2)$ for $\tau\in I_{T^{\mm{\max}}}$,
 and $T^{\mm{max}}$ denotes the maximal time of existence.

Let $\epsilon_0\in (0,1]$ be a constant  that will be determined in \eqref{201907111842} later. We define
\begin{align}\label{times}
T^\delta:=\,&\Lambda^{-1}\mm{ln}({\epsilon_0}/{\delta})>0,\;\; \mbox{ i.e., }\;
 \delta e^{ \Lambda  T^\delta }=\epsilon_0,\\[1mm]
T^*:=\,&\sup\left\{t\in I_{T^{\max}}\left|~
\sup\nolimits_{\tau\in [0,t)}\sqrt{\|\eta(\tau)\|_3^2+\| u (\tau)\|_2^2}\right.\leqslant 2   {c}_3^I \delta_0 \right\}, \label{xfdssdafatimes}\\[1mm]
T^{**}:=\,&  \sup\left\{t\in I_{T^{\max}} \left|~\sup\nolimits_{\tau\in [0,t)}\left\|(\eta,u)(\tau)\right\|_{0}\leqslant 2 {c}_3^I\delta  e^{\Lambda \tau}
 \right\}.\right.  \label{xfdstimes}
\end{align}

Since
\begin{align}
 \left.\sqrt{\|\eta(t)\|_3^2+\| u(t)\|_2^2} \right|_{t=0} = \sqrt{\|\eta_0^\delta\|_3^2+\|u_0^\delta\|_2^2}
\leqslant c_3^I \delta  <  2c_3^I \delta ,
\label{201809121553}
\end{align}
we have $T^{*}>0$, $T^{**}>0$. Obviously,
\begin{align}
\label{0n111} & T^{*}=T^{\max}=\infty\ \  \mbox{ or }\  \ T^{*}<T^{\max}, \\
\label{0502n111} & \left\|(\eta,u) (T^{**}) \right\|_0 =2{c}_3^I\delta e^{\Lambda T^{**}} \ \ \mbox{ if } \ \ T^{**}<T^{\max}, \\
\label{050211} & \sqrt{\|\eta (T^{*})\|_3^2+\| u  (T^{*})\|_2^2} =2   {c}_3^I \delta_0\  \ \mbox{ if }\  \ T^{*}<T^{\max}.
\end{align}

From now on, we define
$$T^{\min}:=\min\{T^\delta,T^*,T^{**}\} .$$
Noting that $\sup\nolimits_{t\in [0,T^{\min} )}\sqrt{\|\eta(t)\|_3^2+\| u (t)\|_2^2} \leqslant \delta_1^I$, we use
Proposition \ref{pro:0301n0845} to see that $(\eta,u,q)$ enjoys the Gronwall-type energy inequality \eqref{2019120521430850} and \eqref{20222121161226}
for any $t\in I_{T^{\min}}$. If we utilize this fact, \eqref{xfdssdafatimes}--\eqref{201809121553}, Lemma \ref{pro:1221} and
the condition $\|\eta\|_3\leqslant \gamma$, we find that for any $t\in [0 , T^{\min})$,
\begin{align}
 & \mathcal{E} (t) +  \int_0^t \mathcal{D}(\tau) \mm{d}\tau  \leqslant c_4^I\delta^2e^{2\Lambda  t},
\label{20191204safda2114} \\
&  \left\|\int_{0}^{\eta_2 }\left(\eta_2 - z\right)\bar{\rho} {''}(y_2+z)\mm{d}z \right\|_{{L^1}}
\leqslant (c_4^I \delta  e^{ \Lambda t})^2 . \label{2022202101204}
\end{align}

Thanks to \eqref{20191204safda2114}, we proceed to bound the errors between $(\eta,u)$ and $(\eta^{\mm{a}}, u^{\mm{a}})$.
\begin{pro}\label{2022202101315}
Let
$(\eta^{\mathrm{d}}, u^{\mathrm{d}} )=(\eta, u )-(\eta^{\mm{a}}, u^{\mm{a}} )$, then there is a constant
${c}_5^I$, such that for any $\delta\in(0,\delta_0]$ and for any $t\in I_{T^{\min}}$,
\begin{align}
\label{ereroe}
\|( {\bar{\rho}'\chi^{\mm{d}}},\eta^{\mm{d}},u^{\mm{d}} ) \|_{\mathfrak{X}} +  \|u_t^{\mm{d}}\|_{0}   &\leqslant {c}_5^I\sqrt{\delta^3e^{3\Lambda t}}, \\
\|(\mathcal{A}_{ik}\partial_k\chi_j
-\partial_i \chi^{\mm{a}}_j)(t)
\|_{{L^1}}&\leqslant  {c}_5^I\sqrt{\delta^3 e^{3\Lambda t}},\label{2018090119021}
\end{align}
where $i,j=1,2$, $\chi=\eta$ or $u$, $ \mathfrak{X}=W^{1,1}$ or $H^1$, and ${c}_5^I$ is independent  of  ${T^{\min}}$.
\end{pro}
\begin{pf}
We can easily obtain Proposition \ref{2022202101315}
by employing the same arguments used in \cite[Proposition 6.1]{JFJSZWC} or \cite[Proposition 3.4]{JFJSZYYO}, and hence, we omit the detailed proof here. \hfill$\Box$
\end{pf}

Let
 $$  \begin{aligned}
\varpi :=& \min_{ {\chi} =\tilde{\eta}^0, \tilde{u}^0} \left\{\|{\chi}_1\|_{L^1}, \|\partial_1{\chi}_1\|_{L^1}, \|\partial_2{\chi}_1\|_{L^1},\|\chi_2\|_{L^1} ,\|\partial_1\chi_2\|_{L^1},\|\partial_2\chi_2\|_{L^1},\|\bar{\rho}'\tilde{\eta}_2^0\|_{L^1} \right\},
\end{aligned}$$
then $\varpi >0$ by \eqref{n05022052}.

Now, denoting
\begin{equation}\label{201907111842}
\epsilon_{0}=\min \left\{ \left(\frac{ {c}_3^I}{2 c_{5}^I }\right)^2, \frac{c_3^I \delta_0}{c_4^I} ,
\frac{\varpi ^2}{4 ({c}_{5}^I+ |{c}_{4}^I|^{2})^2},1 \right\}>0 ,
\end{equation}
we claim that
\begin{equation}\label{201907111840}
T^{\delta}=T^{\mm{min}} = \min \left\{T^{\delta}, T^{*}, T^{* *}\right\} \neq T^{*}\;\;\mbox{or}\;\;T^{**},
\end{equation}
which, in fact, can be shown by contradiction as follows.
\begin{enumerate}[(1)]
  \item If $T^{\mm{min}}=T^{**} $, then $T^{**} < T^{\max} $ by \eqref{0n111}.
Recalling $\sqrt{\epsilon_0}\leqslant {c}_3^I/2{c}_5^I$, we see by
\eqref{0501},  \eqref{201912041727}, \eqref{times} and \eqref{ereroe} that
\begin{equation*}\begin{aligned}
\|  \eta  (T^{**})\|_{0}&\leqslant \| \eta^\mm{a} (T^{**})\|_{0}+\| \eta^\mm{d} (T^{**})\|_{0}\\
&\leqslant  \delta e^{{\Lambda T^{**}}}(c_3^I+ c_5^I\sqrt{\delta e^{\Lambda  T^{**}}})\\
&\leqslant  \delta e^{{\Lambda  T^{**}}}(c_3^I+ c_5^I\sqrt{\epsilon_0})\\
&\leqslant 3c_3^I \delta e^{\Lambda  T^{**}}/2< 2c_3^I\delta e^{\Lambda T^{**}},
\end{aligned} \end{equation*}
 which is in contradiction with \eqref{0502n111}.  Hence, $T^{\mm{min}}\neq T^{**}$.
  \item If $T^{\mm{min}}=T^{*}$, then $T^{*}<T^{**}  $. Recalling
$  \epsilon_0\leqslant {c}_3^I  \delta_0/c_4^{I}$, we deduce from \eqref{20191204safda2114}
that for any $t\in I_{T^{\min}}$,
\begin{equation}\nonumber
 \sqrt{\|\eta(t)\|_3^2+\| u(t)\|_2^2  }  \leqslant c_4^I \delta  e^{ \Lambda T^\delta} \leqslant  {c}_3^I  \delta_0<2 {c}_3^I \delta_0,
\end{equation}
which contradicts  with \eqref{050211}. Hence, $T^{\mm{min}}\neq T^{*}$.
\end{enumerate}

Since $T^{\delta} $ satisfies \eqref{201907111840}, the inequalities \eqref{ereroe} {and} \eqref{2018090119021} hold to $t=T^\delta$.
Using this fact, and \eqref{0501}, \eqref{times}, \eqref{2022202101204} and the condition
$\epsilon_0\leqslant {\varpi ^2}/ 4 ({c}_{5}^I+|{c}_4^I |^{2})^2$, we obtain the following instability relations:
 \begin{align}
\|\mathcal{A}_{ik}\partial_k\chi_j(T^\delta) \|_{{L^1}}\geqslant \,& \| \partial_i \chi^{\mm{a}}_j(T^\delta) \|_{{L^1}}
- \|\mathcal{A}_{ik}\partial_k\chi_j(T^\delta) -\partial_i\chi^{\mm{a}}_j(T^\delta) \|_{{L^1}}\nonumber  \\
\geqslant \, & \delta e^{\Lambda T^\delta }( \|\partial_i\tilde{\chi}^{0}_j\|_{L^1}- {c}_{5}^I\sqrt{\delta e^{\Lambda T^\delta }})
\geqslant \epsilon:= \varpi \epsilon_0 /2,\quad\; i,j=1,2, \nonumber
\end{align}
and
\begin{align}
\|\bar{\rho}-\bar{\rho}(\eta_2(y,T^\delta) +y_2) \|_{{L^1}} & \geqslant
\| \bar{\rho}' \eta^{\mm{a}}_2(T^\delta) \|_{{L^1}}-\| \bar{\rho}' \eta^{\mm{d}}_2(T^\delta) \|_{{L^1}} \nonumber \\
& \quad - \left\|\int_{0}^{\eta_2(y,T^\delta)}\left(\eta_2(y, {T^\delta}) - z\right)\bar{\rho}''(y_2+z)\mm{d}z \right\|_{{L^1}} \nonumber\\
& \geqslant \delta e^{\Lambda T^\delta }( \|\bar{\rho}'\tilde{\eta}^{0}_2\|_{L^1}- {c}_{5}^I\sqrt{\delta e^{\Lambda T^\delta }}-(c_4^I)^2
\delta  e^{ \Lambda T^\delta} ) \geqslant  \varpi \epsilon_0 /2 ,\nonumber
\end{align}
where $\chi=\eta$ or $u$. Similarly, we can also verify that $(\eta,u)$ satisfies the rest instability relations in \eqref{201806012326} by using \eqref{ereroe}.
This completes the proof of  Theorem \ref{thm1}.

\appendix
\section{Analysis tools}\label{sec:09}
\renewcommand\thesection{A}
This appendix is devoted to providing some mathematical results, which have been used in previous sections.
We should point out that  {$\Omega$ and the simplified notations appearing in what follows are  as same as these defined by \eqref{0101a} and  Section \ref{subsec:04}}. In addition, $a\lesssim b$ still denotes $a\leqslant cb$
where the positive constant $c$ depends on the parameters and the domain in the lemmas in which $c$ appears.
\begin{lem}\label{201806171834}
\begin{enumerate}[(1)]
 \item Interpolation inequality in $H^j$
 (see \cite[{Section 5.2}]{ARAJJFF}): Let $D$ be a domain in $\mathbb{R}^2$ satisfying the cone condition or a bounded interval in $\mathbb{R}$, then for any given $0\leqslant j<i$,
\begin{equation}\label{201807291850}
\|f\|_{H^j(D)}\lesssim\|f\|_{L^2(D)}^{(i-j)/i}\|f\|_{H^i(D)}^{j/i} \lesssim\varepsilon^{-j/(i-j)}\|f\|_{L^2(D)} +\varepsilon \|f\|_{H^i(D)},
\quad\forall\;\varepsilon >0,
\end{equation}
where the two estimate constants in \eqref{201807291850} are independent of $\varepsilon$.
\item
Product estimates (see   \cite[Section 4.1]{JFJSNS}): Let $D\subset \mathbb{R}^2$ be a domain satisfying the cone condition,
 and $\varphi$, $\psi$ be functions defined in $D$.  Then
\begin{align}
\label{fgestims}
&  \|\varphi\psi\|_{H^i(D)}\lesssim
  \|\varphi\|_{H^i(D)}\|\psi\|_{H^2(D)}\quad \hbox{ for }\,\,0\leqslant i\leqslant 2.
\end{align}
\item  Anisotropic product estimates
({see} \cite[Lemma 3.1]{jiang2021asymptotic}): Let  the functions $\varphi$ and $\psi$ be defined in $\Omega$. Then
\begin{align}
\label{fgessfdims}
  \|\varphi\psi\|_0\lesssim  \begin{cases}\sqrt{\|\varphi \|_{0} \|\varphi \|_{\underline{1},0}} \|\psi \|_{1},\\
 \|\varphi \|_{0} \sqrt{\|\psi\|_{\underline{1},0} \|\psi\|_{\underline{1},1}} . \end{cases}
\end{align}                \end{enumerate}
\end{lem}

\begin{lem}\label{10220830}
Poincar\'e's inequality (see \cite[Lemma 1.43]{NASII04}): Let $1\leqslant p<\infty$, and $D$ be a bounded Lipchitz domain in $\mathbb{R}^n$ for $n\geqslant 2$ or a finite interval in $\mathbb{R}$. Then for any $w\in W^{1,p}(D)$,
\begin{equation}
\label{poincare}
\|w\|_{L^p(D)}\lesssim \|\nabla w\|_{L^p(D)}^p+\left|\int_{D}w\mathrm{d}y\right|^p.
\end{equation}
\end{lem}
\begin{rem}\label{10220saf830p}
By Poincar\'e's inequality, we have, for any given $i\geqslant  0$,
\begin{align}
&\label{202012241002}
\| w\|_{1,i}\lesssim \|w\|_{2,i}\; \mbox{ for any } w \mbox{ satisfying }\partial_1w,\ \partial_1^2w\in H^i.
\end{align}
 \end{rem}
\begin{lem}\label{pro4a}
Hodge-type elliptic estimates (see \cite[Lemma A.4]{ZHAOYUI}):
If $w\in H^i_{\mm{s}}$ with $i\geqslant1$,
then
\begin{align}
&\label{202005021302}
\|\nabla w\|_{i-1}
\lesssim\|(\mm{curl}w,\mm{div}w)\|_{i-1}.
\end{align}
\end{lem}
\begin{lem}\label{pro:1221}Diffeomorphism mapping theorem (see \cite[Lemma A.8]{ZHAOYUI}):
There is a sufficiently small constant $\gamma \in(0,1)$, depending possibly on $\Omega$, such that for any ${\varsigma}\in H_{\mm{s}}^3$
satisfying $\|\nabla {\varsigma}\|_2\leqslant \gamma $,
$\psi:=\varsigma+y$ (after possibly being redefined on a set of measure zero with respect to variable $y$) satisfies
the same diffeomorphism  properties as $\zeta$ in \eqref{20210301715x} and \eqref{20210301715}, and
$\inf_{y\in \Omega}\det(\nabla \varsigma  +I)\geqslant 1/4$.
\end{lem}

\begin{lem}\label{20222011161567}  Stokes estimates  (see \cite[Lemma A.12]{JFJSZYYO}):
Let $i\geqslant 0$, $(f,\varphi)\in H^i\times H^{2+i}$ and $\varphi_2|_{{\partial\Omega}}=0$,
then the Stokes problem with Navier boundary condition
\begin{equation}
\label{Ston}
\begin{cases}
-  \Delta v + \nabla P = f &\mbox{in }\ \ {\Omega},\\
 \div v=\div {\varphi} &\mbox{in }\ \ {  \Omega}, \\
  (v_2,\partial_2  v_1 )=0&\mbox{on }\ \ {\partial\Omega}
\end{cases}
\end{equation}
admits a unique  solution $(v ,P)\in  {\underline{\mathcal{H}^{2+i }_{\mm{s}}}}$  (or ${^0\mathcal{H}^{2+i }_{\mm{s}}}$) $\times \underline{H}^{1+i }$,  satisfying
\begin{align}
\label{202201122130}
\|v\|_{2+i }+\|P\|_{1+i }\lesssim \|(f,{\varphi} )\|_{i}+\|\mm{div} {\varphi} \|_{1+i }.
\end{align}
\end{lem}
 \begin{lem}\label{261567}  Trace estimate:
\begin{align}
&\label{202108261406}
\int_{\partial\Omega}\partial_1\varphi\psi\mm{d}y_1\lesssim
\sqrt{
\| \varphi\|_{1,0}(\|\varphi\|_0+ \|\partial_2\varphi\|_0)}
\sqrt{
\| \psi\|_{1,0}(\|\psi\|_0+\|\partial_2\psi\|_0) }
\end{align}
for any $\varphi $, $\psi \in H^1$.
\end{lem}
\begin{pf}
We denote by $\overline{f}$  the conjugate function of $f$ and $\widehat{f}(\xi,y_2)$ the horizontal Fourier transform of $f(y) $, i.e.,
$$ \widehat{f}(\xi,y_2)=\int_{2\pi L\mathbb{T}}f( y_1,y_2)e^{-\mm{i}y_1\xi}\mm{d}y_1, $$
 then $\widehat{\partial_2 f} = \partial_{2} \widehat{f}$.  By the Parseval theorem  (see \cite[Proposition 3.1.16]{grafakos2008classical}),
\begin{align}
\int_{\Sigma_0}\partial_1\varphi\psi\mm{d}y_1=\,&\frac{1}{4(\pi {L})^2}\sum_{\xi\in L^{-1}\mathbb{Z}}
 \mm{i}\xi\widehat{\varphi }(\xi,0 ) \overline{ \widehat{\psi }(\xi,0 )} \nonumber \\ \lesssim\,&\frac{1}{4(\pi {L})^2}
\sqrt{\sum_{\xi\in L^{-1}\mathbb{Z}}
|\xi| |\widehat{\varphi }(\xi,0 )|^2 \sum_{\xi\in L^{-1}\mathbb{Z}}|\xi||  \widehat{\psi}(\xi,0 ) |^2}
, \label{20222010221241}
\end{align}
where $\Sigma_0:=2\pi L\mathbb{T}\times \{0\}$.

It is easy to estimate that, for a.e. $y_2\in (0,h)$,
\begin{align}
& \sum_{\xi\in L^{-1}\mathbb{Z}} |\xi| |\widehat{\varphi }(\xi,0 )|^2 \nonumber \\
 {=} \,&
\sum_{\xi\in L^{-1}\mathbb{Z}}
 |\xi| \left( |\widehat{\varphi }(\xi,y_2 )|^2 -2\int_0^{y_2} \widehat{\varphi }(\xi,s ) \partial_s\widehat{\varphi }(\xi,s ) \mm{d}s \right) \nonumber \\
\lesssim\,  & \left(\sum_{\xi\in L^{-1}\mathbb{Z}}
 |\widehat{\varphi }(\xi,y_2 )|^2 \sum_{\xi\in L^{-1}\mathbb{Z}} |\xi\widehat{\varphi }(\xi,y_2 )|^2 \right)^{1/2}\nonumber \\
&\qquad +
\left(
 \sum_{\xi\in L^{-1}\mathbb{Z}}   \|\xi\widehat{\varphi }(\xi,s )\|^2_{L^2(0,h)}
 \sum_{\xi\in L^{-1}\mathbb{Z}}   \| \partial_s\widehat{\varphi }(\xi,s )\|^2_{L^2(0,h)}  \right)^{1/2} \nonumber\\
 =\, &  4(\pi {L})^2( \|\varphi(y_2)\|_{L^2(0,{2\pi L})}\| \partial_1\varphi(y_2)\|_{L^2(0,{2\pi L})}+\| \varphi\|_{1,0}\|\partial_2\varphi\|_0).\nonumber  \end{align}
Similarly, we have
$$ \sum_{\xi\in L^{-1}\mathbb{Z}}|\xi|| { \widehat{\psi }(\xi,0 )}|^2\lesssim 4(\pi {L})^2( \|\psi(y_2)\|_{L^2(0,{2\pi L})}\| \partial_1\psi(y_2)\|_{L^2(0,{2\pi L})}
+\|\psi\|_{1,0}\|\partial_2\psi\|_0). $$

Consequently, putting the above two estimates into \eqref{20222010221241}, and then integrating the resulting inequality over $(0,h)$, we immediately get
$$ \int_{\Sigma_0}\partial_1\varphi\psi\mm{d}y_1\lesssim
\sqrt{
\| \varphi\|_{1,0}(\|\varphi\|_0+ \|\partial_2\varphi\|_0) }
\sqrt{
\| \psi\|_{1,0}( \|\psi\|_0+ \|\partial_2\psi\|_0 )}. $$
Similarly,
$$ \int_{\Sigma_h}\partial_1\varphi\psi\mm{d}y_1\lesssim
\sqrt{
\| \varphi\|_{1,0}( \|\varphi\|_0+ \|\partial_2\varphi\|_0 )}
\sqrt{
\| \psi\|_{1,0}( \|\psi\|_0+ \|\partial_2\psi\|_0) }, $$
where $\Sigma_h:=2\pi L\mathbb{T}\times \{h\}$.
Putting the above two estimates together yields \eqref{202108261406}. This completes the proof.\hfill $\square$
\end{pf}
\begin{lem}\label{261asdas567} Poincar\'e's inequality for vector functions with optimal constant: for any $\varphi\in H^1_\sigma$, it holds that
\begin{align}
&\label{20210asda8261406}\|\varphi_2\|_{0}^2\leqslant  \|\nabla \varphi_2\|_{0}^2/(\pi^2h^{-2}+L^{-2}) ,
\end{align}
where the constant $(\pi^2h^{-2}+L^{-2})^{-1}$ is optimal.
\end{lem}
\begin{pf}
We still denote $\widehat{f}(\xi,y_2)$ be the horizontal Fourier transform of $f(y) $. Since $\mm{div}\varphi=0$, we have
 $$ i\xi \widehat{\varphi_1}+\partial_2\widehat{\varphi_2}=0.$$
 Taking $\xi=0$ in the above identity and then using the boundary condition $\widehat{\varphi_2}(0,0)=0$, we have
 \begin{align}\widehat{\varphi_2}(0,y_2)
 =\partial_2\widehat{\varphi_2}(0,y_2)=0
 .  \label{20222101512150}
  \end{align}

 In addition, it is well-known that there exists a function $\psi_0\in H^1_0(0,h)$ such that (see Lemma 4.4 and  (4.25) in \cite{JFJSARMA2019})
\begin{align}
\frac{\|\psi_0\|_{L^2(0,h)}}{\|\psi'_0\|_{L^2(0,h)}} = \sup_{\psi\in H^1_0(0,h)}\frac{\|\psi\|_{L^2(0,h)}}{\|\psi'\|_{L^2(0,h)}}=\frac{h}{\pi}.
\label{202220101519822}
\end{align}
By the Parseval theorem  (see \cite[Proposition 3.1.16]{grafakos2008classical}),  \eqref{20222101512150} and \eqref{202220101519822}, we have
\begin{align}
\|\nabla \varphi_2\|^2_0=  \,& \frac{1}{4(\pi {L})^2}\sum_{\xi\in L^{-1}\mathbb{Z}}
( \|\xi\widehat{\varphi_2}(\xi,y_2)\|^2_{L^2(0,h)}+ \|\partial_2\widehat{\varphi_2} (\xi,y_2)\|^2_{L^2(0,h)}) \nonumber \\
 \geqslant \,& \frac{( \pi^2 h^{-2}+L^{-2})}{4(\pi {L})^2}\sum_{0\neq \xi\in L^{-1}\mathbb{Z}}
\|\widehat{\varphi_2}(\xi,y_2)\|^2_{L^2(0,h)}= ( \pi^2 h^{-2}+L^{-2}) \| \varphi_2\|_0^2, \nonumber
\end{align}
which implies \eqref{20210asda8261406}.

Now we further prove that the constant $(\pi^2h^{-2}+L^{-2})^{-1}$ is optimal. To this purpose, we define that
$$ \varphi := (-{L} \psi_0'(y_2)\cos (  y_1/L),-\psi_0(y_2)\sin (  y_1 /L) ).$$
It is easy to check that $\varphi \in H^1_\sigma$; moreover, by \eqref{202220101519822},
\begin{align}
\frac{\|\varphi _2\|_0^2}{\|\nabla \varphi _2\|_{0}^2} =  \,&\frac{\|\psi_0(y_2) \|_{L^2(0,h)}^2\| \sin (  y_1 /L)\|_{L^2(0,2\pi L)}^2}{\|\psi_0'(y_2) \|_{L^2(0,h)}^2\|\sin (  y_1/L) \|_{L^2(0,2\pi L)}^2+ \|\psi_0  (y_2) \|_{L^2(0,h)}^2\|\cos (  y_1/L) \|_{L^2(0,2\pi L)}^2 L^{-2}}\nonumber \\
= \,&( \pi^2 h^{-2}+L^{-2})^{-1} .  \label{202221131926}
\end{align}
This means that the constant $(\pi^2h^{-2}+L^{-2})^{-1}$ is optimal, and thus completes the proof.\hfill$\Box$
\end{pf}
\vspace{4mm}
\noindent\textbf{Acknowledgements.}
The research of Fei Jiang was supported by NSFC (Grant Nos.  12022102 and 12231016) and the Natural Science Foundation of Fujian Province of China (Grant Nos. 2020J02013 and 2022J01105), and the research of  Fucai Li by NSFC (Grant No. 12071212)
and a project funded by the priority academic program development of Jiangsu
higher education institutions, and the research of  Zhang by NSFC (Grant No. 12101305).

\renewcommand\refname{References}
\renewenvironment{thebibliography}[1]{%
\section*{\refname}
\list{{\arabic{enumi}}}{\def\makelabel##1{\hss{##1}}\topsep=0mm
\parsep=0mm
\partopsep=0mm\itemsep=0mm
\labelsep=1ex\itemindent=0mm
\settowidth\labelwidth{\small[#1]}
\leftmargin\labelwidth \advance\leftmargin\labelsep
\advance\leftmargin -\itemindent
\usecounter{enumi}}\small
\def\newblock{\ }
\sloppy\clubpenalty4000\widowpenalty4000
\sfcode`\.=1000\relax}{\endlist}
\bibliographystyle{model1b-num-names}

\end{document}